\documentclass{amsart} 
\usepackage{amsmath}
\usepackage{epsfig}
\usepackage{siunitx} 
\usepackage{color,graphics}
\usepackage{booktabs}
\usepackage{caption}
\usepackage{subcaption}
\usepackage{amsthm}
\usepackage{amsmath}
\usepackage{amssymb}
\usepackage{lipsum}
\usepackage{pdflscape,longtable,array}
\usepackage{afterpage}
\usepackage{subcaption}
\usepackage{pdflscape}
\usepackage{fancyhdr} 

\fancypagestyle{mylandscape}{
\fancyhf{} %Clears the header/footer
\fancyfoot{% Footer
\makebox[\textwidth][r]{% Right
  \rlap{\hspace{.75cm}% Push out of margin by \footskip
    \smash{% Remove vertical height
      \raisebox{4.87in}{% Raise vertically
        \rotatebox{90}{\thepage}}}}}}% Rotate counter-clockwise
}
\usepackage{algorithm}
\usepackage{algorithmic}
\usepackage{tabularx} 
\usepackage{lineno}
\usepackage{diagbox} % for table cell with a diagonal line and 2 sub cells 
\numberwithin{equation}{section}
\DeclareMathOperator{\sech}{sech}
\usepackage{xcolor}

\newcolumntype{^}{>{\global\let\currentrowstyle\relax}}% <==========for bold row in table
\newcolumntype{_}{>{\currentrowstyle}}% <==========for bold row in table
\newcommand{\rowstyle}[1]{\gdef\currentrowstyle{#1} #1\ignorespaces} % <==========for bold row in table
\usepackage{multirow}
\usepackage{epstopdf}
\usepackage[all]{xy}
\usepackage{fixmath, bm, amssymb, amsfonts, latexsym}

\newtheorem{theorem}{Theorem}[section]

\newtheorem{example}[theorem]{Example}
\newtheorem{remark}[theorem]{Remark}
\begin{document}
\author{Lavanya V Salian$^{\S}$ and Samala Rathan$^\dagger$}
\title{Exponential approximation space reconstruction WENO scheme for dispersive PDEs}
\thanks{
$^{\S}$Department of Humanities and Sciences, Indian Institute of Petroleum and Energy-Visakhapatnam,
India-530003, (lavanya\_vs@iipe.ac.in)
\newline
$^\dagger$Department of Humanities and Sciences, Indian Institute of Petroleum and Energy-Visakhapatnam,
India-530003 (rathans.math@iipe.ac.in)}
\date{\today}
 \maketitle
%
%\thanks{
%$^{\S}$ Department of Humanities and Sciences, Indian Institute of Petroleum and Energy-Visakhapatnam,
%India-530003 ({lavanya_vs@iipe.ac.in})
%\newline
%$^\dagger$ Department of Humanities and Sciences, Indian Institute of Petroleum and Energy-Visakhapatnam,
%India-530003 ({rathans.math@iipe.ac.in})}
%\date{\today}
% \maketitle
\begin{abstract}
In this work, we construct a fifth-order weighted essentially non-oscillatory (WENO) scheme with exponential approximation space for solving dispersive equations. A conservative third-order derivative formulation is developed directly using WENO spatial reconstruction procedure and third-order TVD Runge-Kutta scheme is used  for the evaluation of time derivative. This exponential approximation space consists a tension parameter that may be optimized to fit the specific feature of the charecteristic data, yielding better results without spurious oscillations compared to the polynomial approximation space.  A detailed formulation is  presented for the the construction of conservative flux approximation, smoothness indicators, nonlinear weights and  verified that the proposed scheme provides the required fifth convergence order. One and two-dimensional numerical examples are presented to support the theoretical claims
\end{abstract}



\pagestyle{myheadings} \thispagestyle{plain} \markboth{Lavaya V Salian, Samala Rathan}{ Exponential approximation space reconstruction WENO scheme for dispersive PDEs}

\section{Introduction}
\label{sec:1}
In this paper, a conservative finite difference fifth-order weighted essentially non-oscillatory (WENO) scheme is designed using exponential polynomials to approximate the third-order derivative term, which evolved in the nonlinear, possibly prototype dispersion equations of the form,
\begin{equation} 
\begin{split}
u_t+f(u)_x + g(u)_{xxx}&=0, \quad (x,t) \in \Omega \times (0,T],\\
u(x,0) &= u_0 (x),
\end{split}
\label{eqn:A}
\end{equation}
where $f(u)$ and $g(u)$ are arbitrary (may smooth) nonlinear functions. The simplest form, where dispersion and nonlinearity occur, is the famous Korteweg de-Vries (KdV) equation,
\begin{equation}
\label{KdV}
u_t+a u u_x + u_{xxx} =0, \quad (x,t) \in \Omega \times (0,T],\\
\end{equation}
where $u_t$ is the time evolution of the wave propagating in one direction,
$uu_x$  is the nonlinearity accounting for the steeping of the wave and
$u_{xxx}$ is the linear dispersion, i.e.\ spreading of the wave.
Solitary waves are localized traveling waves that maintain a constant speed and shape, diminishing to asymptotic zero amplitude at a considerable distance \cite{WAZ}. The solutions for solitary waves are postulated in the form of $u(x,t) = k(x-ct)$, with $c$ representing the propagation speed of the wave and $k(z)$, $k^{'}(z)$, $k^{''}(z)$ $\to 0$  as $z \to \pm \infty$, $z=x-ct$. Pioneering work by Zabusky and Kruskal \cite{ZK}  unveiled that solitary waves engage in nonlinear interactions, resulting in an emerging wave that conserves both its shape and amplitude. These distinctive solitary waves are denoted as \textit{solitons}. There is an increase in research efforts worldwide across a variety of scientific fields to investigate the concept of solitons. This concept of solitons has engendered amplified global research efforts spanning diverse scientific domains. Notably, disciplines encompassing fluid dynamics, astrophysics, plasma physics, magneto-acoustic waves, and numerous others have been captivated by the exploration of solitons due to their pronounced significance. This is owing to the genesis of solitons from the intricate equilibrium between weak nonlinearity and dispersion.
\par Later, a class of solitary waves with compact support that is termed as \textit{compactons} were discovered by Rosenau and Hyman \cite{RH}. The general form of the Rosenau - Hyman equation, also known as $K(n,n)$ equation, is of the form
\begin{equation}
\label{eqn0b}
u_t +(u^n)_x + (u^n)_{xxx} =0, \quad n >1,\\
\end{equation}
where $u(x, t)$ is the wave amplitude as a function of the spatial variable $x$ and time $t$. Compactons represent a distinct subclass of solitons, distinguished by their finite wavelength, absence of exponential tails, lack of infinite wings, and resilient adherence to soliton-like behavioral patterns. Compactons and soliton solutions of the KdV equation exhibit shared attributes. For instance, the velocity of an individual compacton is directly proportional to its amplitude. During the movement and interaction of multiple compactons with varying velocities, nonlinear effects come into play, resulting in an altered phase upon exit. In contrast to conventional KdV solitons, the width of a compacton remains constant regardless of its amplitude.  The $K(n,n)$ equation diverges from the typical energy conservation laws held by the KdV equation and can only be derived from a first-order Lagrangian for $n=1$. While solitons are analytical solutions, compactons are non-analytical solutions distinguished by non-analytical points at their edges. These points of non-analyticity align with genuine nonlinearity points within the differential equations.

\par The resolution of the Rosenau-Hyman equation is challenging due to the concurrent interplay of dispersion effects and nonlinearity. Although the extensively employed pseudo-spectral method \cite{CL}  in spatial domains effectively preserves solution positivity and incorporates high-pass filters to induce artificial dissipation, the post-compacton collision can lead to sign alterations within the solution. A variety of alternate methodologies have been explored in this context. These include the finite difference method with Pade approximation \cite{MCCS}, the local discontinuous Galerkin method \cite{LSY}, the second-order finite difference approach \cite{IT}, the adaptive mesh refinement-based line method \cite{SWSZ}, and and the direct WENO scheme employing polynomial bases for dispersion-type equations \cite{AQ}. The governing equation dictating dispersive waves, as denoted by (\ref{eqn:A}), bears significant resemblances to hyperbolic conservation laws. Conspicuously, both equation classes are susceptible to sharp fronts and wavefronts propagating at finite velocities. Consequently, an avenue of inquiry involves extending numerical techniques originally devised for resolving hyperbolic conservation laws, such as the WENO technique, to contend with the intricacies inherent to the dispersive wave equation. However, this extension mandates a meticulous adaptation of the WENO procedure to ensure the preservation of conservation, accuracy, and nonoscillatory behavior.
\par
In 1994, Liu et al.\cite{LOC} introduced the original version of the WENO scheme within the context of a finite volume framework, employed for the resolution of one-dimensional conservation law equations. Subsequently, in 1996, Jiang and Shu \cite{JS} introduced an enhanced rendition of the WENO scheme within the finite differences framework, exhibiting greater efficiency compared to the method presented by Liu et al.\cite{LOC} for addressing both one-dimensional and two-dimensional conservation law equations. Despite the evident enhancements introduced by Jiang and Shu \cite{JS}, their methodology exhibited certain limitations, primarily manifesting in instances where gradients of higher order became negligible, thereby resulting in a diminution of accuracy order. Numerous researchers have endeavored to mitigate the numerical dissipation near discontinuities and to optimize the computational efficiency of the conventional WENO-JS scheme through various adaptations. The authors of \cite{HAP}, Henrick et. al., proposed a fifth-order WENO method named WENO-M. Notably, WENO-M employs a mapping technique to minimize the deviation of nonlinear weights utilized in the convex combination of stencils from the optimal weights, barring instances involving pronounced discontinuities. A distinct variation of the fifth-order WENO scheme, known as WENO-Z, was proposed by Borges et al. in their work \cite{BCCD}. The weighting formulation employed in this context diverges slightly from that of the WENO-JS formulation. In \cite{HKLY}, Ha et al. presented a new smoothness indicator, referred to as WENO-NS, which utilizes the $L^1$-sum of generalized undivided differences to approximate the derivatives of flux functions. This indicator allows for achieving fifth-order convergence even in smooth regions, including critical points where the first derivatives are zero. Further, Rathan and Naga Raju enhanced the accuracy of WENO scheme at critical points for fifth and seventh order schemes \cite{SR1, SR2, SR3}. Various researchers have proposed different high-order WENO schemes with the aim of enhancing the efficiency of these numerical schemes for solving hyperbolic conservation laws \cite{BS, CCD, HKYY}. In 2016, Ha et al.\cite{HKYY} presented a novel WENO scheme utilizing exponential polynomials for solving hyperbolic conservation laws. 
%After that, Abedian et al.\cite{AD} (2022) developed a WENO scheme using exponential polynomials to solve nonlinear degenerate parabolic equations with high-order accuracy and non-oscillatory solutions.

\par The comparative analysis of numerical outcomes between methodologies based upon exponential, trigonometric and algebraic polynomial constructions for hyperbolic conservation laws \cite{HKYY, HKYY2, HLY, ZQ} and non-linear degenerate parabolic equations \cite{AD, RRJ, RJ} has been studied in literature. In instances involving interpolation of data manifesting rapid gradients or high oscillations, the utilization of exponential or trigonometric polynomial bases confers a superior degree of efficiency than algebraic polynomial bases. These alternative polynomial bases prove to be more suitable for accurately capturing and representing such complex features in the data. The primary objective of this study is to propose a novel fifth-order Weighted Essentially Non-Oscillatory scheme named WENO-E, which employs exponential polynomials for solving the dispersion equation. The key design strategy of the WENO-E scheme is to attain the highest possible approximation order in smooth regions while ensuring accuracy is maintained even at critical points. By utilizing exponential polynomials, the scheme aims to enhance the accuracy and robustness of the numerical solution, enabling more effective and reliable computations for the dispersion equation. A global smoothness indicator based on generalized undivided differences is introduced to aid in the design of the nonlinear weights that play a crucial role in WENO reconstructions. Numerical experiments are conducted and compared with the polynomial-based scheme (WENO-Z) to demonstrate the WENO-E scheme's ability to accurately approximate solutions near singularities. This work is the first of its kind for solving prototype dispersive equations (\ref{eqn:A}) with non polynomial approximation space, as per the authors'  best knowledge.

 \par
The paper's organization is delineated as follows: Section 2 elaborates on a broad framework for the finite differences WENO scheme and introduces the concept of approximating numerical fluxes with exponential polynomials. This approach is utilized to effectively address the third-order derivative term inherent in the prototype dispersion equation. In Section 3, we introduce an innovative technique that leverages exponential polynomials to construct numerical fluxes, both for the large stencil and its sub-stencils. We also discuss the associated ideal weights. In Section 4, the paper outlines the use of $L^1$-norm smoothness indicators to construct non-linear weights and provides an analysis of their accuracy in both smooth regions and at critical points. To support our theoretical claim, we present few one and two dimensional numerical results in Section 5. In Section 6, we provide brief concluding remarks.

%--------------------------------------------------------------------------------------

\section{Finite difference WENO scheme}
\label{sec:2}
In this section, we describe a general framework of finite difference WENO schemes based on exponential polynomials to solve prototype dispersion equations. Without loss of generality, we shall focus on the one-dimensional prototype dispersion equations of the form
\begin{equation}
\label{eqn:1}
\begin{split}
u_t+f(u)_x+g(u)_{xxx}&=0, \quad (x,t) \in \Omega \times (0,T],\\
u(x,0) &= u_0 (x),
\end{split}
\end{equation}
along with periodic boundary conditions. To extend the algorithm to higher dimensions, the 1-D algorithm is applied along each coordinate direction. Assume the uniform spatial mesh as follows:
\begin{equation}
\label{eqn:2}
\left\{
\begin{array}{cc}
x \in [x_l, x_r],\quad x_i=x_l +(i-1)\Delta x,\quad i=1:N, \quad, x_1=x_l, \quad x_N=x_r,\\
I_i=[x_{i-\frac{1}{2}}, x_{i+\frac{1}{2}}], \quad x_{i+\frac{1}{2}} = \frac{(x_{i+1}+x_i)}{2},
\end{array}
\right.
\end{equation}
where $\Delta x$ is the spatial step size and $u_i^n$ is defined as a nodal point value $u(x_i,t^n)$. 
%-------------------------------------------------------------------------------------
\subsection{Formulation of WENO scheme}
We use a conservative finite difference scheme in a method of lines (MOL) approach to writing
\begin{equation*}
    \frac{du_i}{dt} = -\frac{1}{\Delta x} ({\hat{F}_{i+\frac{1}{2}}} - {\hat{F}_{i-\frac{1}{2}}}) -\frac{1}{\Delta x^{3}} ({\hat{G}_{i+\frac{1}{2}}} - {\hat{G}_{i-\frac{1}{2}}}),
\end{equation*}
where $\hat{F}$ and $\hat{G}$ are the numerical flux for convection and dispersion, respectively. The detailed construction procedure for convection flux is given in \cite{JS}. For the dispersion term, the $k^{\text{th}}$-order accurate conservative finite difference scheme is as follows:
\begin{equation}
    \frac{1}{\Delta x^{3}} ({\hat{G}_{i+\frac{1}{2}}} - {\hat{G}_{i-\frac{1}{2}}}) \approx g(u)_{xxx}\big|_{x=x_i} + \mathcal{O}(\Delta x^k),
    \label{eqn:(3)}
\end{equation}
where ${\hat{G}_{i+\frac{1}{2}}}$ is a numerical dispersion flux at the cell boundary $x_{i+\frac{1}{2}}$, ${\hat{G}_{i+\frac{1}{2}}} = {\hat{G} (u_{i-r}, \dots ,u_{i+s})}$ of $(r+s)$ variables. It has to satisfy Lipschitz continuity in each of its arguments and is consistent with the physical flux ${\hat{G} (u, \dots ,u)}=g(u)$. Ahmat et al. in \cite{AQ} implicitly considered a function $h(x)$ as follows to guarantee the conservative property
\begin{equation}
   g(u)=\frac{1}{\Delta x^3}  \int_{x-\frac{\Delta x}{2}}^{x+\frac{\Delta x}{2}} \int_{\eta-\frac{\Delta x}{2}}^{\eta+\frac{\Delta x}{2}}\int_{\zeta-\frac{\Delta x}{2}}^{\zeta+\frac{\Delta x}{2}} h(\theta) \,d\theta d\zeta d\eta ,
   \label{eqn:3int}
\end{equation}
then by triple derivation, we have
\begin{equation}
    g(u)_{xxx} = \frac{h(x+\frac{3}{2}\Delta x) -3h(x+\frac{1}{2}\Delta x)+3h(x-\frac{1}{2}\Delta x)-h(x-\frac{3}{2}\Delta x)}{\Delta x^3}.
\end{equation}
If we define a function $G(x)$, such that
\begin{equation}
    G(x)=h(x+\Delta x) - 2h(x) + h(x-\Delta x),
    \label{eqn:(6)}
\end{equation}
then we have
\begin{equation}
    g(u)_{xxx}\big|_{x=x_i} = \frac{G(x_{i+\frac{1}{2}})-G(x_{i-\frac{1}{2}})}{\Delta x^3}.
    \label{eqn:(7)}
\end{equation}
If we have numerical flux ${\hat{G}_{i+\frac{1}{2}}}$ which is an approximation to $G(x_{i+\frac{1}{2}})$ up to $k^{\text{th}}$-order, then we have equation (\ref{eqn:(3)}). 
\par
We now derive the specific form of the numerical flux ${\hat{G}(u)}$ for WENO5. We first split the flux into positive and negative parts, that is $g(u)=g^+(u) +g^-(u)$ with $\frac{\partial g^+(u)}{\partial u} \ge 0$ and  $\frac{\partial g^-(u)}{\partial u} \le 0$, and $G^+(x)$ and $G^-(x)$ are defined  by (\ref{eqn:(6)}) according to $g^+(u)$ and $g^-(u)$, respectively. The reconstructed numerical fluxes ${\hat{G}^+_{i+\frac{1}{2}}}$ and ${\hat{G}^-_{i+\frac{1}{2}}}$ to approximate $G^+(x_{i+\frac{1}{2}})$ and $G^-(x_{i+\frac{1}{2}})$ up to $k^{\text{th}}$-order, respectively. Finally, we define numerical fluxes ${\hat{G}_{i+\frac{1}{2}}}$ as
\begin{equation}
{\hat{G}_{i+\frac{1}{2}}} ={\hat{G}^+_{i+\frac{1}{2}}}+{\hat{G}^-_{i+\frac{1}{2}}},
\label{eqn:(8)}
\end{equation}
and ${\hat{G}_{i+\frac{1}{2}}}$ is an approximation to $G(x_{i+\frac{1}{2}})$ up to $k^{\text{th}}$-order. The reconstruction procedure for ${\hat{G}^+_{i+\frac{1}{2}}}$ is given below and procedure for ${\hat{G}^-_{i+\frac{1}{2}}}$ is mirror symmetric respect to $x_{i+\frac{1}{2}}$.
\begin{figure}
\begin{center}
\minipage{0.75\textwidth}
  \includegraphics[width=\linewidth]{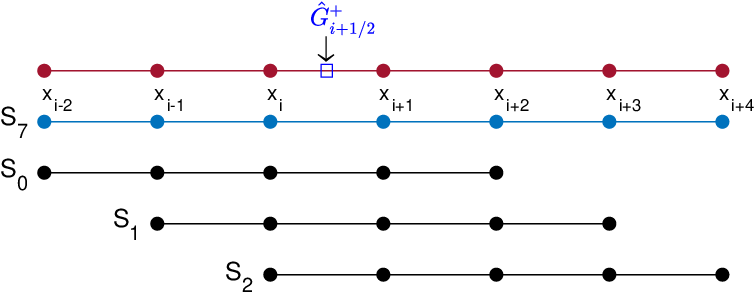}
\endminipage\hfill
\caption{ The numerical flux ${\hat{G}^+_{i+\frac{1}{2}}}$ is constructed on the stencil $S_7=\{x_{i-2}, \dots , x_{i+4}\}$ with seven uniform points and three 5-point substencils $S_0$, $S_1$, $S_2$. }
\end{center}
\end{figure}
%-----------------------------------------------------------------------------------
\subsection{Approximation using exponential polynomials}
The WENO scheme for solving the dispersion type equation involves approximating the value of $G(x_{i+\frac{1}{2}})$ using equation (\ref{eqn:(7)}). The reconstruction process must balance achieving optimal accuracy in smooth regions while maintaining essential non-oscillatory behaviour in non smooth regions. Polynomials are widely used to construct numerical fluxes. However, the limitation of using polynomials is that the approximation space is shift-and-scale invariant. Thus, it cannot be tailored to suit the specific characteristics of the given data. This limitation can result in significant numerical dissipation when interpolating data with rapid gradients, which hinders the ability to produce sharp edges. To overcome this limitation, researchers \cite{HKYY, HKYY2, HLY, ZQ, AD} have explored the use of other types of basis functions, such as exponential polynomials, which have been shown to yield better results in terms of producing sharp edges and reducing numerical dissipation.\\
The general form of exponential polynomials can be written as:
\begin{equation}
    \Phi(x) = x^n e^{\lambda x},
    \label{eqn:phi}
\end{equation}
where $n$ is a non-negative integer and $\lambda \in \mathbb{R}$ or $\lambda \in \iota\mathbb{R}$ $(\iota^2 =-1)$. Let $\{ \Phi_1, \dots, \Phi_r \}$ be a set of exponential polynomials of the form in (\ref{eqn:phi}) and let $\Gamma_r$ be the space defined by
\begin{equation}
    \Gamma_r := \text{span} \{\Phi_1, \dots, \Phi_r \}.
\end{equation}
A necessary condition for the space $\Gamma_r$ is that the
set $\{ \Phi_1, \dots, \Phi_r \}$ is linearly independent so that the determinants of the Wronskian matrix related to them are non-zero, i.e.,
\begin{equation}
    \det(\Phi_n(s_i) : i,n=1, \dots, r) \neq 0,
    \label{eqn:det}
\end{equation}
for any $r$-point stencil $\{s_i : i=1, \dots, r \}$.
The space $\Gamma_r$ needs to satisfy the following basic requirements for the practical computation of the proposed interpolation:
\begin{itemize}
    \item Shift-invariant. The space $\Gamma_r$ should be shift-invariant in the sense that for any $\alpha \in \mathbb{R}$, $f \in \Gamma_r$ implies $f(\cdot - \alpha) \in \Gamma_r$. This ensures that the interpolation kernel is invariant under the shifting of the evaluation location and stencil. It allows for a set of interpolation kernels to be precomputed for a fixed point and then applied to every evaluation position based on the chosen stencil at a given cell boundary. The space $\Gamma_r$ defined in (\ref{eqn:det}) meets this requirement since it is shift-invariant.
    \item Symmetry. The space $\Gamma_r$ should be symmetry means that if a function $f$ is in the space  $\Gamma_r$, then the reflected function $f(-\cdot)$ should also belong to the same space.
\end{itemize}
Including the polynomial $\Phi(x)=1$  in the space $\Gamma_r$ ensures that the sum of the interpolation weights over all basis functions is equal to 1, which is necessary for the interpolation kernel to satisfy the partition of unity property. In this paper, we choose
\begin{equation}\label{S}
 \Gamma_7 := \text{span}\{1, x, x^2, e^{\lambda x}, e^{-\lambda x}, \cos {\lambda x}, \sin {\lambda x}\}
\end{equation}
as the basis functions for global stencil $S_7 := \{x_{i-2},\dots, x_{i+4}\}$ and similarly, 
\begin{equation}\label{Sm}
\Gamma_5 = \text{span}\{1, x, x^2, e^{\lambda x}, e^{-\lambda x}\}
\end{equation}
for the five-point substencils  $ S_m := \{x_{i-2+m},\dots,x_{i+2+m}\}$, $m=0,1,2$. Here, $\Gamma_7$ and $\Gamma_5$ constitute an \textit{extended Tchebysheff systems} on $\mathbb{R}$ so that the non-singularity of the interpolation matrices in (\ref{eqn:det}) is guaranteed as in \cite{KS, HKYY}.
%----------------------------------------------------------------------------------
\section{Fifth-order WENO-E scheme}
Conservative numerical schemes are developed by approximating the function $G(x)$ in (\ref{eqn:(7)}). This approximation, represented as $\hat{G}(x)$, is constructed using a exponential polynomial form with unspecified coefficients. When this polynomial is substituted into (\ref{eqn:3int}), it results in a system of equations. In this system, the flux is known at the nodes surrounding the relevant interface, enabling the determination of a distinct set of coefficients. After obtaining $\hat{G}(x)$, the approximation of the spatial derivative in (\ref{eqn:(3)}) is as follows
\begin{equation}\label{eqn:3.1}
    g(u)_{xxx}\big|_{x=x_i} \approx \frac{1}{\Delta x^{3}} ({\hat{G}_{i+\frac{1}{2}}} - {\hat{G}_{i-\frac{1}{2}}}).
\end{equation}
We will focus on two orders of convergence. The primary concern is the order at which (\ref{eqn:3.1}) is satisfied, as it dictates the spatial convergence rate of the overall scheme. Additionally, we consider the order of individual approximations to the numerical flux, ${\hat{G}_{i\pm\frac{1}{2}}}$, which is significant in establishing criteria for the acceptance of non-oscillatory weights.
%-----------------------------------------------------------------------------------
\subsection{Constructions of numerical flux ${\hat{G}^+_{i+\frac{1}{2}}}$}
The approximations  of ${\hat{G}^+_{i+\frac{1}{2}}}$ are denoted by $p(x)$. We first consider an exponential polynomial approximation to $h(x)$ on the 7-point stencil $S_7$
\begin{equation}
\label{eqn:Epoly}
    h(x) \approx q(x) = a_0 + a_1 x + a_2 x^2 + a_3 e^{\lambda x}+a_4 e^{-\lambda x}+a_5 \cos {\lambda x}+a_6 \sin {\lambda x},
\end{equation}
with undetermined coefficients $a_k$ with $k=0,\dots,6$.\\
Substituting (\ref{eqn:Epoly}) into (\ref{eqn:3int})
\begin{equation}
   g^+(u(x))=\frac{1}{\Delta x^3}  \int_{x-\frac{\Delta x}{2}}^{x+\frac{\Delta x}{2}} \int_{\eta-\frac{\Delta x}{2}}^{\eta+\frac{\Delta x}{2}}\int_{\zeta-\frac{\Delta x}{2}}^{\zeta+\frac{\Delta x}{2}} q(\theta) \,d\theta d\zeta d\eta ,
\end{equation}
and performing the integration  gives
\begin{equation}
\begin{split}
    g^+(u(x))&= a_0 + a_1 x + a_2 \biggl(\frac{\Delta x^2}{4}+x^2\biggr)+a_3 \biggl(-\frac{4 e^{\lambda x} \sinh\left[\frac{\lambda\Delta x}{2}\right]}{\lambda^3\Delta x^3}+\frac{4 e^{\lambda x}\cosh[\lambda\Delta x] \sinh\left[\frac{\lambda\Delta x}{2}\right]}{\lambda^3\Delta x^3} \biggr)\\
   & +a_4 \biggl( \frac{e^{\frac{\lambda}{2}(\Delta x-2x)}(-1+e^{\lambda\Delta x})}{\lambda^3\Delta x^3}+\frac{e^{\frac{-3\lambda\Delta x}{2}-\lambda x}(-1+e^{\lambda\Delta x})}{\lambda^3\Delta x^3}-\frac{4e^{-\lambda x}\sinh\left[\frac{\lambda\Delta x}{2}\right]}{\lambda^3\Delta x^3}\biggr)\\
   &+a_5 \biggl(\frac{8 \cos[\lambda x]\sin{\left[\frac{\lambda\Delta x}{2}\right]}^3}{\lambda^3\Delta x^3}\biggr)+a_6 \biggl(\frac{8 \sin[\lambda x]\sin{\left[\frac{\lambda\Delta x}{2}\right]}^3}{\lambda^3\Delta x^3}\biggr).
\end{split} \label{eqn:11}
\end{equation}
To determine the coefficients $a_0,\dots,a_6$, one can consider (\ref{eqn:11}) as $g^+(u(x_{i}-2\Delta x))=g^+(u_{i-2}),\dots,g^+(u(x_{i}+4\Delta x))=g^+(u_{i+4})$ with $x_i = 0$ and solve the resulting $7 \times 7$ system $AX=B$ are specified in Appendix:B. Substituting the coefficients of vector $X$ into (\ref{eqn:Epoly}) and then calculating $p(x)=q(x+\Delta x) - 2q(x) + q(x-\Delta x)$ at $x=x_{i+\frac{1}{2}}$, we get
\begin{equation}
 {\hat{G}^+_{i+\frac{1}{2}}} = C_0 g^+(u_{i-2}) + C_1 g^+(u_{i-1}) + C_2 g^+(u_{i}) + C_3 g^+(u_{i+1}) + C_4 g^+(u_{i+2}) + C_5 g^+(u_{i+3}) + C_6 g^+(u_{i+4}).
 \label{eqn:G_half}
\end{equation}
After applying the Taylor series to the coefficients $C_j, 0\le j \le 6$ of equation (\ref{eqn:G_half}) is given by
\begin{equation}
    \begin{split}
        C_0 &= -\frac{1}{15}-\frac{229 \lambda^4 \Delta x^4}{226800}+\mathcal{O}(\Delta x^8) , \quad C_1 = \frac{21}{40}+\frac{1493\lambda^4 \Delta x^4}{75600}+\mathcal{O}(\Delta x^8),\\
        C_2 &= \frac{1}{8}-\frac{587 \lambda^4 \Delta x^4}{15120}+\mathcal{O}(\Delta x^8) ,\quad \quad C_3 = -\frac{23}{12}-\frac{65\lambda^4 \Delta x^4}{9072}+\mathcal{O}(\Delta x^8),\\
        C_4 &= \frac{7}{4}+\frac{11\lambda^4 \Delta x^4}{378}+\mathcal{O}(\Delta x^8) ,\quad \quad C_5 = -\frac{19}{40}-\frac{323\lambda^4 \Delta x^4}{18900}+\mathcal{O}(\Delta x^8),\\
        C_6 &= -\frac{7}{120}+\frac{103 \lambda^4 \Delta x^4}{113400}+\mathcal{O}(\Delta x^8).
    \end{split}
\end{equation}
Thus, the equation (\ref{eqn:G_half}) becomes
\begin{equation*}
\begin{aligned}
 {\hat{G}^+_{i+\frac{1}{2}}} =&   \biggl(-\frac{1}{15}-\frac{229 \lambda^4 \Delta x^4}{226800}\biggr) g^+(u_{i-2}) + \biggl(\frac{21}{40}+  
                                               \frac{1493\lambda^4 \Delta x^4}{75600}\biggr) g^+(u_{i-1})+ \biggl(\frac{1}{8}-\frac{587 \lambda^4 \Delta x^4}{15120}\biggr) g^+(u_{i}) \\
                                               & + \biggl(-\frac{23}{12}-\frac{65\lambda^4 \Delta x^4}{9072}\biggr) g^+(u_{i+1}) + \biggl(\frac{7}{4}+\frac{11\lambda^4 \Delta x^4}{378}\biggr) g^+(u_{i+2}) \\
                                               & + \biggl(-\frac{19}{40}-\frac{323\lambda^4 \Delta x^4}{18900}\biggr) g^+(u_{i+3}) + \biggl( -\frac{7}{120}+\frac{103 \lambda^4 \Delta x^4}{113400}\biggr) g^+(u_{i+4}).
   \end{aligned}
 \label{eqn:G_Half}
\end{equation*}
From the definition  of (\ref{eqn:3int}) we also know that
\begin{equation*}
\begin{split}
 {\hat{G}^+_{i+\frac{1}{2}}} &=\frac{1}{\Delta x^3}\Biggl[ \biggl(-\frac{1}{15}-\frac{229 \lambda^4 \Delta x^4}{226800}\biggr) \biggl( \int_{x_{i-\frac{5}{2}}}^{x_{i-\frac{3}{2}}} \int_{\eta-\frac{\Delta x}{2}}^{\eta+\frac{\Delta x}{2}}\int_{\zeta-\frac{\Delta x}{2}}^{\zeta+\frac{\Delta x}{2}} h(\theta) \,d\theta d\zeta d\eta\biggr) \\
 & + \biggl(\frac{21}{40}+\frac{1493\lambda^4 \Delta x^4}{75600}\biggr) \biggl( \int_{x_{i-\frac{3}{2}}}^{x_{i-\frac{1}{2}}} \int_{\eta-\frac{\Delta x}{2}}^{\eta+\frac{\Delta x}{2}}\int_{\zeta-\frac{\Delta x}{2}}^{\zeta+\frac{\Delta x}{2}} h(\theta) \,d\theta d\zeta d\eta\biggr)\\
 &+ \biggl(\frac{1}{8}-\frac{587 \lambda^4 \Delta x^4}{15120}\biggr) \biggl( \int_{x_{i-\frac{1}{2}}}^{x_{i+\frac{1}{2}}} \int_{\eta-\frac{\Delta x}{2}}^{\eta+\frac{\Delta x}{2}}\int_{\zeta-\frac{\Delta x}{2}}^{\zeta+\frac{\Delta x}{2}} h(\theta) \,d\theta d\zeta d\eta\biggr) \\
  & + \biggl(-\frac{23}{12}-\frac{65\lambda^4 \Delta x^4}{9072}\biggr) \biggl( \int_{x_{i+\frac{1}{2}}}^{x_{i+\frac{3}{2}}} \int_{\eta-\frac{\Delta x}{2}}^{\eta+\frac{\Delta x}{2}}\int_{\zeta-\frac{\Delta x}{2}}^{\zeta+\frac{\Delta x}{2}} h(\theta) \,d\theta d\zeta d\eta\biggr)\\ 
 &+ \biggl(\frac{7}{4}+\frac{11\lambda^4 \Delta x^4}{378}\biggr) \biggl( \int_{x_{i+\frac{3}{2}}}^{x_{i+\frac{5}{2}}} \int_{\eta-\frac{\Delta x}{2}}^{\eta+\frac{\Delta x}{2}}\int_{\zeta-\frac{\Delta x}{2}}^{\zeta+\frac{\Delta x}{2}} h(\theta) \,d\theta d\zeta d\eta\biggr) \\
 & + \biggl(-\frac{19}{40}-\frac{323\lambda^4 \Delta x^4}{18900}\biggr) \biggl( \int_{x_{i+\frac{5}{2}}}^{x_{i+\frac{7}{2}}} \int_{\eta-\frac{\Delta x}{2}}^{\eta+\frac{\Delta x}{2}}\int_{\zeta-\frac{\Delta x}{2}}^{\zeta+\frac{\Delta x}{2}} h(\theta) \,d\theta d\zeta d\eta\biggr)\\
 &+ \biggl( -\frac{7}{120}+\frac{103 \lambda^4 \Delta x^4}{113400}\biggr) \biggl( \int_{x_{i+\frac{7}{2}}}^{x_{i+\frac{9}{2}}} \int_{\eta-\frac{\Delta x}{2}}^{\eta+\frac{\Delta x}{2}}\int_{\zeta-\frac{\Delta x}{2}}^{\zeta+\frac{\Delta x}{2}} h(\theta) \,d\theta d\zeta d\eta\biggr) \Biggr].
 \end{split}
\end{equation*}
To analyze the accuracy of the approximation to the third derivative for smooth solutions, we assume $h(x)$ has sufficient local regularity. Substituting the Taylor series expansion at the grid point $x_i$,
\begin{equation*}
    h(\theta) = h(x_i) + \sum_{j=1}^7 \frac{(\theta-x_i)^j}{j!}\frac{d^j h}{dx^j}\bigg|_{x=x_i} + \mathcal{O}(\Delta x^8),
\end{equation*}
and performing the integration, we have
\begin{equation}
\begin{split}
    {\hat{G}^+_{i+\frac{1}{2}}} &= \frac{d^2 h}{dx^2}\bigg|_{x=x_i}\Delta x^2 +\frac{1}{2}\frac{d^3 h}{dx^3}\bigg|_{x=x_i}\Delta x^3 +\frac{5}{24}\frac{d^4 h}{dx^4}\bigg|_{x=x_i}\Delta x^4 +\frac{1}{16}\frac{d^5 h}{dx^5}\bigg|_{x=x_i}\Delta x^5+\frac{91}{5760}\frac{d^6 h}{dx^6}\bigg|_{x=x_i}\Delta x^6\\
    & +\biggl(-\frac{7}{240}\lambda^4 \frac{d^3 h}{dx^3}\bigg|_{x=x_i}+\frac{25}{768} \frac{d^7 h}{dx^7}\bigg|_{x=x_i}\biggr)\Delta x^7+\mathcal{O}(\Delta x^8).
\end{split}  
\label{eqn:15}
\end{equation}
Comparing  (\ref{eqn:15}) with Taylor series expansion
\begin{equation}
\begin{split}
    G(x_{i+\frac{1}{2}})=&h(x_{i+\frac{3}{2}}) - 2h(x_{i+\frac{1}{2}}) + h(x_{i-\frac{1}{2}})\\
    =&\frac{d^2 h}{dx^2}\bigg|_{x=x_i}\Delta x^2 +\frac{1}{2}\frac{d^3 h}{dx^3}\bigg|_{x=x_i}\Delta x^3 +\frac{5}{24}\frac{d^4 h}{dx^4}\bigg|_{x=x_i}\Delta x^4 +\frac{1}{16}\frac{d^5 h}{dx^5}\bigg|_{x=x_i}\Delta x^5\\
    & +\frac{91}{5760}\frac{d^6 h}{dx^6}\bigg|_{x=x_i}\Delta x^6+\frac{13}{3840}\frac{d^7 h}{dx^7}\bigg|_{x=x_i}\Delta x^7+\mathcal{O}(\Delta x^8),
    \label{eqn:16}
\end{split}
\end{equation}
we obtain 
\begin{equation}
\begin{split}
    {\hat{G}^+_{i+\frac{1}{2}}}&= G(x_{i+\frac{1}{2}})+\biggl(-\frac{7}{240}\lambda^4 \frac{d^3 h}{dx^3}\bigg|_{x=x_i}+\frac{7}{240} \frac{d^7 h}{dx^7}\bigg|_{x=x_i}\biggr)\Delta x^7+\mathcal{O}(\Delta x^8),\\
 &= G(x_{i+\frac{1}{2}})+\mathcal{A}^+\Delta x^7+\mathcal{O}(\Delta x^8).   
 \end{split} \label{eqn:17}
\end{equation}
Similarly, we calculate the numerical flux ${\hat{G}^+_{i-\frac{1}{2}}}$
\begin{equation}
\begin{split}
{\hat{G}^+_{i-\frac{1}{2}}} &= \biggl(-\frac{1}{15}-\frac{229 \lambda^4 \Delta x^4}{226800}\biggr) g^+(u_{i-3}) + \biggl(\frac{21}{40}+\frac{1493\lambda^4 \Delta x^4}{75600}\biggr) g^+(u_{i-2}) + \biggl(\frac{1}{8}-\frac{587 \lambda^4 \Delta x^4}{15120}\biggr) g^+(u_{i-1}) \\
 &+ \biggl(-\frac{23}{12}-\frac{65\lambda^4 \Delta x^4}{9072}\biggr) g^+(u_{i})+ \biggl(\frac{7}{4}+\frac{11\lambda^4 \Delta x^4}{378}\biggr) g^+(u_{i+1}) + \biggl(-\frac{19}{40}-\frac{323\lambda^4 \Delta x^4}{18900}\biggr) g^+(u_{i+2})\\
 & + \biggl( -\frac{7}{120}+\frac{103 \lambda^4 \Delta x^4}{113400}\biggr) g^+(u_{i+3}) \\
 {\hat{G}^+_{i-\frac{1}{2}}}&= G(x_{i-\frac{1}{2}})+\biggl(-\frac{7}{240}\lambda^4 \frac{d^3 h}{dx^3}\bigg|_{x=x_i}+\frac{7}{240} \frac{d^7 h}{dx^7}\bigg|_{x=x_i}\biggr)\Delta x^7+\mathcal{O}(\Delta x^8)= G(x_{i-\frac{1}{2}})+\mathcal{A}^-\Delta x^7+\mathcal{O}(\Delta x^8),
 \end{split}
 \label{eqn:18}
\end{equation}
where the values $\mathcal{A}^{\pm}$ are independent of $\Delta x$, but dependent on $\lambda$. Therefore, substituting (\ref{eqn:17}) and (\ref{eqn:18}) into (\ref{eqn:(7)}), we have the fifth order approximation 
\begin{equation}
    g(u)_{xxx}\big|_{x=x_i} = \frac{{\hat{G}^+_{i+\frac{1}{2}}}-{\hat{G}^+_{i-\frac{1}{2}}}}{\Delta x^3}+\mathcal{O}(\Delta x^5).
    \label{eqn:(19)}
\end{equation}
%-----------------------------------------------------------------------------------
\subsection{Constructions of numerical flux ${\hat{G}^{m}_{i+\frac{1}{2}}}$ for $m=0,1,2$}
When the stencil $S_7$ contains a discontinuity, the numerical flux ${\hat{G}^+_{i+\frac{1}{2}}}$ formulation on the big stencil $S_7$ can result in oscillations because all seven nodes are affected. To mitigate this issue, the WENO procedure considers smaller stencils  $ S_m = \{x_{i-2+m},\dots,x_{i+2+m}\}$, $m=0,1,2$, of the big stencil $S_7$.
Following a similar argument, on each small stencils $ S_m,$ $m=0,1,2$, by taking approximation to $h(x)$ as
\begin{equation}
\label{eqn:Epolysub}
    h(x) \approx q_m(x) = a_0^m + a_1^m x + a_2^m x^2 + a_3^m e^{\lambda x}+a_4^m e^{-\lambda x},
\end{equation}
with undetermined coefficients $a_k^m$ with $k=0, \dots, 4$, $m=0,1,2$. Substituting (\ref{eqn:Epolysub}) into (\ref{eqn:3int})
\begin{equation}
\begin{aligned}
   \frac{1}{\Delta x^3}  \int_{x_j-\frac{\Delta x}{2}}^{x_j+\frac{\Delta x}{2}} \int_{\eta-\frac{\Delta x}{2}}^{\eta+\frac{\Delta x}{2}}\int_{\zeta-\frac{\Delta x}{2}}^{\zeta+\frac{\Delta x}{2}} q_m(\theta) \,d\theta d\zeta d\eta &= g^+(u_j), \quad j=i-2+m, \dots, i+2+m, \quad m=0,1,2,\\
   \end{aligned}
\end{equation}
and performing the integration  gives
\begin{equation}
\begin{split}
    g^+(u(x))&= a_0^m + a_1^m x + a_2^m \biggl(\frac{\Delta x^2}{4}+x^2\biggr)+a_3^m \biggl(\frac{e^{\frac{-3\lambda \Delta x}{2}+x\lambda}(-1+e^{\lambda \Delta x})^3}{\lambda^3 \Delta x^3}\biggr) + a_4^m \biggl(\frac{e^{\frac{-3\lambda \Delta x}{2}-\lambda x}(-1+e^{\lambda \Delta x})^3}{ \lambda^3 \Delta x^3}\biggr).
\end{split} \label{eqn:11b}
\end{equation}
To determine the coefficients $a_0^m,\dots,a_4^m$ we solve the resulting $5 \times 5$ systems $A_m X_m=B_m$ is given by,
\begin{enumerate}
 \item For $m=0$,\\
\begin{eqnarray*}
A_0=
  \left[\begin{matrix} 1 & -2\Delta x & \frac{17}{4} \Delta x^2 & \frac{e^{-\frac{7\lambda\Delta x}{2}}(-1+e^{\lambda\Delta x})^3}{\lambda^3 \Delta x^3 } & \frac{e^{\frac{\lambda \Delta x}{2}}(-1+e^{\lambda \Delta x})^3}{\lambda^3\Delta x^3} \\
1 & -\Delta x & \frac{5}{4} \Delta x^2 & \frac{e^{-\frac{5\lambda\Delta x}{2}}(-1+e^{\lambda \Delta x})^3}{\lambda^3 \Delta x^3} & \frac{e^{-\frac{\lambda \Delta x}{2}}(-1+e^{\lambda \Delta x})^3}{\lambda^3 \Delta x^3} \\
1 & 0 & \frac{1}{4} \Delta x^2 & \frac{e^{-\frac{3\lambda \Delta x}{2}}(-1+e^{\lambda \Delta x})^3}{\lambda^3 \Delta x^3} & \frac{e^{-\frac{3\lambda\Delta x}{2}}(-1+e^{\lambda\Delta x})^3}{\lambda^3 \Delta x^3} \\
1 & \Delta x & \frac{5}{4} \Delta x^2 & \frac{e^{-\frac{\lambda \Delta x}{2}}(-1+e^{\lambda \Delta x})^3}{\lambda^3 \Delta x^3} & \frac{e^{-\frac{5\lambda\Delta x}{2}}(-1+e^{\lambda \Delta x})^3}{\lambda^3 \Delta x^3} \\
1 & 2\Delta x & \frac{17}{4} \Delta x^2 & \frac{e^{\frac{\lambda \Delta x}{2}}(-1+e^{\lambda \Delta x})^3}{\lambda^3 \Delta x^3} & \frac{e^{-\frac{7\lambda \Delta x}{2}}(-1+e^{\lambda \Delta x})^3}{\lambda^3 \Delta x^3} \\
\end{matrix}\right] ,
\end{eqnarray*}
\[B_0=
  \left[\begin{matrix}g^+(u_{i-2}) & g^+(u_{i-1}) & g^+(u_{i}) & g^+(u_{i+1}) & g^+(u_{i+2}) \end{matrix}\right]^T, \]
\[ X_0=
  \left[\begin{matrix}a_0^0 & a_1^0 & a_2^0 & a_3^0 & a_4^0  \end{matrix}\right]^T.
\]
\item  For $m=1$,\\
\begin{eqnarray*}
A_1=
  \left[\begin{matrix} 
1 & -\Delta x & \frac{5}{4} \Delta x^2 & \frac{e^{-\frac{5\lambda \Delta x}{2}}(-1+e^{\lambda \Delta x})^3}{\lambda^3 \Delta x^3} & \frac{e^{-\frac{\lambda \Delta x}{2}}(-1+e^{\lambda \Delta x})^3}{\lambda^3 \Delta x^3} \\
1 & 0 & \frac{1}{4} \Delta x^2 & \frac{e^{-\frac{3\lambda \Delta x}{2}}(-1+e^{\lambda \Delta x})^3}{\lambda^3 \Delta x^3} & \frac{e^{-\frac{3\lambda\Delta x}{2}}(-1+e^{\lambda \Delta x})^3}{\lambda^3 \Delta x^3} \\
1 & \Delta x & \frac{5}{4} \Delta x^2 & \frac{e^{-\frac{\lambda \Delta x}{2}}(-1+e^{\lambda \Delta x})^3}{\lambda^3 \Delta x^3} & \frac{e^{-\frac{5\lambda \Delta x}{2}}(-1+e^{\lambda \Delta x})^3}{\lambda^3 \Delta x^3} \\
1 & 2\Delta x & \frac{17}{4} \Delta x^2 & \frac{e^{\frac{\lambda \Delta x}{2}}(-1+e^{\lambda \Delta x})^3}{\lambda^3 \Delta x^3} & \frac{e^{-\frac{7\lambda \Delta x}{2}}(-1+e^{\lambda \Delta x})^3}{\lambda^3 \Delta x^3} \\
1 & 3\Delta x & \frac{37}{4} \Delta x^2 & \frac{e^{\frac{3\lambda\Delta x}{2}}(-1+e^{\lambda\Delta x})^3}{\lambda^3 \Delta x^3} & \frac{e^{-\frac{9\lambda \Delta x}{2}}(-1+e^{\lambda \Delta x})^3}{\lambda^3 \Delta x^3} \\
\end{matrix}\right] ,
\end{eqnarray*}
\[B_1=
  \left[\begin{matrix} g^+(u_{i-1}) & g^+(u_{i}) & g^+(u_{i+1}) & g^+(u_{i+2}) & g^+(u_{i+3}) \end{matrix}\right]^T,\]
 \[X_1=
  \left[\begin{matrix}a_0^1 & a_1^1 & a_2^1 & a_3^1 & a_4^1  \end{matrix}\right]^T.
\]
\item  For $m=2$,\\
\begin{eqnarray*}
A_2=
  \left[\begin{matrix} 
1 & 0 & \frac{1}{4} \Delta x^2 & \frac{e^{-\frac{3\lambda \Delta x}{2}}(-1+e^{\lambda \Delta x})^3}{\lambda^3 \Delta x^3} & \frac{e^{-\frac{3\lambda \Delta x}{2}}(-1+e^{\lambda \Delta x})^3}{\lambda^3 \Delta x^3} \\
1 & \Delta x & \frac{5}{4} \Delta x^2 & \frac{e^{-\frac{\lambda \Delta x}{2}}(-1+e^{\lambda \Delta x})^3}{\lambda^3 \Delta x^3} & \frac{e^{-\frac{5\lambda \Delta x}{2}}(-1+e^{\lambda \Delta x})^3}{\lambda^3 \Delta x^3} \\
1 & 2\Delta x & \frac{17}{4} \Delta x^2 & \frac{e^{\frac{\lambda \Delta x}{2}}(-1+e^{\lambda \Delta x})^3}{\lambda^3 \Delta x^3} & \frac{e^{-\frac{7\lambda \Delta x}{2}}(-1+e^{\lambda \Delta x})^3}{\lambda^3 \Delta x^3} \\
1 & 3\Delta x & \frac{37}{4} \Delta x^2 & \frac{e^{\frac{3\lambda \Delta x}{2}}(-1+e^{\lambda \Delta x})^3}{\lambda^3 \Delta x^3} & \frac{e^{-\frac{9\lambda \Delta x}{2}}(-1+e^{\lambda \Delta x})^3}{\lambda^3 \Delta x^3} \\
1 & 4\Delta x & \frac{65}{4} \Delta x^2 & \frac{e^{\frac{5\lambda \Delta x}{2}}(-1+e^{\lambda \Delta x})^3}{\lambda^3 \Delta x^3} & \frac{e^{-\frac{11\lambda \Delta x}{2}}(-1+e^{\lambda \Delta x})^3}{\lambda^3 \Delta x^3} \\
\end{matrix}\right] ,
\end{eqnarray*}
\[B_2=
  \left[\begin{matrix}  g^+(u_{i}) & g^+(u_{i+1}) & g^+(u_{i+2}) & g^+(u_{i+3}) & g^+(u_{i+4}) \end{matrix}\right]^T,\]
   \[X_2=
  \left[\begin{matrix}a_0^2 & a_1^2 & a_2^2 & a_3^2 & a_4^2  \end{matrix}\right]^T.
\]
\end{enumerate}
Substituting the obtained coefficients of $X_0,$ $X_1$ and $X_2$ into (\ref{eqn:Epolysub}) and then calculating $p_m(x)=q_m(x+\Delta x) - 2q_m(x) + q_m(x-\Delta x)$ at $x=x_{i+\frac{1}{2}}$, we get
\begin{equation}
\begin{aligned}
 {\hat{G}^{(0)}_{i+\frac{1}{2}}} &= C_{0}^0 g^+(u_{i-2}) + C_{1}^0  g^+(u_{i-1}) + C_{2}^0 g^+(u_{i})+ C_{3}^0 g^+(u_{i+1}) + C_{4}^0 g^+(u_{i+2}),\\
 {\hat{G}^{(1)}_{i+\frac{1}{2}}} &= C_{0}^1 g^+(u_{i-1}) + C_{1}^1  g^+(u_{i}) + C_{2}^1 g^+(u_{i+1})+ C_{3}^1 g^+(u_{i+2}) + C_{4}^1 g^+(u_{i+3}),\\
 {\hat{G}^{(2)}_{i+\frac{1}{2}}} &= C_{0}^2 g^+(u_{i}) + C_{1}^2  g^+(u_{i+1}) + C_{2}^2 g^+(u_{i+2})+ C_{3}^2 g^+(u_{i+3}) + C_{4}^2 g^+(u_{i+4}).
\end{aligned}
\label{eqn:21}
 \end{equation}
 After applying the Taylor series to the coefficients $C_j^m, 0\le j \le 4, 0\le m \le 2$ of equation (\ref{eqn:21}), we get
\begin{eqnarray}
\begin{aligned}
 {\hat{G}^{(0)}_{i+\frac{1}{2}}} &= \biggl(-\frac{1}{4}+\frac{7 \lambda^2 \Delta x^2}{120}-\frac{263\lambda^4 \Delta x^4}{30240}\biggr) g^+(u_{i-2}) + \biggl(\frac{3}{2}-\frac{13 \lambda^2 \Delta x^2}{120}+\frac{97\lambda^4 \Delta x^4}{6048}\biggr) g^+(u_{i-1})\\
 & + \biggl(-2-\frac{ \lambda^2 \Delta x^2}{40}+\frac{41\lambda^4 \Delta x^4}{10080}\biggr) g^+(u_{i})+ \biggl(\frac{1}{2}+\frac{17 \lambda^2 \Delta x^2}{120}-\frac{649\lambda^4 \Delta x^4}{30240}\biggr) g^+(u_{i+1}) \\
 &+ \biggl(\frac{1}{4}-\frac{ \lambda^2 \Delta x^2}{15}+\frac{19\lambda^4 \Delta x^4}{1890}\biggr) g^+(u_{i+2}),\\
 {\hat{G}^{(1)}_{i+\frac{1}{2}}} &= \biggl(\frac{1}{4}-\frac{ \lambda^2 \Delta x^2}{15}+\frac{19\lambda^4 \Delta x^4}{1890}\biggr) g^+(u_{i-1}) + \biggl(\frac{1}{2}+\frac{17 \lambda^2 \Delta x^2}{120}-\frac{649\lambda^4 \Delta x^4}{30240}\biggr) g^+(u_{i}) \\
 & + \biggl(-2-\frac{ \lambda^2 \Delta x^2}{40}+\frac{41\lambda^4 \Delta x^4}{10080}\biggr) g^+(u_{i+1})+ \biggl(\frac{3}{2}-\frac{13 \lambda^2 \Delta x^2}{120}+\frac{97\lambda^4 \Delta \Delta x^4}{6048}\biggr) g^+(u_{i+2})\\
 & + \biggl(-\frac{1}{4}+\frac{ 7\lambda^2 h^2}{120}-\frac{263\lambda^4 \Delta x^4}{30240}\biggr) g^+(u_{i+3}),\\
  {\hat{G}^{(2)}_{i+\frac{1}{2}}} &= \biggl(\frac{7}{4}+\frac{7 \lambda^2 \Delta x^2}{120}-\frac{103\lambda^4 \Delta x^4}{6048}\biggr) g^+(u_{i}) + \biggl(-\frac{9}{2}-\frac{13 \lambda^2 \Delta x^2}{120}+\frac{1241\lambda^4 \Delta x^4}{30240}\biggr) g^+(u_{i+1}) \\
  & + \biggl(4-\frac{ \lambda^2 \Delta x^2}{40}-\frac{211\lambda^4 \Delta x^4}{10080}\biggr) g^+(u_{i+2})+ \biggl(-\frac{3}{2}+\frac{17 \lambda^2 \Delta x^2}{120}-\frac{397\lambda^4 \Delta x^4}{30240}\biggr) g^+(u_{i+3}) \\
  &+ \biggl(\frac{1}{4}-\frac{ \lambda^2 \Delta x^2}{15}+\frac{19\lambda^4 \Delta x^4}{1890}\biggr) g^+(u_{i+4}).
\end{aligned}
 \label{eqn:22}
\end{eqnarray}
By shifting each index by -1, we obtain the flux ${\hat{G}^{(m)}_{i-\frac{1}{2}}}$. Hence the Taylor series expansions of (\ref{eqn:22}) gives
\begin{equation}
\begin{aligned}
    {\hat{G}^{(0)}_{i\pm \frac{1}{2}}}&= G(x_{i\pm\frac{1}{2}})+\biggl(-\frac{1}{8}\lambda^2 \frac{d^3 h}{dx^3}\bigg|_{x=x_i}+\frac{1}{8} \frac{d^5 h}{dx^5}\bigg|_{x=x_i}\biggr)\Delta x^5+\mathcal{O}(\Delta x^6),\\  
    {\hat{G}^{(1)}_{i\pm \frac{1}{2}}}&= G(x_{i\pm\frac{1}{2}})+\biggl(\frac{1}{8}\lambda^2 \frac{d^3 h}{dx^3}\bigg|_{x=x_i}-\frac{1}{8} \frac{d^5 h}{dx^5}\bigg|_{x=x_i}\biggr)\Delta x^5+\mathcal{O}(\Delta x^6),\\     
    {\hat{G}^{(2)}_{i\pm \frac{1}{2}}}&= G(x_{i\pm\frac{1}{2}})+\biggl(-\frac{1}{8}\lambda^2 \frac{d^3 h}{dx^3}\bigg|_{x=x_i}+\frac{1}{8} \frac{d^5 h}{dx^5}\bigg|_{x=x_i}\biggr)\Delta x^5+\mathcal{O}(\Delta x^6), \\
{\hat{G}^{(m)}_{i\pm \frac{1}{2}}}&= G(x_{i\pm\frac{1}{2}})+\mathcal{B}_m\Delta x^5+\mathcal{O}(\Delta x^6),     
\end{aligned} 
 \label{eqn:23}
\end{equation}
where the values $\mathcal{B}_{m}$ are independent of $\Delta x$, but dependent on $\lambda$. Hence, we have third order approximations
\begin{eqnarray}
\begin{aligned}
    g(u)_{xxx}\big|_{x=x_i} &= \frac{{\hat{G}^{(0)}_{i+\frac{1}{2}}}-{\hat{G}^{(0)}_{i-\frac{1}{2}}}}{\Delta x^3}+\mathcal{O}(\Delta x^3),\\
    g(u)_{xxx}\big|_{x=x_i} &= \frac{{\hat{G}^{(1)}_{i+\frac{1}{2}}}-{\hat{G}^{(1)}_{i-\frac{1}{2}}}}{\Delta x^3}+\mathcal{O}(\Delta x^3),\\
    g(u)_{xxx}\big|_{x=x_i} &= \frac{{\hat{G}^{(2)}_{i+\frac{1}{2}}}-{\hat{G}^{(2)}_{i-\frac{1}{2}}}}{\Delta x^3}+\mathcal{O}(\Delta x^3).
\end{aligned}
\end{eqnarray}
%-----------------------------------------------------------------------------
\subsection{Ideal weights based on exponential polynomials}
The final WENO-E approximation is defined by a convex combination of local fluxes
with non-linear weights $\omega_m$:
\begin{equation}\label{Gnw}
    {\hat{G}_{i \pm \frac{1}{2}}} = \sum_{m=0}^2 \omega_m^{\pm} {\hat{G}^{(m)}_{i+ \frac{1}{2}}},
\end{equation}
To derive the non-linear weights $\omega_m$, we initially determine the $d_m$, referred to as ideal (or optimal) weights.  These $d_m$ values are chosen such that their linear combination with ${\hat{G}^{(m)}_{i+ \frac{1}{2}}}$ retains the fifth convergence order to ${G(x_{i+ \frac{1}{2}})}$. That is,  
\begin{equation}\label{Ghalf}
    {\hat{G}^{+}_{i+ \frac{1}{2}}} = \sum_{m=0}^2 d_m {\hat{G}^{(m)}_{i+ \frac{1}{2}}},
\end{equation}
satisfying $\sum_{m=0}^2 d_m = 1$. The ideal weights $d_m$, for the proposed WENO-E scheme can be obtained as
\begin{equation}
    d_0 = \frac{C_0}{C_{0}^0}, \quad d_1 = \frac{C_1 - d_0 C_{1}^0}{C_{0}^1}, \quad d_2 = \frac{C_2 - d_0 C_{2}^0 -d_1 C_{1}^1}{C_{0}^2}.
\end{equation}
In contrast to the classical WENO scheme, the optimal weights { $d_m$} in the proposed WENO method may exhibit variation depending on the parameter $\lambda$'s choice, but they converge towards the original ideal weights as $\Delta x \to 0$. Adding and subtracting $\sum_{m=0}^2 d_m {\hat{G}^{(m)}_{i\pm \frac{1}{2}}}$ to (\ref{Gnw}), gives:
\begin{equation}\label{eqn:Apm}
    {\hat{G}_{i \pm \frac{1}{2}}} = \sum_{m=0}^2 d_m {\hat{G}^{(m)}_{i\pm \frac{1}{2}}}+\sum_{m=0}^2 (\omega_m^{\pm}-d_m) {\hat{G}^{(m)}_{i\pm \frac{1}{2}}}=\bigl[G(x_{i\pm \frac{1}{2}})++\mathcal{A}^{\pm}\Delta x^7+\mathcal{O}(\Delta x^8)\bigr]+\sum_{m=0}^2 (\omega_m^{\pm}-d_m) {\hat{G}^{(m)}_{i\pm \frac{1}{2}}}.
\end{equation}
(The superscripts $\pm$ corresponds to the $\pm$ in the ${\hat{G}^{(m)}_{i\pm \frac{1}{2}}}$). Expanding the second term with the help of (\ref{eqn:23}) we obtain
\begin{equation}\label{eqn:Bpm}
\begin{split}
 \sum_{m=0}^2 (\omega_m^{\pm}-d_m) {\hat{G}^{(m)}_{i\pm \frac{1}{2}}} &= \sum_{m=0}^2 (\omega_m^{\pm}-d_m)\biggl( G(x_{i\pm \frac{1}{2}})+\mathcal{B}_{m}\Delta x^5+\mathcal{O}(\Delta x^6) \biggr)\\
 &= G(x_{i\pm \frac{1}{2}}) \sum_{m=0}^2 (\omega_m^{\pm}-d_m) + \Delta x^5 \sum_{m=0}^2 \mathcal{B}_{m}(\omega_m^{\pm}-d_m)+\sum_{m=0}^2 (\omega_m^{\pm}-d_m)\mathcal{O}(\Delta x^6) 
\end{split}
\end{equation}
Substituting the result above at a finite difference formula for the exponential polynomial approximation ${\hat{G}_{i \pm \frac{1}{2}}}$:
\begin{equation}
\begin{split}
 \frac{\hat{G}_{i+\frac{1}{2}}-\hat{G}_{i-\frac{1}{2}}}{\Delta x^3} &= \frac{G(x_{i+ \frac{1}{2}})-G(x_{i- \frac{1}{2}})}{\Delta x^3}+ \mathcal{O}(\Delta x^5)+ \Biggl[\frac{\sum_{m=0}^2 (\omega_m^{+}-d_m) {\hat{G}^{(m)}_{i+ \frac{1}{2}}}-\sum_{m=0}^2 (\omega_m^{-}-d_m) {\hat{G}^{(m)}_{i- \frac{1}{2}}}}{\Delta x^3}\Biggr],\\
 &= g(u)_{xxx}\big|_{x=x_i}+ \mathcal{O}(\Delta x^5)+\Biggl[\frac{G(x_{i+ \frac{1}{2}})\sum_{m=0}^2 (\omega_m^{+}-d_m) -G(x_{i- \frac{1}{2}}) \sum_{m=0}^2 (\omega_m^{-}-d_m) }{\Delta x^3}\Biggr]\\
 &+ \Delta x^2 \sum_{m=0}^2 \mathcal{B}_{m}(\omega_m^{+}-\omega_m^{-}) + \Biggl[\sum_{m=0}^2 (\omega_m^{+}-d_m) - \sum_{m=0}^2 (\omega_m^{-}-d_m) \Biggr] \mathcal{O}(\Delta x^3). 
\end{split}
\end{equation}
The $\mathcal{O}(\Delta x^5)$ term remains after division by $\Delta x^3$ because $\mathcal{A}^+ = \mathcal{A}^-$ in (\ref{eqn:Apm}). Thus, necessary and sufficient conditions for fifth-order convergence in (\ref{eqn:(7)}) are given by
\begin{equation}
\begin{split}
 \sum_{m=0}^2 (\omega_m^{\pm}-d_m) &= \mathcal{O}(\Delta x^8),\\
 \sum_{m=0}^2 \mathcal{B}_{m}(\omega_m^{+}-\omega_m^{-}) &= \mathcal{O}(\Delta x^3),\\
 \omega_m^{\pm}-d_m &= \mathcal{O}(\Delta x^2).
\end{split}
\end{equation}
Note that the first constraint is always satisfied due to normalization $\sum_{m=0}^2 \omega_m^{\pm}=\sum_{m=0}^2 d_m$  and, from (\ref{eqn:Bpm}), we see that a sufficient condition for fifth-order of convergence is given by
\begin{equation}\label{eqn:Omega_d}
\omega_m^{\pm}-d_m = \mathcal{O}(\Delta x^3).
\end{equation}
\section{Smoothness indicators and non-linear weights} \label{sec:3}
The smoothness indicator plays a pivotal role in WENO reconstruction, as it serves as a fundamental factor in determining the non-linear weights. These weights are derived by assessing the smoothness of the local solution within each sub-stencil $S_m$. In this section, we introduce a novel set of non-linear weights that enhances existing fifth-order WENO schemes. We constructs the local and global smoothness indicators employing an $L^1$-norm  approach \cite{HLY}. This global indicator gauges the approximate magnitude of the derivatives of the local solution within each sub-stencil. Moreover, it is demonstrated that within smooth regions, the nonlinear weights closely approximate the linear weights at a rate of $\mathcal{O}(\Delta x^3)$, even near the critical points. The proof presented in this section establishes the fulfilment of condition (\ref{eqn:Omega_d}), thereby ensuring that the newly proposed scheme achieves fifth-order accuracy in smooth areas.
\subsection{Development of global and local smoothness indicators}
For the construction of a smoothness indicator, we use $n^{\text{th}}$-order generalized undivided differences $\mathcal{D}_m^n g(u)$ ($n=3,4$) of $g(u)$ on the stencil $S_m$, $m=0,1,2$ is given by
\begin{equation}
    \mathcal{D}_m^n g(u_{i+\frac{1}{2}}):=\sum _{x_j \in S_m} a_{m,j}^{[n]}g(u(x_j)).
    \label{eqn:Dnm}
\end{equation}
Let $n_m$ denote the number of points inside the stencil $S_m$ and define the coefficient vector $\mathbf{a}_{m}^{[n]}:=(a_{m,j}^{[n]}:x_j \in S_m)^T$ in (\ref{eqn:Dnm}) by solving the linear system
\begin{equation*}
    \mathbf{V} \cdot \mathbf{a}^{[n]}_m = \mathbf{\delta}^{[n]},
\end{equation*}
for the non-singular matrix
\begin{equation*}
    \mathbf{V} := \biggl( \frac{(x_j - x_{i+\frac{1}{2}})^l}{\Delta x^{l}l!} : x_j \in S_m, l=0, \dots, n_m -1 \biggr), \quad \text{and} \quad \mathbf{\delta}^{[n]}:= (\delta_{n,l}: l=0, \dots, n_m -1)^T.
\end{equation*}
Note that the coefficients in (\ref{eqn:Dnm}) are independent of $\Delta x$ and evaluation point $x_{i+\frac{1}{2}}$. The operators $\mathcal{D}_m^3 g(u)$ and $\mathcal{D}_m^4 g(u)$ for $m=0,1,2,$ can be written as
\begin{equation}\label{Dm}
\begin{aligned}
    \mathcal{D}_0^3 g(u) &= -g(u_{i-1})+3g(u_{i})-3g(u_{i+1})+g(u_{i+2}),  \\
    \mathcal{D}_1^3 g(u) &= -g(u_{i-1})+3g(u_{i})-3g(u_{i+1})+g(u_{i+2}),\\
    \mathcal{D}_2^3 g(u) &= -2g(u_{i})+7g(u_{i+1})-9g(u_{i+2})+5g(u_{i+3})-g(u_{i+4}),\\
    \mathcal{D}_0^4 g(u) &= g(u_{i-2})-4g(u_{i-1})+6g(u_{i})-4g(u_{i+1})+g(u_{i+2}),\\
    \mathcal{D}_1^4 g(u) &= g(u_{i-1})-4g(u_{i})+6g(u_{i+1})-4g(u_{i+2})+g(u_{i+3}),\\
    \mathcal{D}_2^4 g(u) &= g(u_{i})-4g(u_{i+1})+6g(u_{i+2})-4g(u_{i+3})+g(u_{i+4}).
\end{aligned}
\end{equation}
Then a simple calculation with Taylor expansion shows that they can approximate the derivatives with higher accuracy than classical undivided differences.
\begin{theorem}
Let the stencil $S_m$ be a stencil around $x_{i+\frac{1}{2}}$ with $\#S_m = n_m$, and assume that $g\in \mathcal{C}^{n_m}(\Omega)$, where $\Omega$ is an open interval containing $S_m$. Then, the functional $\mathcal{D}_m^n g(u_{i+\frac{1}{2}})$ in (\ref{eqn:Dnm}) has the convergence property
\begin{equation}
    \mathcal{D}_m^n g(u_{i+\frac{1}{2}}) = \frac{d^n}{dx^n}g(u)\bigg|_{x=x_{i+\frac{1}{2}}} \Delta x^n +\mathcal{O}(\Delta x^5), \quad n=3,4. 
\end{equation}
\end{theorem}
We define the smoothness indicator $\beta_m$ in each substencil by 
\begin{equation}\label{beta}
    \beta_m = |\mathcal{D}_m^3 g(u)|+|\mathcal{D}_m^4 g(u)|, \quad m=0,1,2,
\end{equation}
and the global smoothness indicator $\zeta$ is simply defined as the absolute difference between $\beta_0$ and $\beta_2$, i.e.,
\begin{equation}\label{zeta}
    \zeta = |\beta_0 - \beta_2|.
\end{equation}
%-------------------------------------------------------------------------------------
\subsection{Construction of non-linear weights and analysis of convergence order}
The non-linear weights for the scheme are defined as \cite{ABBM}
\begin{equation}
    \omega_m = \frac{\alpha_m}{\sum_{j=0}^2 \alpha_j}, \quad \alpha_m = d_m \biggl( 1+ \frac{\zeta}{\beta_m + \Delta x^2}\biggr) , m= 0,1,2.
    \label{eqn:nonlinearw}
\end{equation}
To attain the optimal order approximation for $G(x_{i+\frac{1}{2}})$ in smooth regions, the weights should converge appropriately towards the ideal weights as $\Delta x$ approaches zero. Conversely, in regions where a discontinuity is present, the weights should effectively eliminate the contribution of stencils that contain the discontinuity.\\
Using the concept of Taylor expansion, the operators $\mathcal{D}_m^3 g(u)$ and $\mathcal{D}_m^4 g(u)$ for $m=0,1,2$ can be represented as follows
\begin{equation}
    \begin{aligned}
    \mathcal{D}_0^3 g(u) &= \Delta x^3 g^{(3)}_{i+\frac{1}{2}} + \frac{1}{8}\Delta x^5 g^{(5)}_{i+\frac{1}{2}} +\mathcal{O}(\Delta x^6),  \\
    \mathcal{D}_1^3 g(u) &= \Delta x^3 g^{(3)}_{i+\frac{1}{2}} + \frac{1}{8}\Delta x^5 g^{(5)}_{i+\frac{1}{2}} +\mathcal{O}(\Delta x^6), \\
    \mathcal{D}_2^3 g(u) &=  \Delta x^3 g^{(3)}_{i+\frac{1}{2}} - \frac{7}{8}\Delta x^5 g^{(5)}_{i+\frac{1}{2}} +\mathcal{O}(\Delta x^6),  \\
    \mathcal{D}_0^4 g(u) &= \Delta x^4 g^{(4)}_{i+\frac{1}{2}} - \frac{1}{2}\Delta x^5 g^{(5)}_{i+\frac{1}{2}} +\mathcal{O}(\Delta x^6),  \\
    \mathcal{D}_1^4 g(u) &= \Delta x^4 g^{(4)}_{i+\frac{1}{2}} + \frac{1}{2}\Delta x^5 g^{(5)}_{i+\frac{1}{2}} +\mathcal{O}(\Delta x^6),  \\
    \mathcal{D}_2^4 g(u) &= \Delta x^4 g^{(4)}_{i+\frac{1}{2}} + \frac{3}{2}\Delta x^5 g^{(5)}_{i+\frac{1}{2}} +\mathcal{O}(\Delta x^6).  
    \end{aligned}
\end{equation}
By definition of each $\beta_m$ with $m=0,1,2,$ it is straightforward that the truncation error of the local smooth indicators $\beta_m$ is of the form
\begin{equation}
\begin{aligned}
    \beta_0 &= |g^{(3)}_{i+\frac{1}{2}}| \Delta x^3 + |g^{(4)}_{i+\frac{1}{2}}| \Delta x^4 -\frac{3}{8}|g^{(5)}_{i+\frac{1}{2}}| \Delta x^5+\frac{7}{24}|g^{(6)}_{i+\frac{1}{2}}| \Delta x^6- \frac{187}{1920}|g^{(7)}_{i+\frac{1}{2}}| \Delta x^7 +\mathcal{O}(\Delta x^{8}), \\
    \beta_1 &= |g^{(3)}_{i+\frac{1}{2}}| \Delta x^3 + |g^{(4)}_{i+\frac{1}{2}}| \Delta x^4 +\frac{5}{8}|g^{(5)}_{i+\frac{1}{2}}| \Delta x^5+\frac{7}{24}|g^{(6)}_{i+\frac{1}{2}}| \Delta x^6+ \frac{71}{640}|g^{(7)}_{i+\frac{1}{2}}| \Delta x^7 +\mathcal{O}(\Delta x^{8}), \\
    \beta_2 &= |g^{(3)}_{i+\frac{1}{2}}| \Delta x^3 + |g^{(4)}_{i+\frac{1}{2}}| \Delta x^4 +\frac{5}{8}|g^{(5)}_{i+\frac{1}{2}}| \Delta x^5+\frac{7}{24}|g^{(6)}_{i+\frac{1}{2}}| \Delta x^6+ \frac{71}{640}|g^{(7)}_{i+\frac{1}{2}}| \Delta x^7 +\mathcal{O}(\Delta x^{8}). \\
    \zeta &= |g^{(5)}_{i+\frac{1}{2}}| \Delta x^5+\frac{5}{24}|g^{(7)}_{i+\frac{1}{2}}| \Delta x^7+\mathcal{O}(\Delta x^{9}). 
\end{aligned}
\end{equation}
It implies that
\begin{equation}
    \alpha_m = d_m \biggl(1+ \frac{\zeta}{\beta_m +\Delta x^2} \biggr) =  \begin{cases}
			d_m \bigl( 1+\mathcal{O}(\Delta x^3)\bigr), & \text{if} \quad g^+_{xxx}\neq 0, \\
			d_m \bigl( 1+\mathcal{O}(\Delta x^3)\bigr), & \text{if} \quad g^+_{xxx}= 0, g^+_{xxxx}\neq 0.
		\end{cases}.
\end{equation}
Then, using the fact that $\sum_{m=0}^2 d_m = 1$ and $\omega_m = \frac{\alpha_m}{\sum_{j=0}^2\alpha_j} $
\begin{equation}
    d_m = \frac{1}{1+ \frac{\zeta}{\beta_m +\Delta x^2}} \omega_m \bigl(1+\mathcal{O}(\Delta x^3)\bigr)\sum_{m=0}^2 d_m = \omega_m + \mathcal{O}(\Delta x^3).
\end{equation}
Thus, the non-linear weights satisfy the sufficient condition (\ref{eqn:Omega_d}).

\begin{theorem} 
Assume that $d_m, m=0,1,2$, is the linear weights and that the function $g(u)_{xxx}$ is smooth around the big stencil $S_7$, so the nonlinear weights (\ref{eqn:nonlinearw}) satisfy the conditions (\ref{eqn:Omega_d}) even in the presence of critical points where the third and higher derivatives vanish.
\end{theorem}
\par Now we summarise the proposed scheme in algorithmic fashion. We iteratively solve equation (\ref{eqn:1}) until the final time $t=T$. The numerical solution at the $n^{\text{th}}$ time step $t=t^n$ is denoted by $\{u^n_i:=u^n(x_i)\}$. Commencing with $n=0$ and $\{u^0_i\}$ being a given initial condition, the following steps are taken:
\begin{algorithm}\caption{}\label{alg:Alg_WenoE}
\begin{enumerate}
\item Establishing a uniform mesh distribution.
\begin{enumerate}
\item Define the spatial domain as $x \in [x_l, x_r]$.
\item  Set the number of grid points as $N$.
\item  Calculate the spatial step size as $\Delta x=\frac{(x_r-x_l)}{N}$.
\end{enumerate} 
\item  Finite difference WENO scheme for the third derivatives.
\begin{enumerate}
\item The conservative form of the equation within the MOL framework is defined, employing the numerical fluxes $\hat{F}$ and $\hat{G}$ for convection and dispersion, correspondingly.
\item   For the convection term, the flux construction procedure outlined in reference \cite{JS} is followed.
\item  The approximation of fluxes for the dispersion term is achieved by applying the WENO reconstruction approach as detailed in equation (\ref{eqn:(7)}).
 \end{enumerate}
\item  Construct exponential WENO-E reconstruction at the cell interface $x_{i+\frac{1}{2}}$.
\begin{enumerate}
\item  Split the numerical flux $\hat{G}_{i+\frac{1}{2}}$ into positive part $\hat{G}^+_{i+\frac{1}{2}}$ and negative part $\hat{G}^-_{i+\frac{1}{2}}$ as in (\ref{eqn:(8)}).
\item Construct the approximations $\hat{G}^+_{i+\frac{1}{2}}$ based on exponential polynomial (\ref{S}) on stencil $S_7$ and local approximations $\hat{G}^{(m)}_{i+\frac{1}{2}}$ based on exponential polynomial (\ref{Sm}) on each substencil $S_m$ for $m=0,1,2.$
\item   Determine the ideal weights $d_m$ and form the approximation given by (\ref{Ghalf}).
\item   Analogously, replication of steps 3(b) and 3(c) to obtain  $\hat{G}^-_{i+\frac{1}{2}}$.
 \end{enumerate} 
\item   Compute non-linear weights.
\begin{enumerate}
\item   Calculate undivided differences of $g(u)$ on the associated substencils : $\mathcal{D}^n_{S^+_m}$, $\mathcal{D}^n_{S^-_m}$ for $n=3,4$, $m=0,1,2.$
\item   Evaluate the smoothness indicators $\beta_m$ on each substencil $S_m$ as in (\ref{beta}) and determine the value $\zeta$ as in (\ref{zeta}).
\item  Compute the unnormalize weights $\alpha_m$ for each $m=0,1,2$ and using $\zeta$ and $\alpha_m$, compute the nonlinear weights $\omega_m$ as in (\ref{eqn:nonlinearw}).
 \end{enumerate} 
\item   Repeat step 3 and 4, to construct numerical flux $\hat{G}_{i-\frac{1}{2}}$ by shifting each index by -1.
\item   Time discretization.
 \begin{enumerate}
\item  Update the time step from $t^n$ to $t^{n+1}=t^n+\Delta t^n$, by applying third order SSP Runge-Kutta scheme (\ref{eqn:TVD}).
\item  If $t^{n+1}<T$, set the time step size $\Delta t^{n+1}$ by (\ref{timestep}) with updated wave speeds. In case $t^{n+1}+\Delta t^{n+1}>T$, $\Delta t^{n+1}$ is set $T-t^{n+1}$.
\item   Repeat steps beginning with step 3 for each time step until time $T$.
\end{enumerate} 
\end{enumerate} 
\end{algorithm}
%---------------------------------------------------------------------------------------
\begin{remark}
\normalfont
When $\lambda \Delta x$ is close to zero, the WENO-E scheme reduces to the WENO scheme based on algebraic polynomials as given in \cite{AQ}. In the polynomial case, the numerical flux on the big stencil $S_7$ is given by
\begin{equation}
  {\hat{G}^+_{i+\frac{1}{2}}} = 
  \biggl[\frac{-1}{15}g^+(u_{i-2}) +\frac{21}{40} g^+(u_{i-1}) +\frac{1}{8} g^+(u_{i}) -\frac{23}{12} g^+(u_{i+1}) +\frac{7}{4} g^+(u_{i+2}) -\frac{19}{40} g^+(u_{i+3}) +\frac{7}{120} g^+(u_{i+4}) \biggr].
\end{equation}
 Also the numerical fluxes on the small stencils $ S^m = \{x_{i-2+m},\dots,x_{i+2+m}\}$ with $m=0,1,2$ are given by
\begin{equation}
\begin{aligned}
 {\hat{G}^{(0)}_{i+\frac{1}{2}}} &= -\frac{1}{4} g^+(u_{i-2}) +\frac{3}{2} g^+(u_{i-1}) -2 g^+(u_{i}) +\frac{1}{2} g^+(u_{i+1}) + \frac{1}{4} g^+(u_{i+2}),\\
 {\hat{G}^{(1)}_{i+\frac{1}{2}}} &= \frac{1}{4}g^+(u_{i-1}) + \frac{1}{2} g^+(u_{i})-2 g^+(u_{i+1}) + \frac{3}{2} g^+(u_{i+2})-\frac{1}{4} g^+(u_{i+3}),\\
  {\hat{G}^{(2)}_{i+\frac{1}{2}}} &= \frac{7}{4}  g^+(u_{i}) -\frac{9}{2}  g^+(u_{i+1}) + 4 g^+(u_{i+2}) -\frac{3}{2} g^+(u_{i+3}) + \frac{1}{4} g^+(u_{i+4}).
\end{aligned}
\end{equation}
The smoothness indicator $\beta_m$ of the interpolation polynomials on the smaller stencils $S_m$ are given by
\begin{equation}
    \beta_m^Z = \sum_{\kappa =1}^{r} \Delta x^{2\kappa -1} \int_{I_j}\biggl(\frac{d^{\kappa}}{dx^{\kappa}} \Tilde{p}_m(x)\biggr)^2 dx.
\end{equation}
In order to increase the accuracy of the nonlinear weights, an alternative to mapped nonlinear weights \cite{AQ}, we chose to use the WENO-Z nonlinear weights introduced in \cite{BCCD}. The nonlinear weights $\omega_m^Z$ of WENO-Z are defined by
\begin{equation}
    \omega_m^Z = \frac{\alpha_m^Z}{\sum_{k=0}^2 \alpha_k^Z}, \quad \alpha_m^Z = d_m \biggl(1+\frac{\tau_5}{\beta_m^Z + \epsilon} \biggr), \quad \tau_5 := |\beta_0^Z - \beta_2^Z|, \quad m=0,1,2,
    \label{eqn:WENOZ}
\end{equation}
where $d_0 = \frac{4}{15}$, $d_1 = \frac{1}{2}$, $d_2 = \frac{7}{30}$ are linear weights, $\alpha_m^Z$ stands for unnormalized weights.
\end{remark}
\begin{remark}
\normalfont
    For the convection term, we use the fifth-order finite difference WENO scheme with WENO-Z nonlinear weights as (\ref{eqn:WENOZ}) and its corresponding linear weights in \cite{BCCD}.
\end{remark}
\begin{remark}
\normalfont   
In later numerical experiments, we use $\epsilon = \Delta x^2$ to ensure the fifth-order convergence for WENO-Z and WENO-E schemes, even at critical points, while keeping ENO property near discontinuity \cite{ABBM, DB}.
\end{remark}
\begin{remark}
\normalfont   
The extension of the present scheme to the two-dimensional case can be done by discretizing the spatial dimensions one at a time.
\end{remark}
\section{Numerical results}
\label{sec:4}
In this section, we illustrate several examples to test the proposed WENO-E-$\lambda \Delta x$ scheme (henceforth referred to as the WENO-E scheme) derived by the exponential approximation space for linear and non-linear dispersion-type equations with various initial conditions. We compared the performance of the proposed scheme with the WENO-Z scheme, which is based on polynomial approximations. In the WENO-E scheme, the tension parameter $\lambda$ plays a significant role. The numerical value for $\lambda$ can be tuned according to the characteristics of the initial data; which however is time-consuming and resource intense. To simplify, the $\lambda$ value is chosen based on the initial cell size such that $0 \leq \lambda \Delta x \leq 0.1$. When $\lambda \Delta x$ is close to zero, the WENO-E scheme reduces to the standard WENO-Z scheme based on algebraic polynomials. On the other hand, when $\lambda \Delta x$ is large, the exponential polynomials can better capture the rapid variations in the solution near discontinuities, resulting in improved interpolation results. However, choosing a very large value of $\lambda \Delta x$ can lead to numerical instability and accuracy issues. Thus, an optimal value of $\lambda$ is chosen depending on the given data feature and the required accuracy.  For the time integration, we use the third-order explicit strong stability preserving (SSP) Runge-Kutta method \cite{GKS}. For the ordinary differential equation of the form 
\begin{equation}
    \frac{du}{dt} =RHS(u),
\end{equation}
the SSP Runge-Kutta method is given by\\
\begin{equation}
\label{eqn:TVD}
\begin{split}
u^{(1)} &= u^n + \Delta t RHS(u^n),\\
u^{(2)} &= \frac{3}{4} u^n + \frac{1}{4} u^{(1)} + \frac{1}{4} \Delta t RHS(u^{(1)}),\\
u^{n+1} &= \frac{1}{3} u^n + \frac{2}{3} u^{(2)} + \frac{2}{3} \Delta t RHS(u^{(2)}).\\
\end{split}
\end{equation}
To ensure the $CFL$ stability condition, the time step is given by \\
\begin{equation}\label{timestep}
    \Delta t \leq \min \Biggl(\frac{CFL \cdot \Delta x^{\frac{5}{3}} }{\max(|f'(u)|)}, \frac{CFL \cdot \Delta x^{3} }{\max(|g'(u)|)} \Biggr).
\end{equation}
We have chosen $CFL = 0.3$ as obtained by Fourier analysis for the linear dispersion equation \cite{AQ} for further analysis. 
%------------------------------------------------------------------------------------
% Example 1
\begin{example}\label{example:1}
\normalfont Consider the linear KdV equations (Airy equation),
\begin{equation}
\left\{
\begin{aligned}
u_t+u_{xxx} &= 0, \quad (x,t) \in[0,2\pi] \times [0,T], \\
u_0(x)&= \sin(x), \quad x \in[0,2\pi],
\end{aligned}
\right.  
 \end{equation}
 with periodic boundary conditions and the exact solution is given by $u(x,t)=\sin(x+t)$. We solve this linear equation up to the final time $T=1$ by the WENO-E scheme and the WENO-Z scheme. The following error norms are used to compute the accuracy of schemes.
  \begin{equation*}
L^{\infty} = \max_{0\leq i \leq N} |u_e - u_{a}|,\,\,\,\,\,\,\,\,\,\, L^{1} = \dfrac{1}{N+1} \sum_{i=0}^{N} |u_e - u_{a}|,
 \end{equation*}
where $u_e$ and $u_{a}$ denote the exact and approximate solutions of the PDE.
A comprehensive analysis of error norms, specifically $L^{1}$ and $L^{\infty}$ errors, is conducted across various combinations of the parameter $\lambda \Delta x $ and grid resolution $N$. The results are detailed in Appendix 1, where we have considered 15 distinct $\lambda \Delta x $ values carefully chosen to span the entire error spectrum within the $[0,1]$ range. $\lambda \Delta x $ can be categorized into three classes based on the observed outcomes. Firstly, there exists a category where $\lambda \Delta x $ consistently yields precise and accurate results. Secondly, values in this category demonstrate reasonably accurate results for lower resolution and remains stagnant over larger values of $N$. Lastly, there is a subset of $\lambda \Delta x $ values significantly deviating from the optimal $\lambda \Delta x $, resulting in inconsistent outcomes.
\par In this example, we observe proximate convergence within the range of $0.02 \leq \lambda \Delta x  \leq 0.04$, with $\lambda \Delta x = 0.02$ exhibiting the lowest errors among these values. Additionally, we deliberately included $\lambda \Delta x $ values of $0.06$ and $0.1$ in our analysis to illustrate the substantial incongruity in results that arises when $\lambda \Delta x $ values significantly diverge from the optimal choice.
Figure \ref{eqn:F_1a} provides a comparison of errors between the WENO-Z scheme and WENO-E scheme with $\lambda \Delta x $ values of $0.02$, $0.04$, $0.06$, and $0.1$ in terms of the $L^{\infty}$- and $L^{1}$-norms, respectively. Furthermore, Table \ref{eqn:T_1} presents a summary of the $L^{\infty}$-error, $L^{1}$-error, and order of convergence for all schemes under consideration.
%\clearpage
\begin{table}[htbp!]
\centering
\captionof{table}{Comparison of WENO-Z and WENO-E schemes in terms of $L^{\infty}$- and $L^1$- errors along with their convergence rate for Example \ref{example:1} over the domain $\Omega = [0,2\pi]$ at time $T=1$.}\label{eqn:T_1}
%\bigskip\bigskip
\setlength{\tabcolsep}{0pt}
\begin{tabular*}{\textwidth}{@{\extracolsep{\fill}} l *{10}{c} }
\toprule
$\boldsymbol{N}$ &
\multicolumn{2}{c}{\textbf{WENO-Z}} &
\multicolumn{2}{c}{\textbf{WENO-E-0.02}} &
\multicolumn{2}{c}{\textbf{WENO-E-0.04}} &
\multicolumn{2}{c}{\textbf{WENO-E-0.06}} &
\multicolumn{2}{c}{\textbf{WENO-E-0.1}}\\
\cmidrule{2-3} \cmidrule{4-5} \cmidrule{6-7} \cmidrule{8-9}\cmidrule{10-11}
& $\boldsymbol{L^{\infty}}$\textbf{-error} & \textbf{Rate}
& $\boldsymbol{L^{\infty}}$\textbf{-error} & \textbf{Rate}
& $\boldsymbol{L^{\infty}}$\textbf{-error} & \textbf{Rate}
& $\boldsymbol{L^{\infty}}$\textbf{-error} & \textbf{Rate}
& $\boldsymbol{L^{\infty}}$\textbf{-error} & \textbf{Rate}\\
\midrule
10 & 2.6042e-03  & -  & 2.5610e-03 & - & 2.5611e-03 & - & 2.5611e-03 & - & 2.5602e-03 & - \\
20 & 8.704e-05  & 4.9029  & 8.7186e-05 & 4.8765  & 8.7170e-05 & 4.8768 & 8.7085e-05 & 4.8782 & 8.6335e-05 & 4.8902 \\
40 & 2.7752e-06  & 4.9711  & 2.7735e-06 & 4.9743 & 2.7627e-06 & 4.9797 & 2.7154e-06 & 5.0032 & 2.3194e-06 & 5.2181  \\ 
80 & 8.7052e-08 & 4.9946  & 8.6670e-08 & 5.0001 & 8.1182e-08 & 5.0888 & 5.7395e-08 & 5.5641 & 1.4167e-07 & 4.0331 \\ 
160 & 2.7262e-09  & 4.9969  & 2.5429e-09 & 5.0910 & 2.0500e-10 & 8.6294 & 1.2112e-08 & 2.2445 & 1.1176e-07 & 3.4211e-01 \\ 
320 & 1.1102e-10 & 4.6180  & 2.1082e-11 & 6.9143 & 1.3553e-09& -2.7250 & 7.3103e-09 & 7.2848e-01 & 5.7148e-08 & 9.6769e-01 \\
\bottomrule
$\boldsymbol{N}$ &
\multicolumn{2}{c}{\textbf{WENO-Z}} &
\multicolumn{2}{c}{\textbf{WENO-E-0.02}} &
\multicolumn{2}{c}{\textbf{WENO-E-0.04}} &
\multicolumn{2}{c}{\textbf{WENO-E-0.06}} &
\multicolumn{2}{c}{\textbf{WENO-E-0.1}}\\
\cmidrule{2-3} \cmidrule{4-5} \cmidrule{6-7} \cmidrule{8-9}\cmidrule{10-11}
& $\boldsymbol{L^{1}}$\textbf{-error} & \textbf{Rate}
& $\boldsymbol{L^1}$\textbf{-error} & \textbf{Rate}
& $\boldsymbol{L^{1}}$\textbf{-error} & \textbf{Rate}
& $\boldsymbol{L^{1}}$\textbf{-error} & \textbf{Rate}
& $\boldsymbol{L^{1}}$\textbf{-error} & \textbf{Rate}\\
\midrule
10 & 1.7456e-03 & -  & 1.7519e-03 & - & 1.7520e-03 & - & 1.7520e-03 & - & 1.7516e-03 & - \\
20 & 5.6856e-05 & 4.9403  & 5.7100e-05 & 4.9393 & 5.7089e-05 & 4.9396 & 5.7031e-05 & 4.9411 & 5.6536e-05 & 4.9534 \\
40 & 1.7815e-06 & 4.9962  & 1.7820e-06 & 5.0019 & 1.7751e-06 & 5.0072 & 1.7447e-06 & 5.0307 & 1.4902e-06 & 5.2456 \\ 
80 & 5.5652e-08 & 5.0005  & 5.5421e-08 & 5.0070 & 5.1911e-08 & 5.0957 & 3.6701e-08 & 5.5710 & 9.0592e-08 & 4.0400\\ 
160 & 1.7390e-09 & 5.0001  & 1.6221e-09 & 5.0945 & 1.3076e-10 & 8.6330 & 7.7265e-09 & 2.2479 & 7.1295e-08 & 3.4559e-01\\ 
320 & 7.0742e-11 & 4.6195 & 1.3412e-11 & 6.9182 & 8.6372e-10 & -2.7236 & 4.6589e-09 & 7.2983e-01 & 3.6420e-08 & 9.6905e-01 \\
\bottomrule
\end{tabular*}
\end{table}

\begin{figure}[ht!]
\begin{center}
\minipage{0.45\textwidth}
  \includegraphics[width=\linewidth]{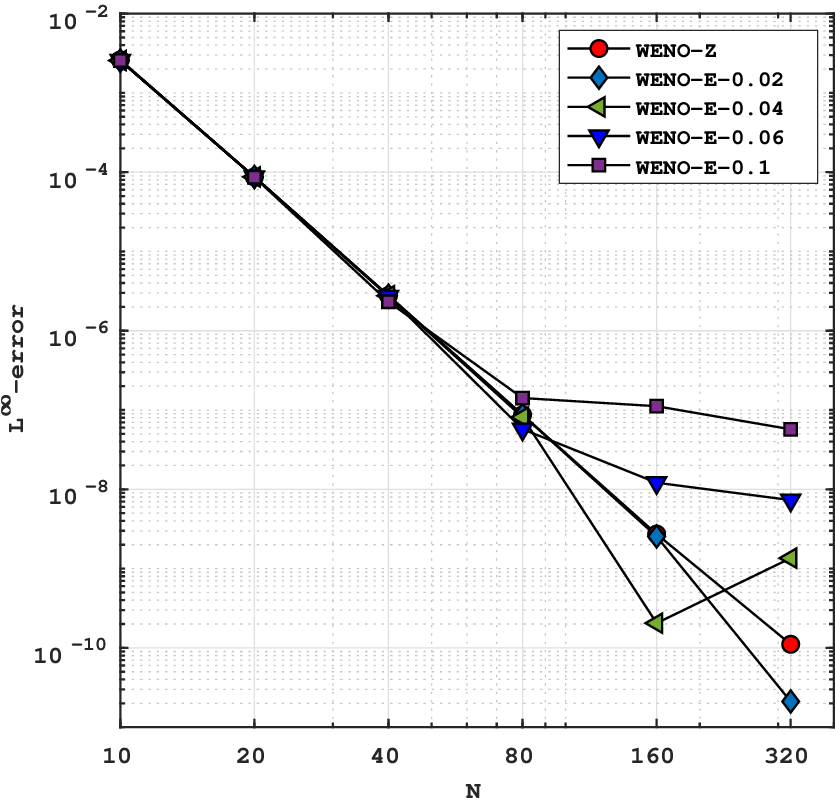}
  \subcaption*{\normalsize{\centering (a) Convergence plot: $L^{\infty}$- error}}
\endminipage\hfill
\minipage{0.45\textwidth}%
  \includegraphics[width=\linewidth]{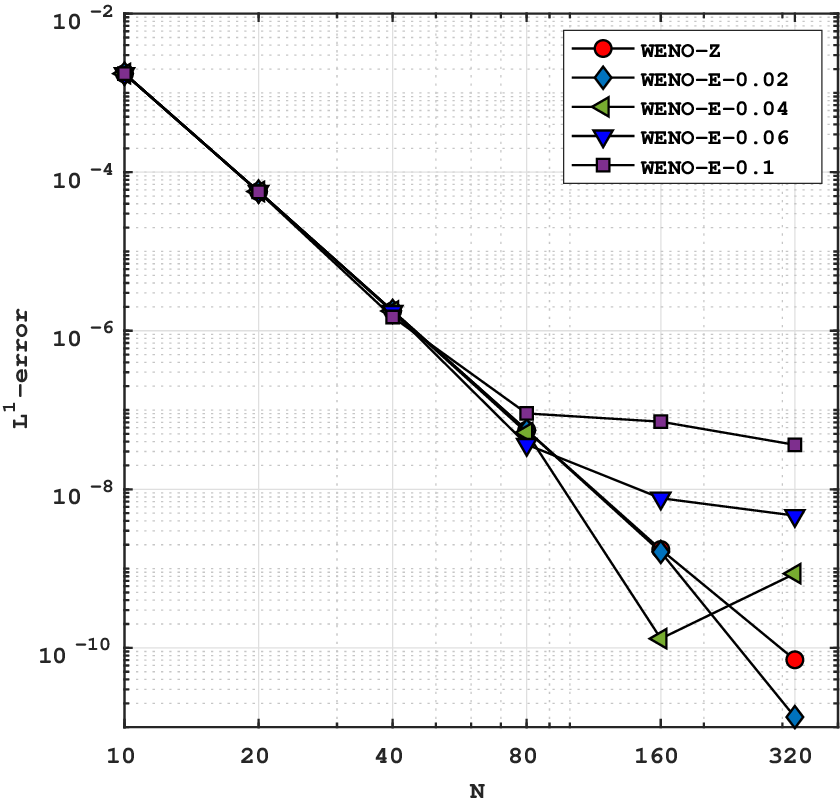}
  \subcaption*{\normalsize{\centering (b) Convergence plot: $L^{1}$- error}}
\endminipage
\caption{Comparison WENO-Z and WENO-E schemes in terms of $L^{\infty}$ and $L^{1}$ errors (in $\log_{10}$ scale) for Example \ref{example:1} at $T=1$}.\label{eqn:F_1a}
\end{center}
\end{figure}
\end{example}
%-----------------------------------------------------------------------------------
% Example 2
\begin{example}\label{example:2D}
\normalfont
We compute the solution of the two-dimensional linear  dispersion equation 
\begin{equation}
    u_t + u_{xxx} + u_{yyy} = 0, \quad (x,y,t) \in (0,2\pi)\times(0,2\pi) \times (0,T],
\end{equation}
subject to initial data $u(x,y,0)=\sin(x+y)$ with periodic boundary conditions in both directions. The exact solution is given by $u(x,y,t)=\sin(x+y+2t)$. The time step is set to a fixed value of $0.1 \Delta x^3$, and numerical solutions are computed at a specific time $T = 1$. In Figure \ref{eqn:F_2D_a} and Table \ref{eqn:T_2D}, shows a comparison of the errors of the WENO-Z scheme and WENO-E scheme with $\lambda \Delta x = 0.02, 0.04, 0.06, 0.1$ in the $L^{\infty}$- and $L^{1}$-norms, respectively. As we refined the grid,  WENO-E scheme with $\lambda \Delta x = 0.04 $ gives more accurate results and converge to the exact solution with the order of five. Figure \ref{eqn:F_2D_b} shows the numerical solutions on a grid with $80 \times 80$ grid points for WENO-Z and WENO-E scheme with $\lambda \Delta x = 0.04 $.
 \begin{table}[ht!]
\centering
\captionof{table}{Comparison of WENO-Z and WENO-E schemes in terms of $L^{\infty}$- and $L^1$- errors along with their convergence rate for Example \ref{example:2D} over the domain $\Omega = (0,2\pi)\times(0,2\pi)$ at time $T=1$.}\label{eqn:T_2D}
%\bigskip\bigskip
\setlength{\tabcolsep}{0pt}
\begin{tabular*}{\textwidth}{@{\extracolsep{\fill}} l *{10}{c} }
\toprule
$\boldsymbol{N_x \times N_y}$ &
\multicolumn{2}{c}{\textbf{WENO-Z}} &
\multicolumn{2}{c}{\textbf{WENO-E-0.02}} &
\multicolumn{2}{c}{\textbf{WENO-E-0.04}} &
\multicolumn{2}{c}{\textbf{WENO-E-0.06}} &
\multicolumn{2}{c}{\textbf{WENO-E-0.1}}\\
\cmidrule{2-3} \cmidrule{4-5} \cmidrule{6-7} \cmidrule{8-9}\cmidrule{10-11}
& $\boldsymbol{L^{\infty}}$\textbf{-error} & \textbf{Rate}
& $\boldsymbol{L^{\infty}}$\textbf{-error} & \textbf{Rate}
& $\boldsymbol{L^{\infty}}$\textbf{-error} & \textbf{Rate}
& $\boldsymbol{L^{\infty}}$\textbf{-error} & \textbf{Rate}
& $\boldsymbol{L^{\infty}}$\textbf{-error} & \textbf{Rate}\\
\midrule
10 $\times$ 10 & 5.4139e-03  & - & 5.3032e-03 & - & 5.3035e-03 & - & 5.3035e-03  & -  & 5.3020e-03
  & - \\
20 $\times$ 20 & 1.7482e-04 & 4.9527 & 1.7505e-04 & 4.9210 & 1.7502e-04 & 4.9214 & 1.7484e-04 & 4.9228 & 1.7333e-04 & 4.9349 \\ 
40 $\times$ 40 & 5.5487e-06 & 4.9776 & 5.5462e-06 & 4.9801 & 5.5245e-06 & 4.9855 & 5.4300e-06 & 5.0090 &  4.6380e-06  & 5.2239 \\ 
80 $\times$ 80 & 1.7411e-07 & 4.9940 & 1.7336e-07 & 4.9997 & 1.6238e-07 & 5.0884 & 1.1480e-07 & 5.5637  & 2.8337e-07 & 4.0327 \\ 
160 $\times$ 160 & 5.4551e-09 & 4.9963 & 5.0885e-09 & 5.0904 & 4.0691e-10 & 8.6405 & 2.4219e-08 & 2.2450  & 2.2350e-07 & 0.3425 \\ 

\bottomrule
$\boldsymbol{N_x \times N_y}$ &
\multicolumn{2}{c}{\textbf{WENO-Z}} &
\multicolumn{2}{c}{\textbf{WENO-E-0.02}} &
\multicolumn{2}{c}{\textbf{WENO-E-0.04}} &
\multicolumn{2}{c}{\textbf{WENO-E-0.06}} &
\multicolumn{2}{c}{\textbf{WENO-E-0.1}}\\
\cmidrule{2-3} \cmidrule{4-5} \cmidrule{6-7} \cmidrule{8-9}\cmidrule{10-11}
& $\boldsymbol{L^{1}}$\textbf{-error} & \textbf{Rate}
& $\boldsymbol{L^{1}}$\textbf{-error} & \textbf{Rate}
& $\boldsymbol{L^{1}}$\textbf{-error} & \textbf{Rate}
& $\boldsymbol{L^{1}}$\textbf{-error} & \textbf{Rate}
& $\boldsymbol{L^{1}}$\textbf{-error} & \textbf{Rate}\\
\midrule
10 $\times$ 10 & 3.4560e-03 & - & 3.4272e-03 & - & 3.4273e-03 & - & 3.4273e-03  & -  & 3.4261e-03 & - \\
20 $\times$ 20 & 1.1104e-04 & 4.9599 & 1.1149e-04 & 4.9420 & 1.1147e-04 & 4.9423 & 1.1136e-04 & 4.9438 & 1.1039e-04  & 4.9559  \\ 
40 $\times$ 40 & 3.5318e-06 & 4.9746 & 3.5329e-06 & 4.9800 & 3.5191e-06 & 4.9853 & 3.4588e-06 & 5.0088 & 2.9543e-06  & 5.2237   \\ 
80 $\times$ 80 & 1.1080e-07 & 4.9944 & 1.1033e-07 & 5.0009 & 1.0335e-07 & 5.0896 & 7.3069e-08 & 5.5649 & 1.8036e-07 & 4.0339 \\
160 $\times$ 160 & 3.4732e-09 & 4.9955 & 3.2400e-09 & 5.0898 & 2.5903e-10 & 8.6402 & 1.5421e-08 & 2.2443 & 1.4231e-07 & 0.3418 \\  
\bottomrule
\end{tabular*}
\end{table}
\begin{figure}[ht!]
\begin{center}
\minipage{0.45\textwidth}
  \includegraphics[width=\linewidth]{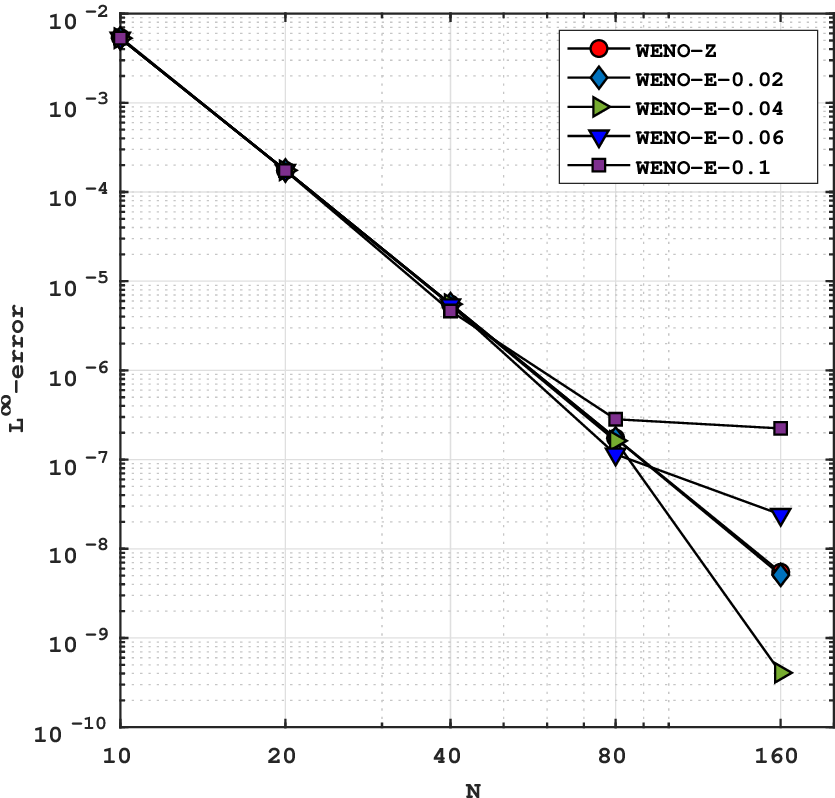}
    \subcaption*{\normalsize{\centering (a) Convergence plot: $L^{\infty}$- error}}
\endminipage\hfill
\minipage{0.45\textwidth}%
  \includegraphics[width=\linewidth]{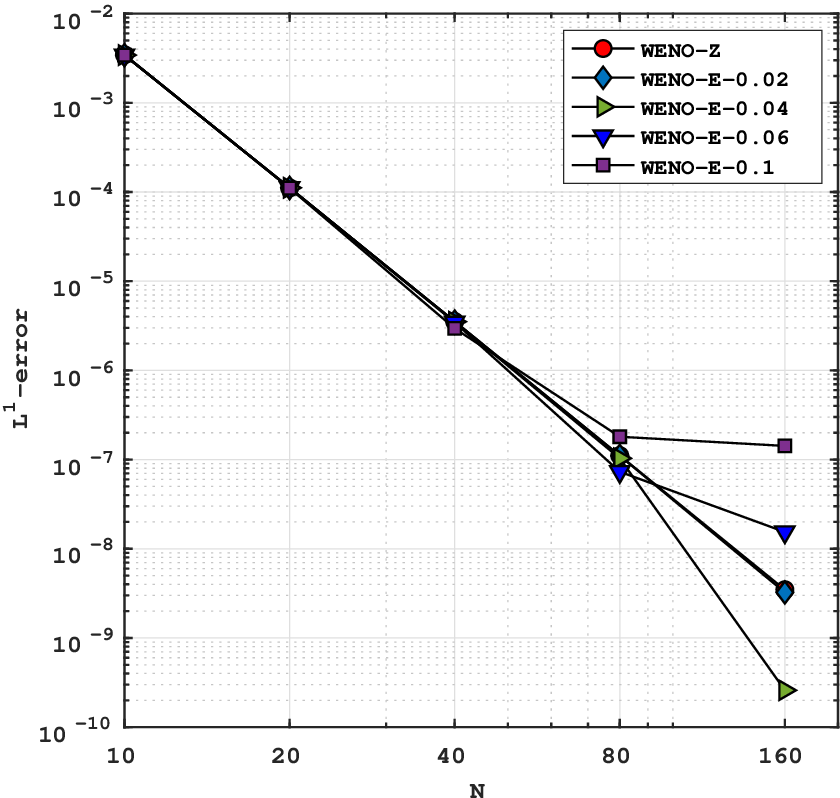}
    \subcaption*{\normalsize{\centering (b) Convergence plot: $L^{1}$- error}}
\endminipage
\caption{Comparison WENO-Z and WENO-E schemes in terms of $L^1$ and $L^{\infty}$ errors (in $\log_{10}$ scale) for Example \ref{example:2D} at $T=1$.}\label{eqn:F_2D_a}
\end{center}
\end{figure}

\begin{figure}[ht!]
  \begin{minipage}[b]{0.5\linewidth}
    \centering
    \includegraphics[trim=0.1cm 0cm 0cm -0.6cm, clip=true,width=\linewidth]{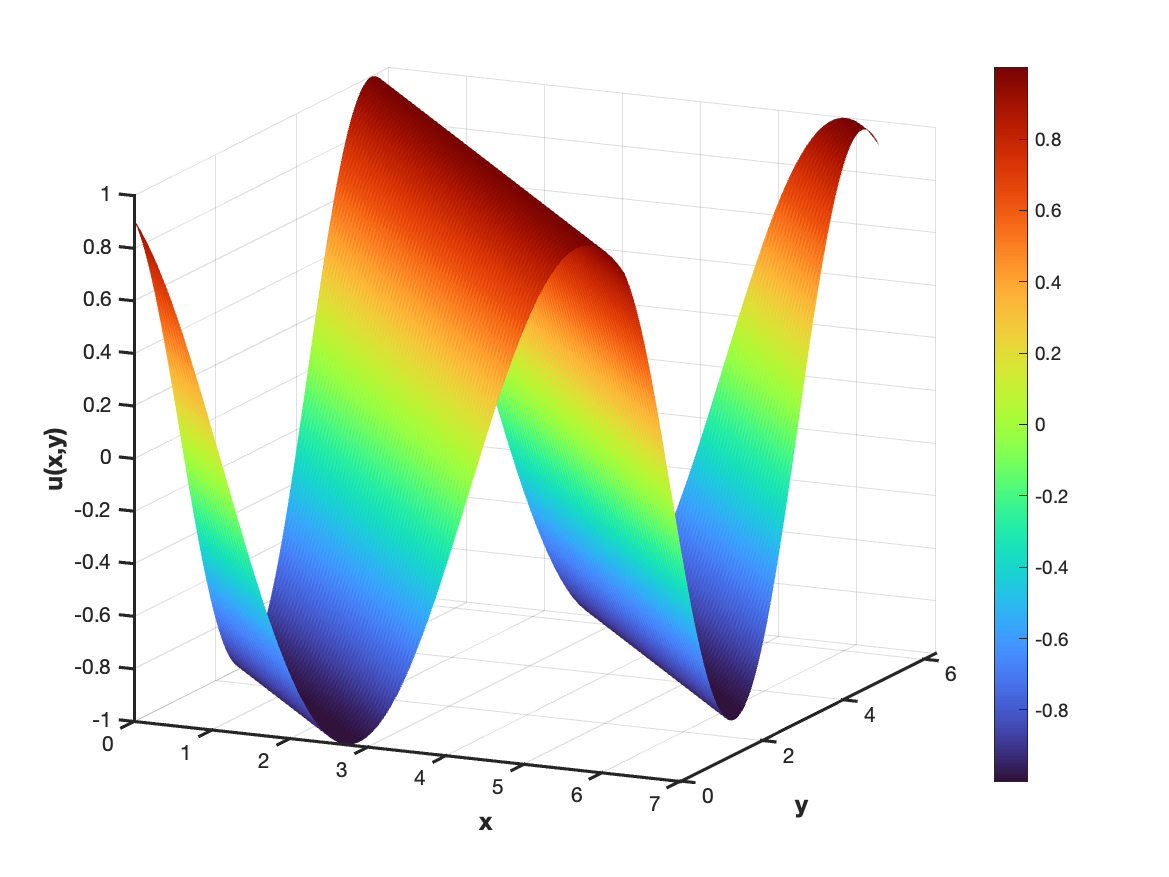}
    \subcaption*{\normalsize{\centering (a) Numerical solution: WENO-Z}}
  \end{minipage}\hfill
  \begin{minipage}[b]{0.5\linewidth}
    \centering
   \includegraphics[trim=0.1cm 0cm 0cm -0.6cm, clip=true,width=\linewidth]{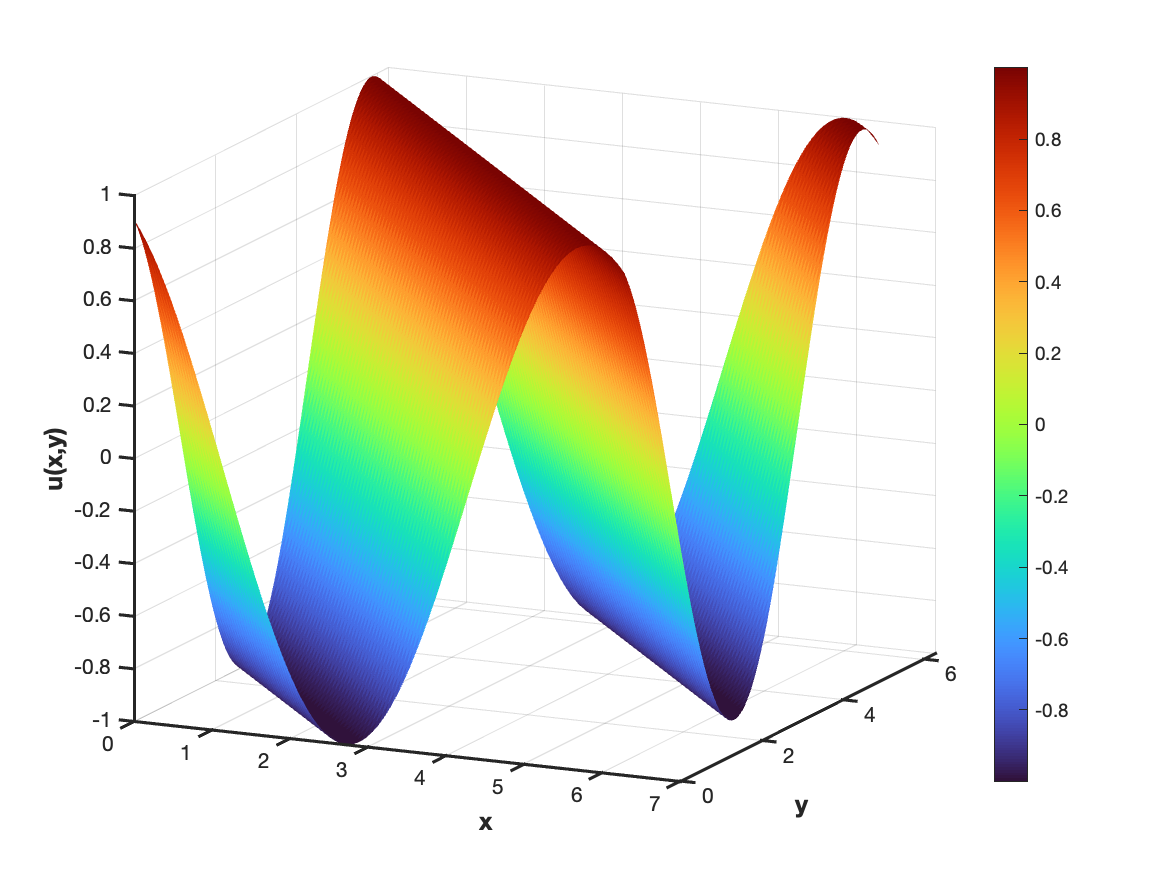}
       \subcaption*{\normalsize{\centering (b) Numerical solution: WENO-E ($\lambda \Delta x =0.04$)}}
  \end{minipage} 
  \caption{ Comparison of numerical solutions obtained from using WENO-Z and WENO-E schemes with $N_x \times N_y=80 \times 80$ at $T=1$ in $(0,2\pi) \times (0,2\pi)$, for Example, \ref{example:2D}.}\label{eqn:F_2D_b}
\end{figure}
\end{example}

 %-------------------------------------------------------------------------------
 % Example 3
\begin{example}\label{example:2}
\normalfont
Consider the nonlinear KdV equation 
\begin{equation}
\left\{
\begin{aligned}
u_t - 3(u^2)_x+ u_{xxx}&=0, \quad x \in[-10,10], \quad t \ge 0, \\
u_0(x)&= -2\sech^2(x), \quad x \in[-10,10].
\end{aligned}
\right.  
 \end{equation}
 The exact solution is given by $u(x,t)=-2\sech^2(x-4t)$.  We embark on the numerical computation of the solution within the spatial domain $x \in[-10,10]$ at a specific time $T=0.5$ while employing periodic boundary conditions. The $L^{\infty}$-error, $L^{1}$-error and convergence order are reported in Table \ref{eqn:T_2} and Figure \ref{eqn:F_2a}. A comparison of numerical solutions obtained from WENO-Z and WENO-E schemes with $N=640$ at $T=0.5$ in $x \in [-10,10]$ is given in Figure \ref{eqn:F_2b}. The result of WENO-E with $\lambda \Delta x = 0.04$ is close to the exact solution and better compared to the polynomial case. As the number of grid points increases, the error converges to machine epsilon.  
\begin{table}[ht!]
\centering
\captionof{table} {Comparison of WENO-Z and WENO-E schemes in terms of $L^{\infty}$- and $L^1$- errors along with their convergence rate for Example \ref{example:2} over the domain $\Omega =[-10,10]$ at time $T=0.5$.}\label{eqn:T_2}
%\bigskip\bigskip
\setlength{\tabcolsep}{0pt}
\begin{tabular*}{\textwidth}{@{\extracolsep{\fill}} l *{10}{c} }
\toprule
$\boldsymbol{N}$  &
\multicolumn{2}{c}{\textbf{WENO-Z}} &
\multicolumn{2}{c}{\textbf{WENO-E-0.02}} &
\multicolumn{2}{c}{\textbf{WENO-E- 0.04}} &
\multicolumn{2}{c}{\textbf{WENO-E-0.06}} &
\multicolumn{2}{c}{\textbf{WENO-E-0.1}}\\
\cmidrule{2-3} \cmidrule{4-5} \cmidrule{6-7} \cmidrule{8-9} \cmidrule{10-11}
& $\boldsymbol{L^{\infty}}$\textbf{-error} & \textbf{Rate}
& $\boldsymbol{L^{\infty}}$\textbf{-error} & \textbf{Rate}
& $\boldsymbol{L^{\infty}}$\textbf{-error} & \textbf{Rate}
& $\boldsymbol{L^{\infty}}$\textbf{-error} & \textbf{Rate}
& $\boldsymbol{L^{\infty}}$\textbf{-error} & \textbf{Rate}\\
\midrule
80 & 8.3944e-02 & - & 1.2351e-02 & - & 1.2351e-02 & - & 1.2352e-02 & -  & 1.2354e-02 & - \\ 
160 & 3.6399e-04 & 7.8494 & 4.1591e-04 & 4.8922 & 4.1591e-04 & 4.8922 & 4.1583e-04 & 4.8926 & 4.1486e-04 & 4.8962  \\
320 & 1.2547e-05 & 4.8586 & 1.3086e-05 & 4.9901 & 1.3072e-05 & 4.9918 & 1.3004e-05 & 4.9990 & 1.2429e-05 & 5.0608 \\ 
640 & 8.5693e-07 & 3.8720 &  8.5676e-07 & 3.9330 & 8.5721e-07 & 3.9306 & 8.5918e-07 & 3.9199 & 8.7566e-07 & 3.8272  \\ 
1280 & 8.9345e-07 & -0.0602 &  8.9347e-07 & -0.0605 & 8.9369e-07 & -0.0601 & 8.9467e-07 & -0.0584 & 9.0287e-07 & -0.0441  \\ 
\bottomrule
$\boldsymbol{N}$ &
\multicolumn{2}{c}{\textbf{WENO-Z}} &
\multicolumn{2}{c}{\textbf{WENO-E-0.02}} &
\multicolumn{2}{c}{\textbf{WENO-E-0.04}} &
\multicolumn{2}{c}{\textbf{WENO-E-0.06}} &
\multicolumn{2}{c}{\textbf{WENO-E-0.1}}\\
\cmidrule{2-3} \cmidrule{4-5} \cmidrule{6-7} \cmidrule{8-9} \cmidrule{10-11}
& $\boldsymbol{L^{1}}$\textbf{-error} & \textbf{Rate}
& $\boldsymbol{L^{1}}$\textbf{-error} & \textbf{Rate}
& $\boldsymbol{L^{1}}$\textbf{-error} & \textbf{Rate}
& $\boldsymbol{L^{1}}$\textbf{-error} & \textbf{Rate}
& $\boldsymbol{L^{1}}$\textbf{-error} & \textbf{Rate}\\
\midrule
80 & 2.3854e-02 & - & 2.1970e-03 & - & 2.1972e-03 & - & 2.1974e-03 & - & 2.1979e-03 & - \\ 
160 & 5.9718e-05 & 8.6419 & 6.8676e-05 & 4.9996 & 6.8676e-05 & 4.9997 & 6.8663e-05 & 5.0001 & 6.8506e-05 & 5.0038  \\
320 & 1.9878e-06 & 4.9089 & 2.0648e-06 & 5.0557 & 2.0624e-06 & 5.0574 & 2.0510e-06 & 5.0651 & 1.9547e-06 & 5.1312 \\ 
640 & 8.4259e-08 & 4.5602 & 8.4750e-08 & 4.6067 & 8.3420e-08 & 4.6278 & 7.7660e-08 & 4.7230 & 4.4961e-08 & 5.4421 \\ 
1280 & 2.5115e-08 & 1.7463 & 2.5080e-08 & 1.7567 & 2.4566e-08 & 1.7637 & 2.5510e-08 & 1.6061 & 4.8688e-08 & -0.1149 \\ 
\bottomrule
\end{tabular*}
\end{table}

\begin{figure}[htb!]
\begin{center}
\minipage{0.45\textwidth}
  \includegraphics[width=\linewidth]{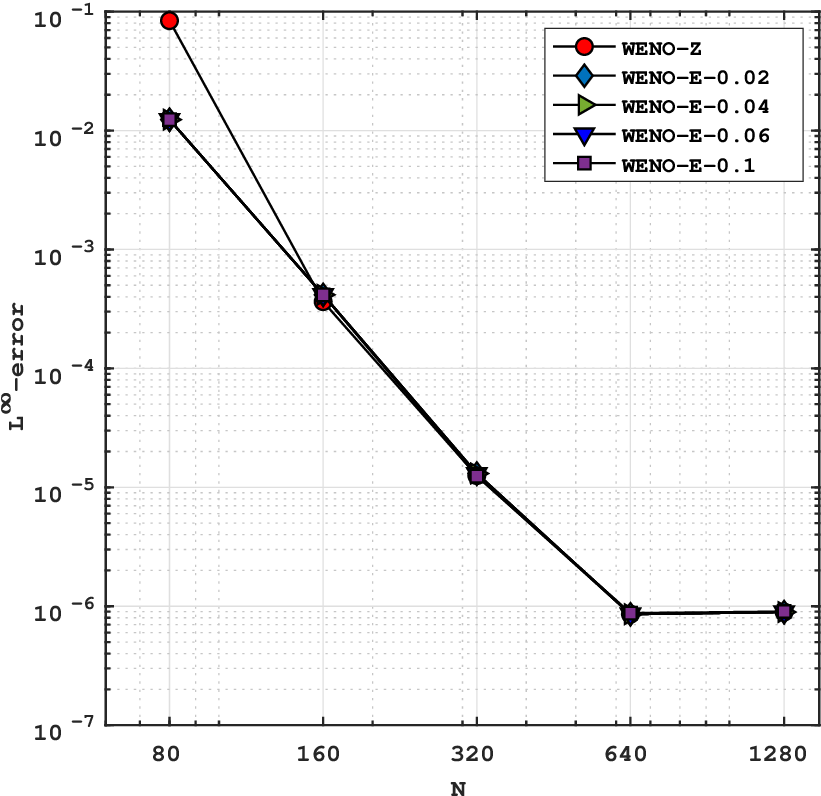}
      \subcaption*{\normalsize{\centering (a) Convergence plot: $L^{\infty}$- error}}
\endminipage\hfill
\minipage{0.45\textwidth}%
  \includegraphics[width=\linewidth]{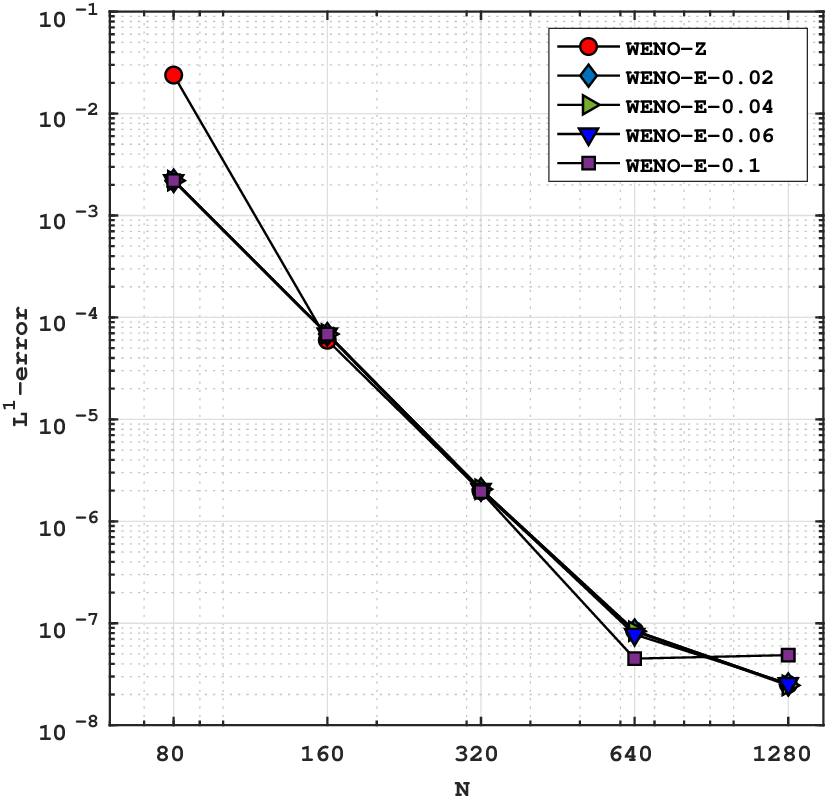}
      \subcaption*{\normalsize{\centering (b) Convergence plot: $L^{1}$- error}}
\endminipage
\caption{Comparison WENO-Z and WENO-E schemes in terms of $L^1$ and $L^{\infty}$ errors (in $\log_{10}$ scale) for Example \ref{example:2} at $T=0.5$.}\label{eqn:F_2a}
\end{center}
\end{figure}
\begin{figure}[htb!]
\captionsetup[subfigure]{justification=centering}
    \centering
      \begin{subfigure}{0.41\textwidth}
        \includegraphics[width=\textwidth]{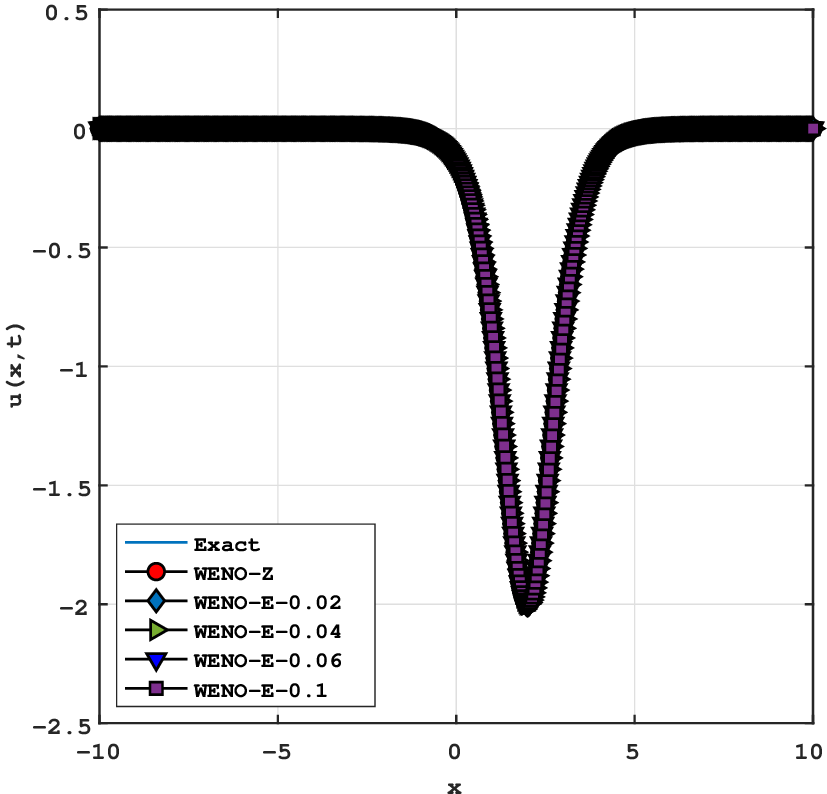}
          \caption*{\normalsize{(a) Numerical solution}}
      \end{subfigure}\hfill
      \begin{subfigure}{0.49\textwidth}
        \includegraphics[trim=0cm 0cm 0cm 0cm, clip=true,width=\textwidth]{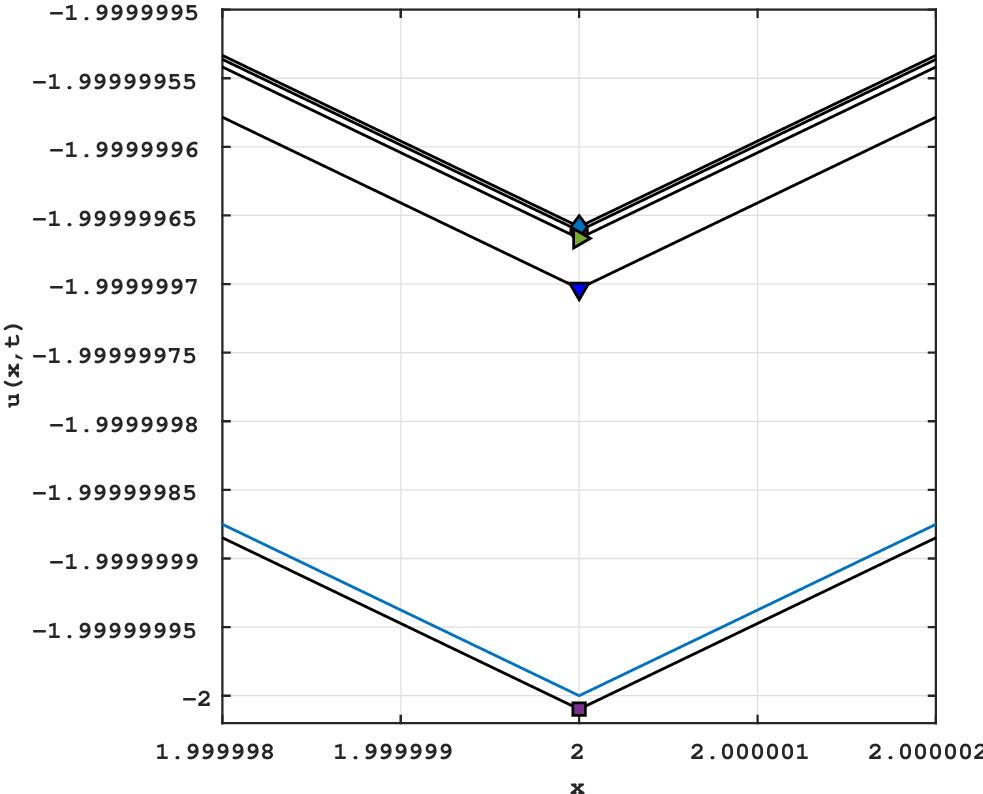}
          \caption*{\normalsize{(b) Zoom in near smooth extrema}}
      \end{subfigure}
\caption{ Comparison of numerical solutions obtained from using WENO-Z and WENO-E schemes with $N=640$ at $T=0.5$ in $x \in [-10,10]$, for Example \ref{example:2}.}\label{eqn:F_2b}
\end{figure}
\end{example}
 \begin{table}[ht!]
\centering
\captionof{table} {Comparison of WENO-Z and WENO-E schemes in terms of $L^{\infty}$- and $L^1$- errors along with their convergence rate for Example \ref{example:4a} over the domain $\Omega =[0,2\pi]$ at time $T=\pi/2$.} \label{eqn:T_4a}
%\bigskip\bigskip
\setlength{\tabcolsep}{0pt}
\begin{tabular*}{\textwidth}{@{\extracolsep{\fill}} l *{10}{c} }
\toprule
$\boldsymbol{N}$ &
\multicolumn{2}{c}{\textbf{WENO-Z}} &
\multicolumn{2}{c}{\textbf{WENO-E-0.02}} &
\multicolumn{2}{c}{\textbf{WENO-E-0.04}} &
\multicolumn{2}{c}{\textbf{WENO-E-0.06}} &
\multicolumn{2}{c}{\textbf{WENO-E-0.1}}\\
\cmidrule{2-3} \cmidrule{4-5} \cmidrule{6-7} \cmidrule{8-9} \cmidrule{10-11}
& $\boldsymbol{L^{\infty}}$\textbf{-error} & \textbf{Rate}
& $\boldsymbol{L^{\infty}}$\textbf{-error} & \textbf{Rate}
& $\boldsymbol{L^{\infty}}$\textbf{-error} & \textbf{Rate}
& $\boldsymbol{L^{\infty}}$\textbf{-error} & \textbf{Rate}
& $\boldsymbol{L^{\infty}}$\textbf{-error} & \textbf{Rate}\\
\midrule
40 & 6.3724e-04 & - & 6.3724e-04 & - & 6.3724e-04 & - & 6.3718e-04  & -  & 6.3673e-04 & - \\
80 & 1.7932e-05 & 5.1513 & 1.7927e-05 & 5.1516 & 1.7920e-05 & 5.1522 & 1.7890e-05 & 5.1544  & 1.7642e-05 & 5.1736\\ 
160 & 7.0911e-07 & 4.6604 & 7.0899e-07 & 4.6602 & 7.0660e-07 & 4.6646 & 6.9624e-07 & 4.6835  & 6.6316e-07 & 4.7336 \\ 
320 & 4.6784e-08 & 3.9219 & 4.6796e-08 & 3.9213 & 4.6983e-08 & 3.9107 & 4.7792e-08 & 3.8647  & 6.4315e-08 & 3.3661\\
640 & 5.0636e-09 & 3.2078 & 5.0663e-09 & 3.2074 & 5.1073e-09 & 3.2015 & 5.3168e-09 & 3.1682  & 3.4122e-08 & 0.9144\\
1280 & 4.1144e-10 & 3.6214 & 4.0924e-10 & 3.6299 & 5.0185e-10 & 3.3472 & 2.3221e-09 & 1.1951  & 1.7924e-08 & 0.9288\\

\bottomrule
$\boldsymbol{N}$ &
\multicolumn{2}{c}{\textbf{WENO-Z}} &
\multicolumn{2}{c}{\textbf{WENO-E-0.02}} &
\multicolumn{2}{c}{\textbf{WENO-E-0.04}} &
\multicolumn{2}{c}{\textbf{WENO-E-0.06}} &
\multicolumn{2}{c}{\textbf{WENO-E-0.1}}\\
\cmidrule{2-3} \cmidrule{4-5} \cmidrule{6-7} \cmidrule{8-9} \cmidrule{10-11}
& $\boldsymbol{L^{1}}$\textbf{-error} & \textbf{Rate}
& $\boldsymbol{L^{1}}$\textbf{-error} & \textbf{Rate}
& $\boldsymbol{L^{1}}$\textbf{-error} & \textbf{Rate}
& $\boldsymbol{L^{1}}$\textbf{-error} & \textbf{Rate}
& $\boldsymbol{L^{1}}$\textbf{-error} & \textbf{Rate}\\
\midrule
40 & 3.1608e-04 & - & 3.1654e-04 & - & 3.1654e-04 & - & 3.1644e-04 & -  & 3.1628e-04 & - \\
80 & 1.1146e-05 & 4.8256 & 1.1149e-05 & 4.8273 & 1.1146e-05 & 4.8278 & 1.1091e-05 & 4.8345  & 1.1005e-05 & 4.8449\\ 
160 & 4.5757e-07 & 4.6064 & 4.5746e-07 & 4.6072 & 4.5581e-07 & 4.6120 & 4.2927e-07 & 4.6914  & 3.8846e-07 & 4.8243 \\ 
320 & 2.0152e-08 & 4.5050 & 2.0104e-08 & 4.5081 & 1.9391e-08 & 4.5550 & 1.3932e-08 & 4.9454  & 2.5432e-08 & 3.9330\\ 
640 & 1.4181e-09 & 3.8290 & 1.4092e-09 & 3.8346 & 1.3336e-09 & 3.8620 & 2.1872e-09 & 2.9288  & 1.7251e-08 & 0.5600\\
1280 & 1.5958e-10 & 3.1516 & 1.5226e-10 & 3.2102 & 2.1720e-10 & 2.6182 & 1.0973e-09 & 0.9951  & 8.7921e-09 & 0.9724\\
\bottomrule
\end{tabular*}
\end{table}
\begin{example} \label{example:3}
\normalfont
Consider the following classical KdV equation with zero dispersion limit of conservation law
\begin{equation}
\begin{split}
u_t + \left(\frac{u^2}{2}\right)_x+\epsilon u_{xxx}&=0, \quad x \in[0,1], \quad t \ge 0, \\
\end{split} 
 \end{equation}
 with continuous initial condition
 \begin{equation}
u(x,0) = 2+0.5 \sin(2\pi x), \quad x \in[0,1],
\label{eqn:3a}
\end{equation}
 and discontinuous initial condition
  \begin{equation}
     u(x,0)=\begin{cases}
			1, & \text{if} \quad 0.25<x<4, \\
			0, & \text{else}.
		\end{cases}
  \label{eqn:3b}
 \end{equation}
 
\begin{figure}[htb!]
\captionsetup[subfigure]{justification=centering}
    \centering
      \begin{subfigure}{0.32\textwidth}
        \includegraphics[width=\textwidth]{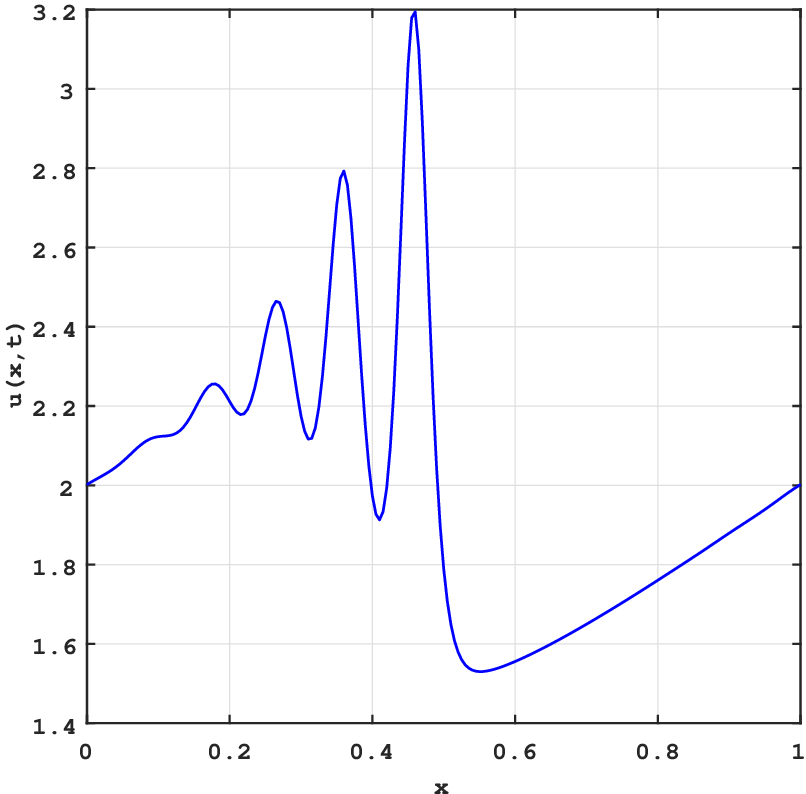}
          \caption*{\normalsize{(a) $T=0.5$, $N=200$, $\epsilon=10^{-4}$}}
      \end{subfigure}\hfill
      \begin{subfigure}{0.32\textwidth}
        \includegraphics[width=\textwidth]{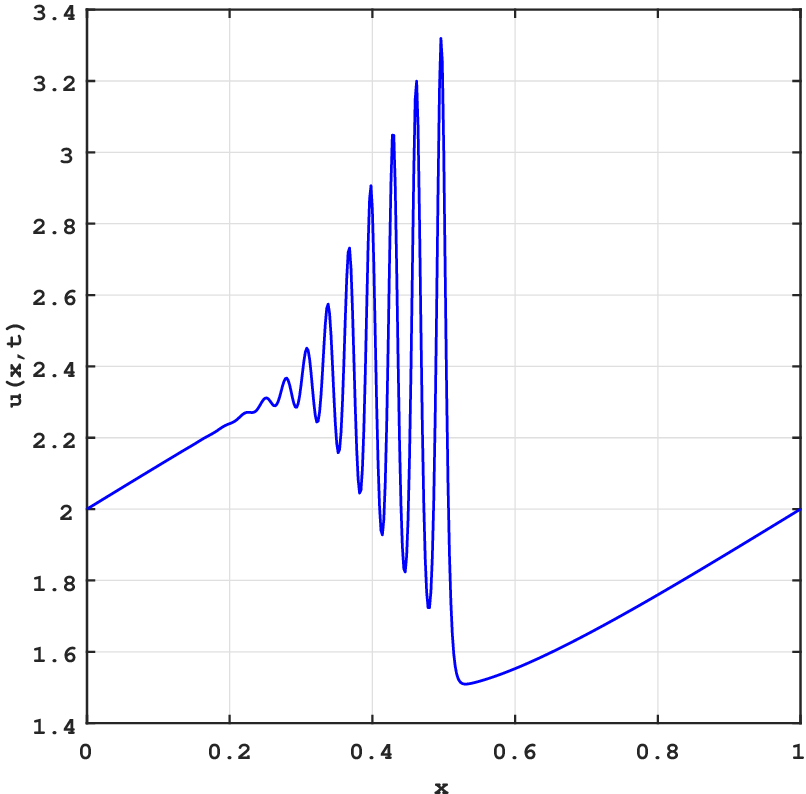}
          \caption*{\normalsize{(b) $T=0.5$, $N=500$, $\epsilon=10^{-5}$}}
      \end{subfigure}\hfill
      \begin{subfigure}{0.32\textwidth}
        \includegraphics[width=\textwidth]{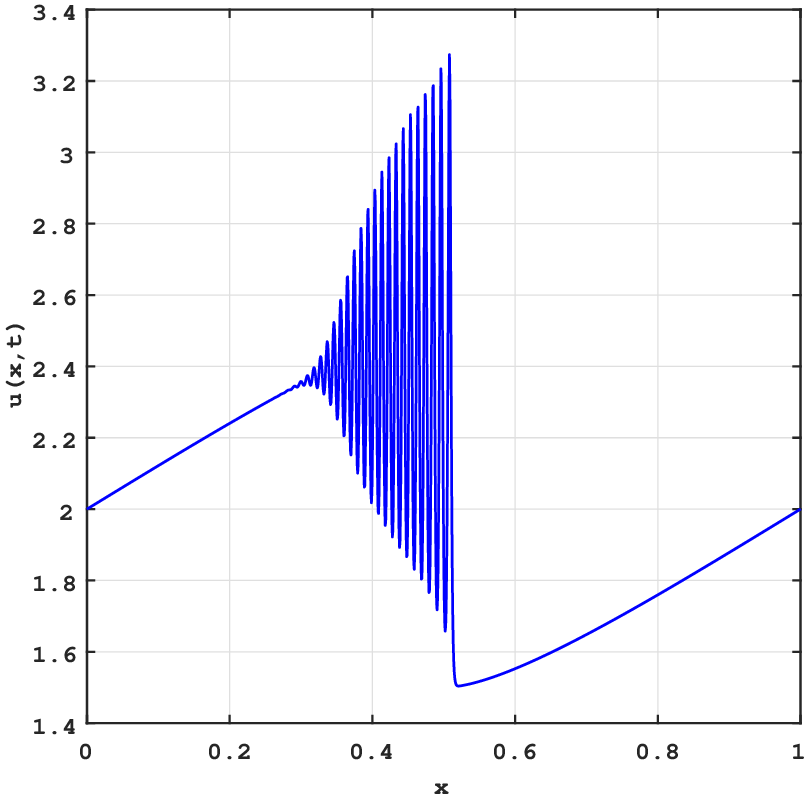}
          \caption*{\normalsize{(c) $T=0.5$, $N=1500$, $\epsilon=10^{-6}$}}
      \end{subfigure}
      \caption{Numerical solution of KdV equation with $N=200,500,1000$ at $T=0.5$ and $\lambda \Delta x = 0.02$ in $x \in [0,1]$, WENO-E for initial condition (\ref{eqn:3a}) of Example \ref{example:3}.}\label{eqn:F_3a}
\end{figure}

\begin{figure}[htb!]
\captionsetup[subfigure]{justification=centering}
    \centering
      \begin{subfigure}{0.32\textwidth}
        \includegraphics[width=\textwidth]{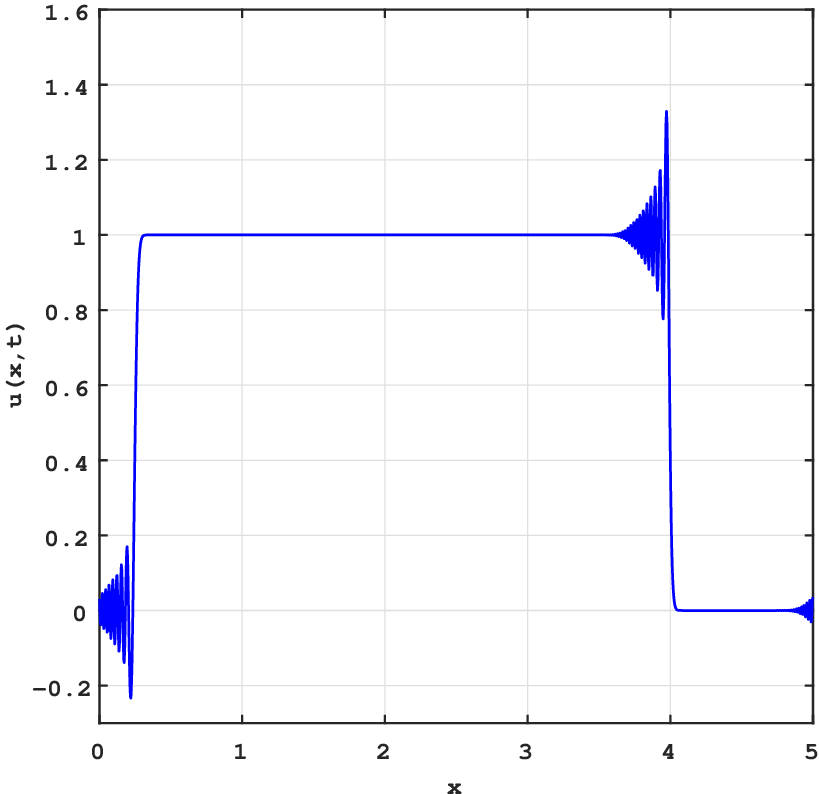}
          \caption*{\normalsize{(a) $T=0.01$, $N=1500$, $\epsilon=10^{-4}$}}
      \end{subfigure}\hfill
      \begin{subfigure}{0.32\textwidth}
        \includegraphics[width=\textwidth]{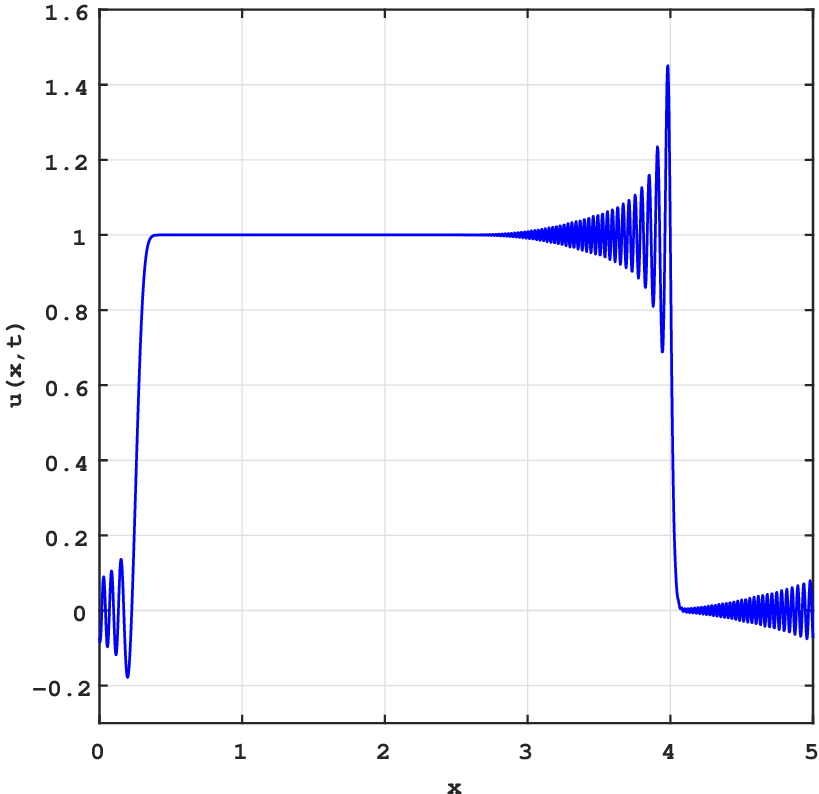}
          \caption*{\normalsize{(b) $T=0.05$, $N=1500$, $\epsilon=10^{-4}$}}
      \end{subfigure}\hfill
      \begin{subfigure}{0.32\textwidth}
        \includegraphics[width=\textwidth]{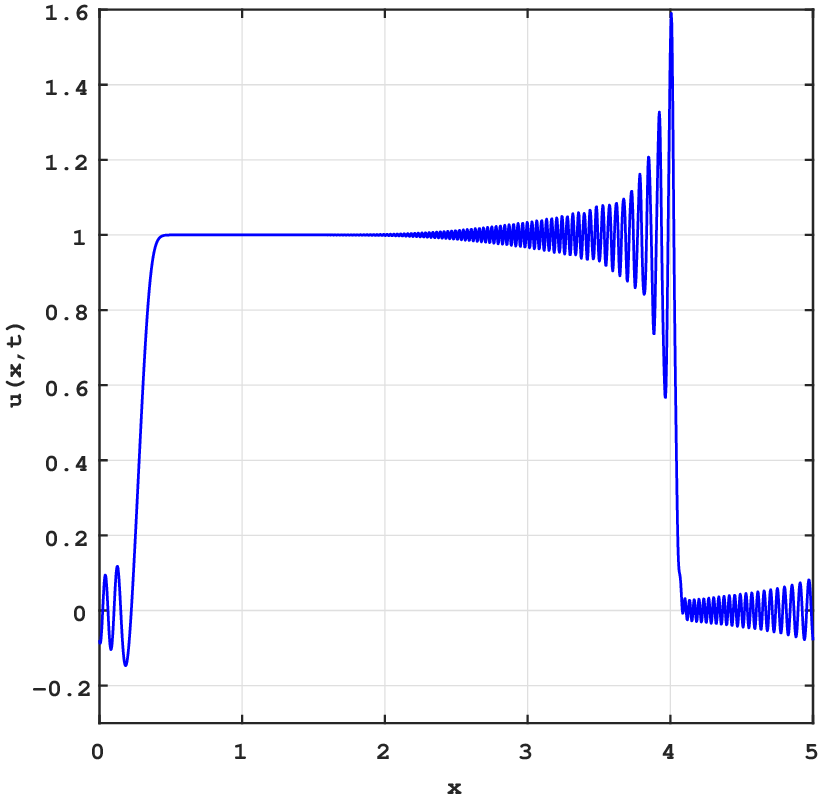}
          \caption*{\normalsize{(c) $T=0.1$, $N=1500$, $\epsilon=10^{-4}$}}
      \end{subfigure}\hfill
            \begin{subfigure}{0.32\textwidth}
        \includegraphics[width=\textwidth]{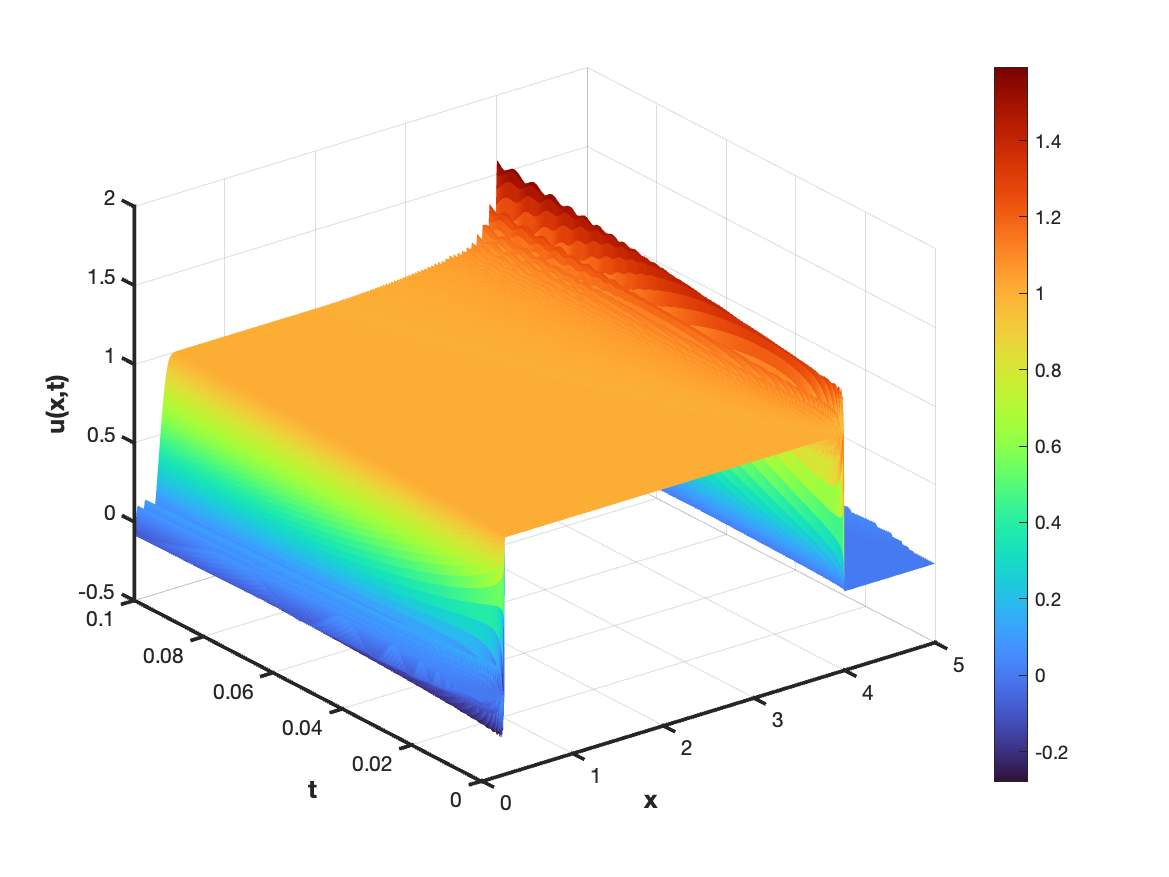}
          \caption*{\normalsize{(d) Wave simulation across various T}}
      \end{subfigure}\hfill
      \begin{subfigure}{0.32\textwidth}
        \includegraphics[width=\textwidth]{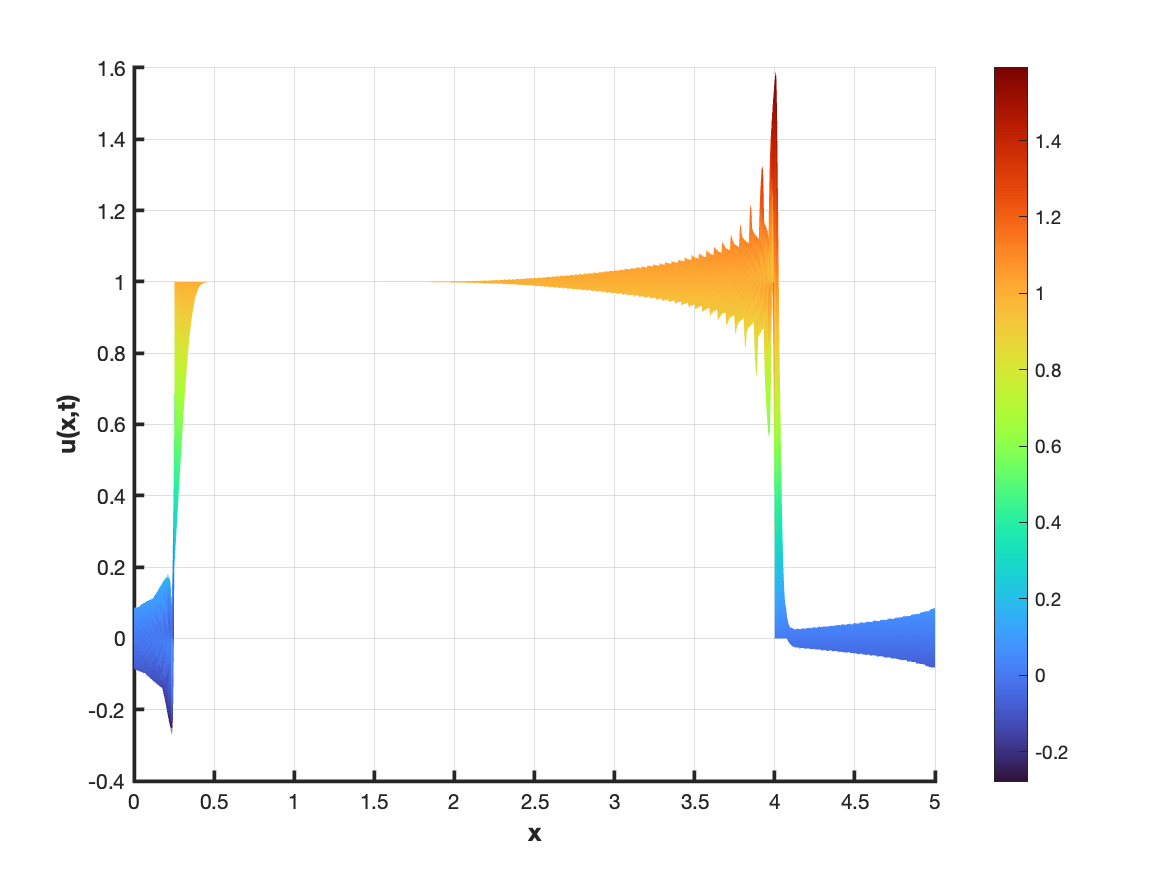}
          \caption*{\normalsize{(e) Wave simulation in 2D}}
      \end{subfigure}\hfill
           \begin{subfigure}{0.32\textwidth}
        \includegraphics[width=\textwidth]{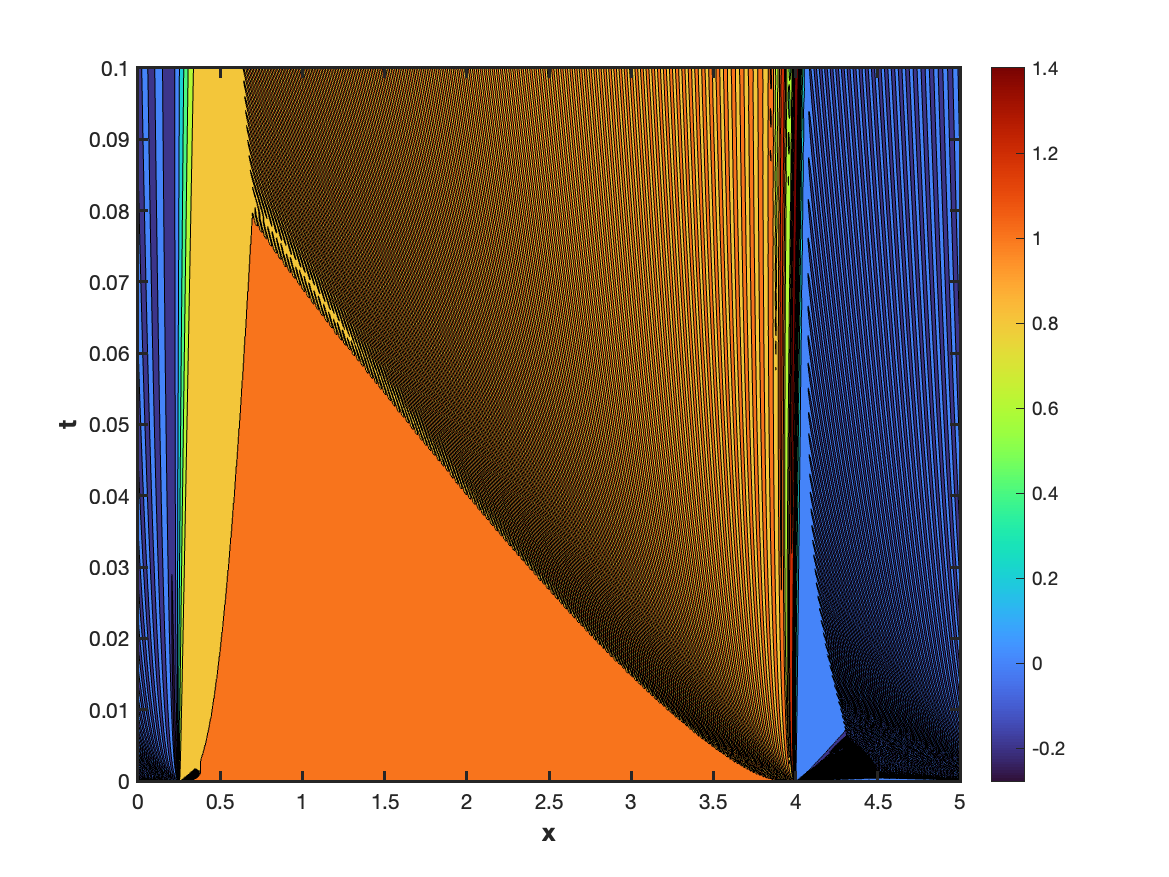}
          \caption*{\normalsize{(f) Contour for wave simulation}}
      \end{subfigure}
      \caption{Numerical solution of KdV equation with $\epsilon = 10^{-4}$, $N=1500$ at $T=0.01,0.05,0.1$ and $\lambda \Delta x = 0.02$ in $x \in [0,5]$, WENO-E for initial condition (\ref{eqn:3b}) of Example \ref{example:3}}\label{eqn:F_3b}
\end{figure}
 The dispersive non-linear KdV equation with initial condition  (\ref{eqn:3a}) with periodic boundary condition has a dispersive shock wave behavior and produces continuous wavelets in the vicinity of the discontinuity for small $\epsilon$. This example demonstrates the effectiveness of our scheme in accurately resolving high-frequency wavelets even for very small values of $\epsilon$. The solution is computed with $\epsilon=10^{-4},10^{-5},10^{-6}$, respectively at $T=0.5$ and plotted in Figure \ref{eqn:F_3a}. The findings indicate that the simulations accurately capture physical oscillations, and the solutions are devoid of noise, especially before and after the occurrence of dispersive shocks. Stable numerical oscillations can be observed in both downstream and upstream of the continuous wavelets in solutions obtained with significantly coarser grids.\par
Upon considering the aforementioned equation with a discontinuous initial condition and inflow-outflow boundary condition, we can see a dispersion-shock wave propagating to the left at each discontinuous interface. The amplitude of the wave increases over time. Figure \ref{eqn:F_3b} demonstrates how a top-hat initial solution transforms into a series of traveling waves. Even with physically discontinuous initial data, the zero-dispersion limit solutions do not exhibit discontinuity. Instead, they develop continuous fine-scale wavelets and eventually break up into solitary waves. The dispersion-shock wave plotted with the current scheme is more accurate and flexible than those plotted in previous studies \cite{MN, MRN, AQ} with the same mesh.
\end{example}
%--------------------------------------------------------------------------------------
% Example 5
\begin{example}\label{example:4a}
\normalfont
The canonical traveling wave solution for the $K(2,2)$ equation 
\begin{equation}
u_t + (u^2)_x+(u^2)_{xxx} =0,
\end{equation}
is given by the compacton
\begin{equation*}
u(x,t)=\begin{cases}
			\frac{4c}{3} \cos^2\left(\frac{x-ct}{4}\right), & \text{if}  \quad |x-ct| \le 2\pi, \\
			0, & \text{else.}  
		\end{cases} 
 \end{equation*}
 The accuracy of the scheme is measured away from the interference of the non-smooth interfaces in the interval $[0,2\pi]$ at $T=\pi/2$ with the periodic boundary condition. Figure \ref{eqn:F_4a} and Table \ref{eqn:T_4a} present a comparison of the errors in the $L^{\infty}$- and $L^{1}$-norms between the WENO-Z scheme and WENO-E scheme with different values of  $\lambda \Delta x$ $(0.02, 0.04, 0.06, 0.1)$. WENO-E scheme with $\lambda \Delta x =0.02$ produces better accuracy as plotted in Figure \ref{eqn:F_4b}.
  \begin{figure}[h!]
\captionsetup[subfigure]{justification=centering}
    \centering
      \begin{subfigure}{0.3\textwidth}
        \includegraphics[trim=0.1cm 0cm 0cm -0.6cm, clip=true,width=\linewidth]{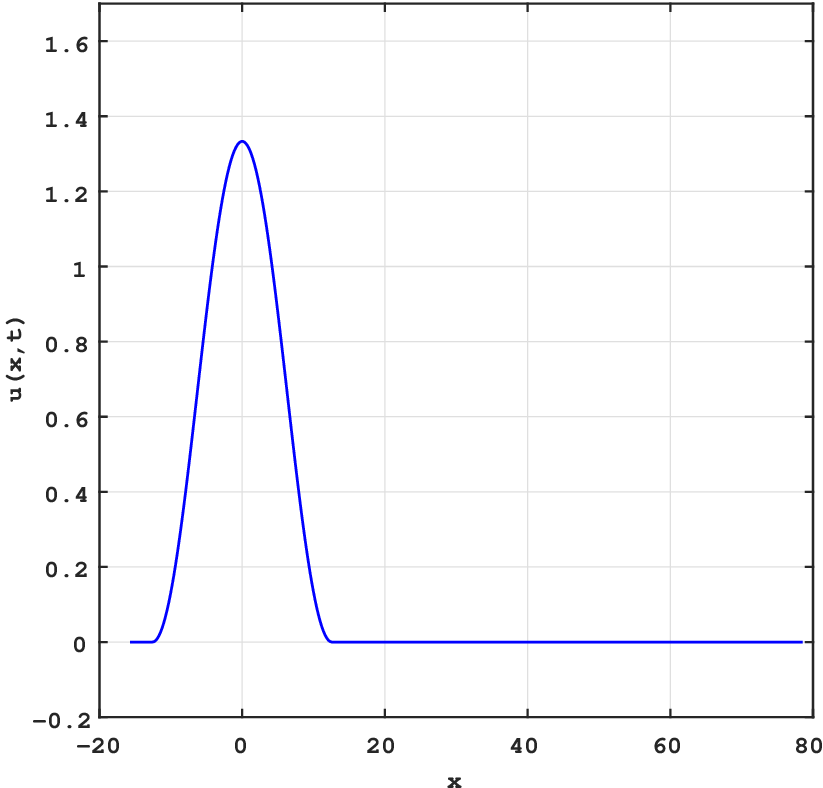} 
          \caption*{\normalsize{(a) $T=0$}}
      \end{subfigure}\hfill
      \begin{subfigure}{0.3\textwidth}
        \includegraphics[trim=2.6cm 0.3cm 3cm 0.9cm, clip=true,width=\linewidth]{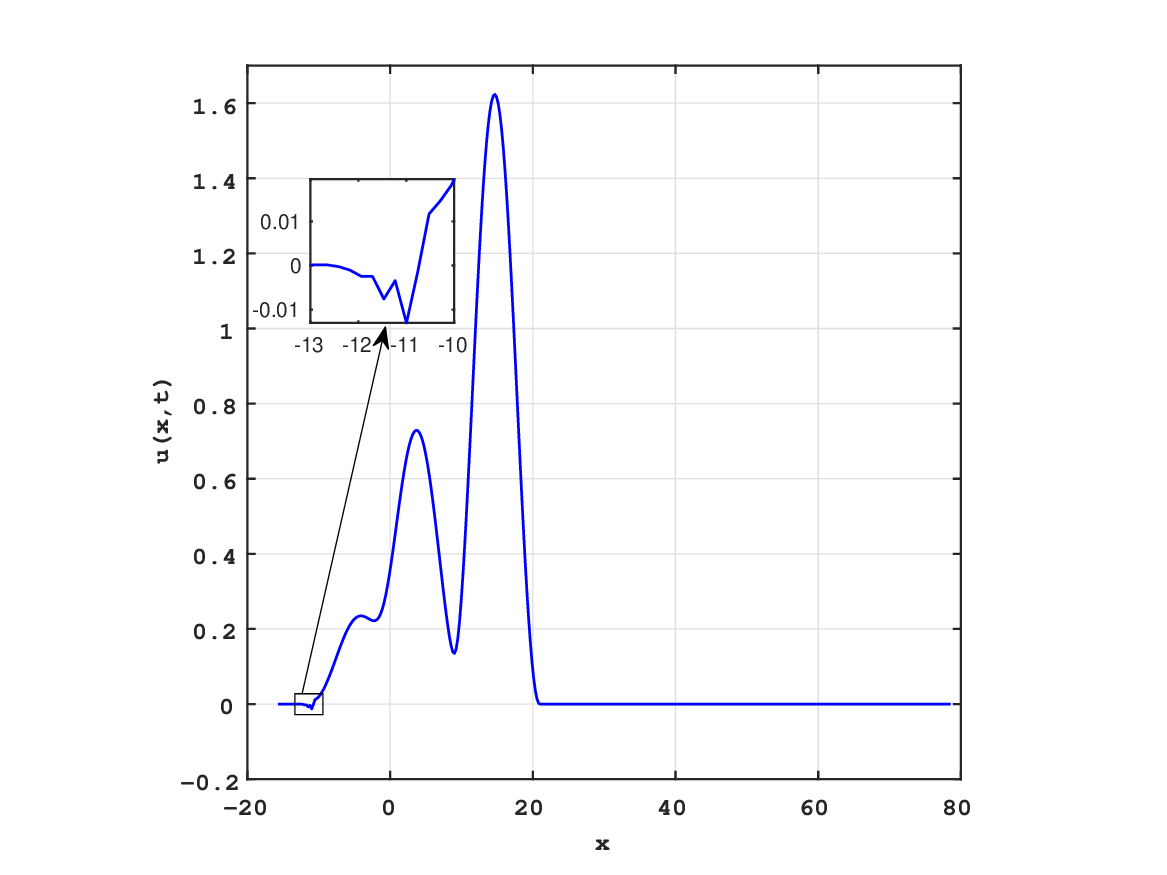}
          \caption*{\normalsize{(b) $T=10$}}
      \end{subfigure}\hfill
      \begin{subfigure}{0.3\textwidth}
        \includegraphics[trim=2.6cm 0.3cm 3cm 1.1cm, clip=true,width=\linewidth]{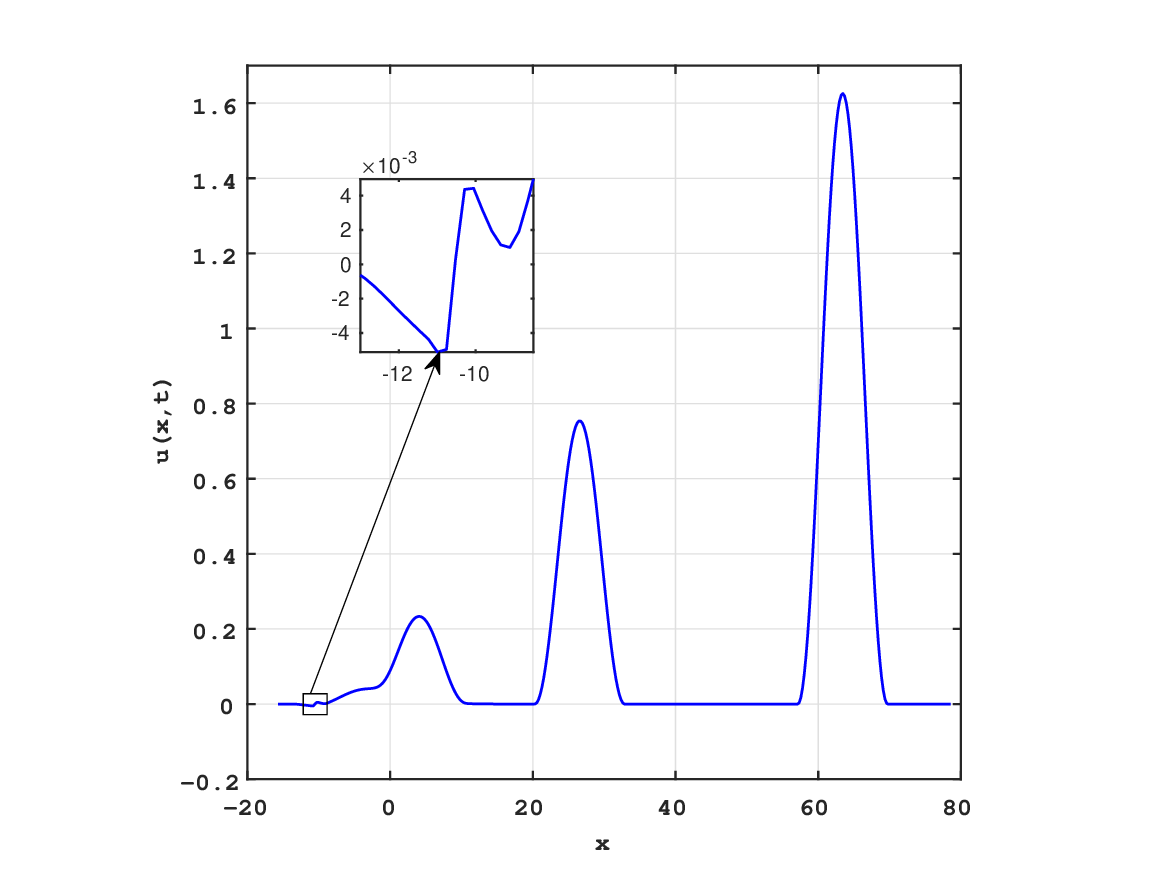}
          \caption*{\normalsize{(c) $T=50$}}
      \end{subfigure}\hfill
      \begin{subfigure}{0.3\textwidth}
        \includegraphics[trim=2.6cm 0.3cm 3cm 1.1cm, clip=true,width=\linewidth]{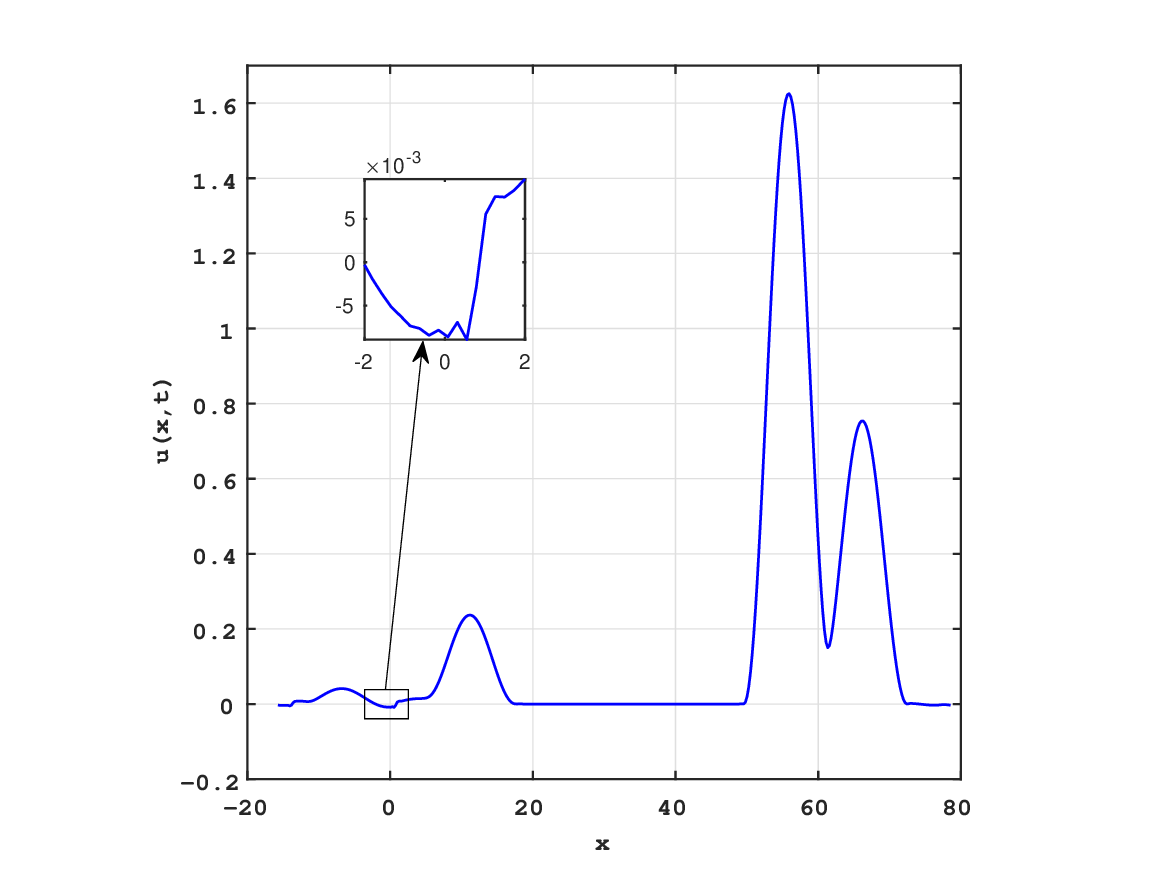}
          \caption*{\normalsize{(d) $T=120$}}
      \end{subfigure}\hfill
            \begin{subfigure}{0.3\textwidth}
        \includegraphics[trim=0.5cm 0.5cm 1cm 1cm, clip=true,width=\linewidth]{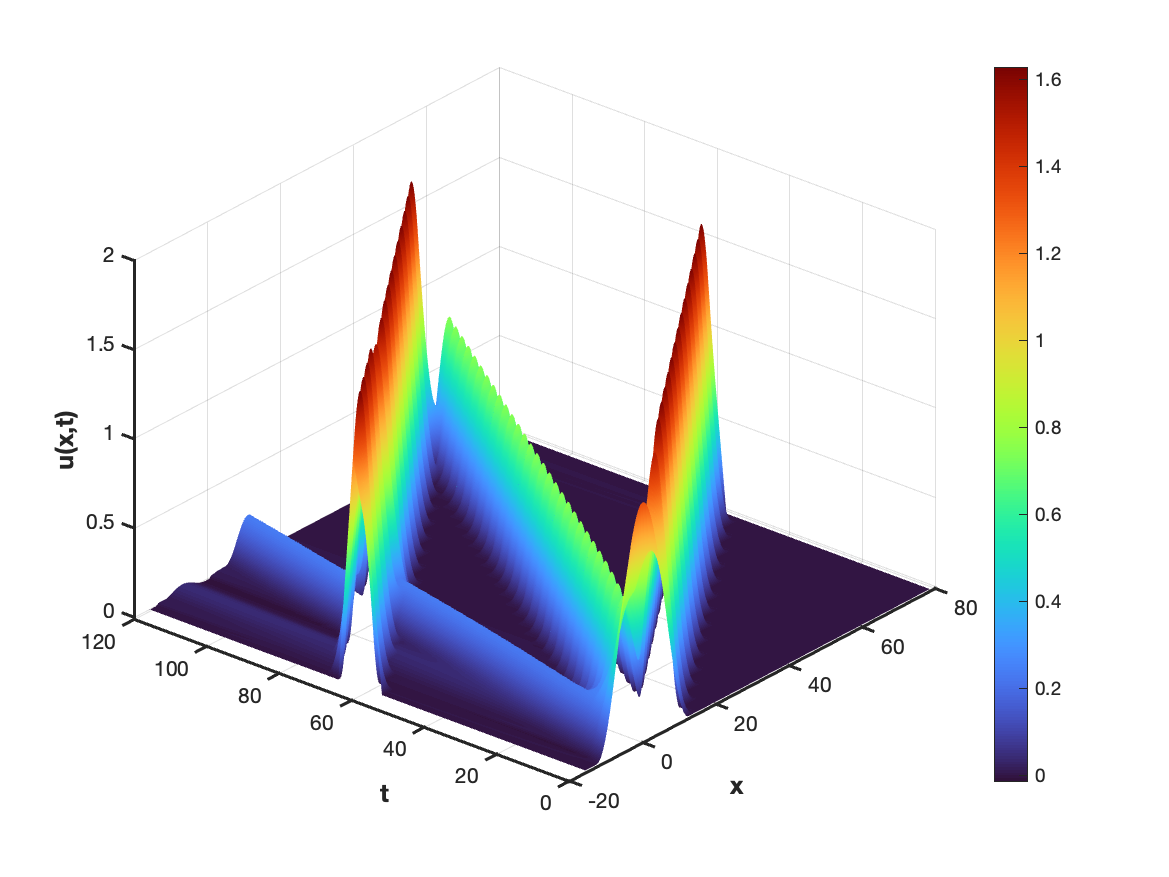}
          \caption*{\normalsize{(e) Wave simulation across various T}}
      \end{subfigure}\hfill
           \begin{subfigure}{0.3\textwidth}
        \includegraphics[width=\linewidth]{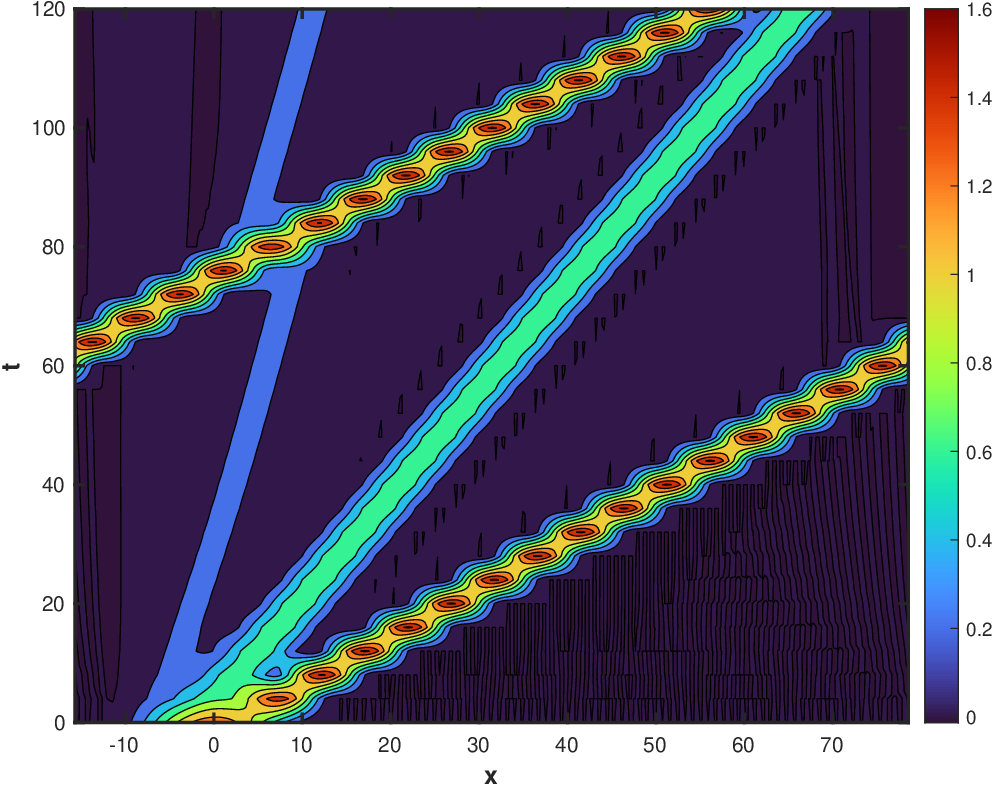}
          \caption*{\normalsize{(f) Contour for wave simulation}}
      \end{subfigure}
      \caption{Numerical solution of K(2,2) equation with $N = 400$ at $T=0,10,50,120$ and $\lambda \Delta x =0.02$ WENO-E in $x \in [-5\pi,25\pi]$ for initial condition (\ref{example:4b}) of Example \ref{example:4a}.}\label{eqn:F_4c}
\end{figure} 
\begin{figure}[h!]
\captionsetup[subfigure]{justification=centering}
    \centering
      \begin{subfigure}{0.3\textwidth}
        \includegraphics[width=\linewidth]{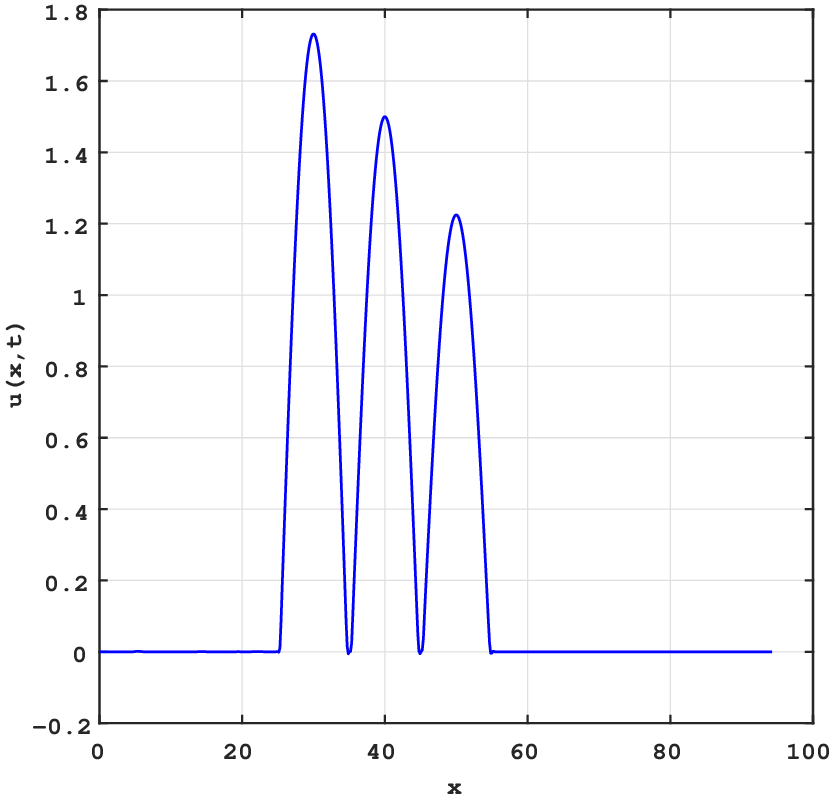}
   \caption*{\normalsize{(a) $T=10$}}
      \end{subfigure}\hfill
      \begin{subfigure}{0.3\textwidth}
       \includegraphics[width=\linewidth]{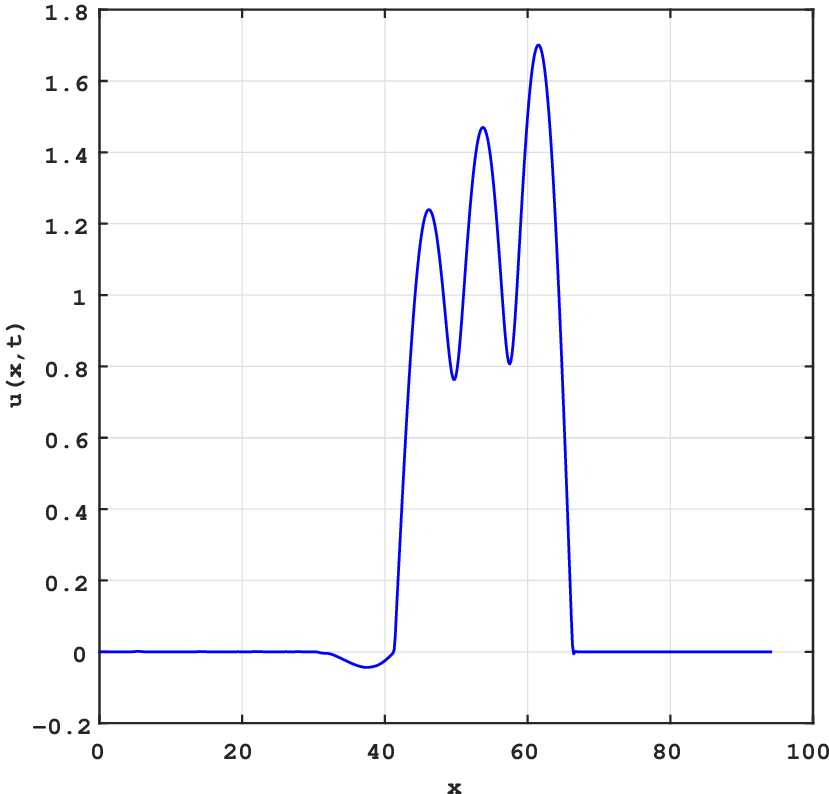}
  \caption*{\normalsize{(b) $T=25$}}
      \end{subfigure}\hfill
      \begin{subfigure}{0.3\textwidth}
        \includegraphics[width=\linewidth]{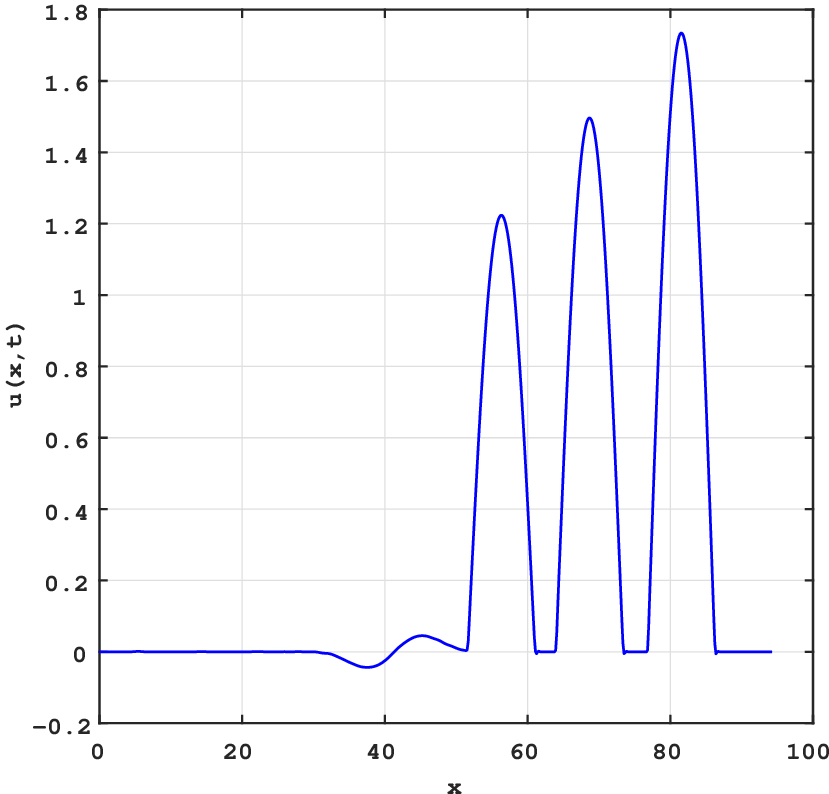}
  \caption*{\normalsize{(c) $T=50$}}
      \end{subfigure}\hfill
      \begin{subfigure}{0.45\textwidth}
        \includegraphics[trim=0cm 0cm 1cm 0cm, clip=true,width=\linewidth]{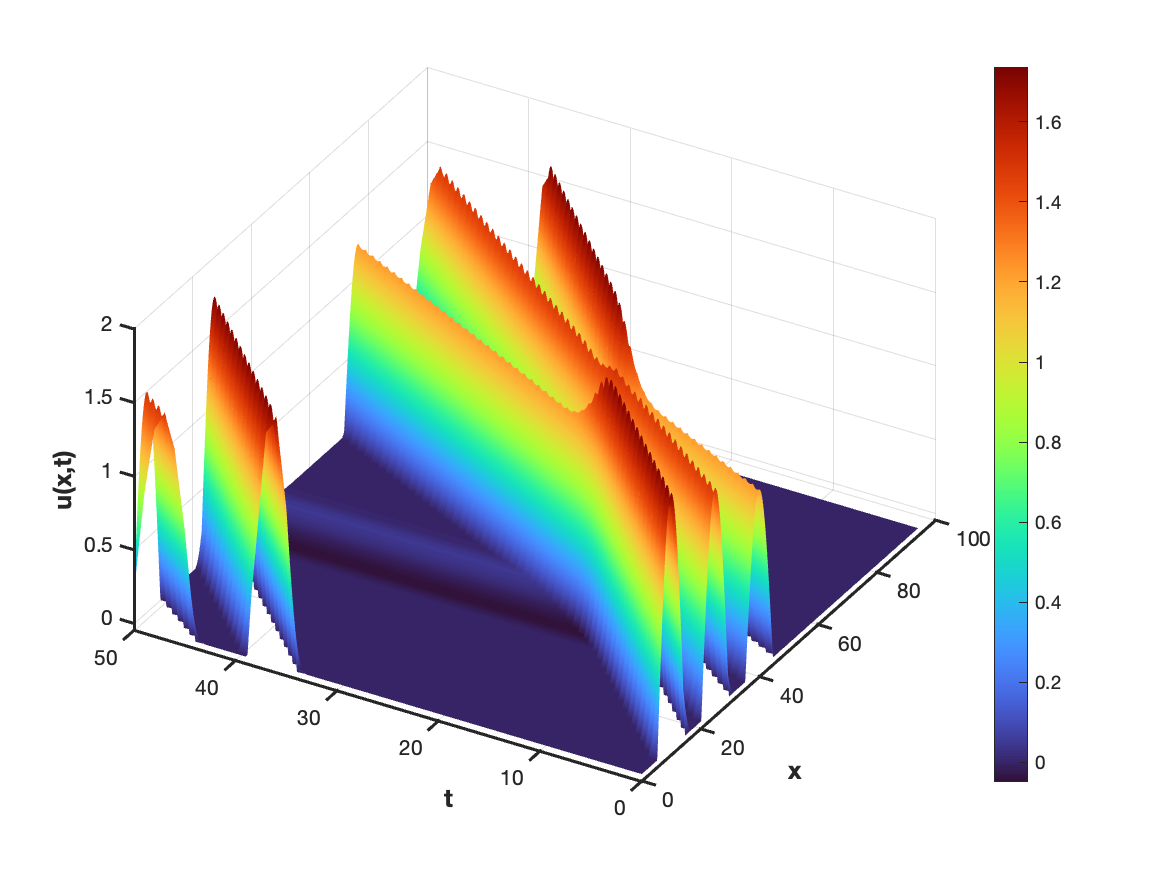} 
    \caption*{\normalsize{(d) Wave simulation across various T}}
      \end{subfigure}\hfill
            \begin{subfigure}{0.45\textwidth}
        \includegraphics[width=\linewidth]{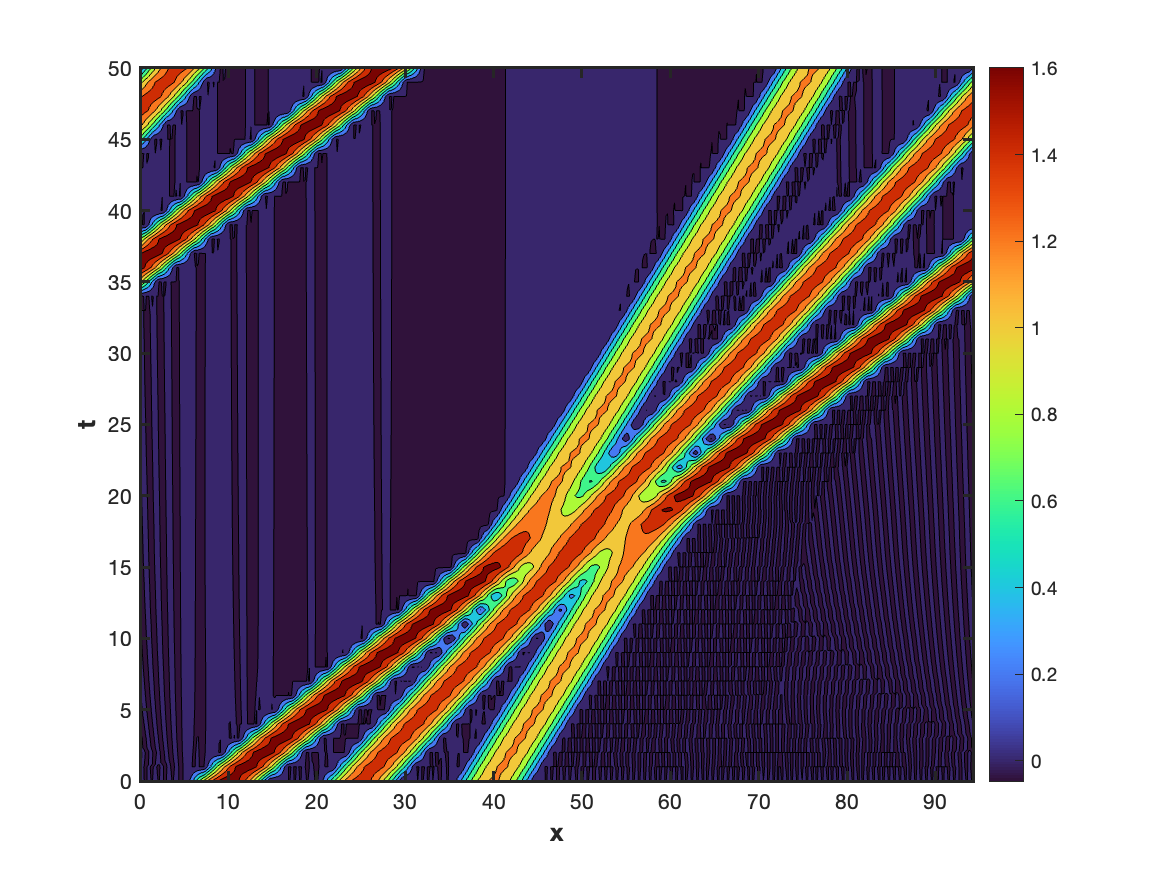}
  \caption*{\normalsize{(e) Contour for wave simulation}}
      \end{subfigure}
     \caption{Numerical solution of K(3,3) equation with $N = 600$ at $T=10,25,50$  and $\lambda \Delta x =0.02$ WENO-E in $x \in [0,30\pi]$ for Example \ref{example:5}.}\label{eqn:F_5a}
\end{figure} 
\begin{figure}[ht!]
\begin{center}
\minipage{0.42\textwidth}
\includegraphics[trim=0cm 0cm 0cm 0cm, clip=true,width=\linewidth]{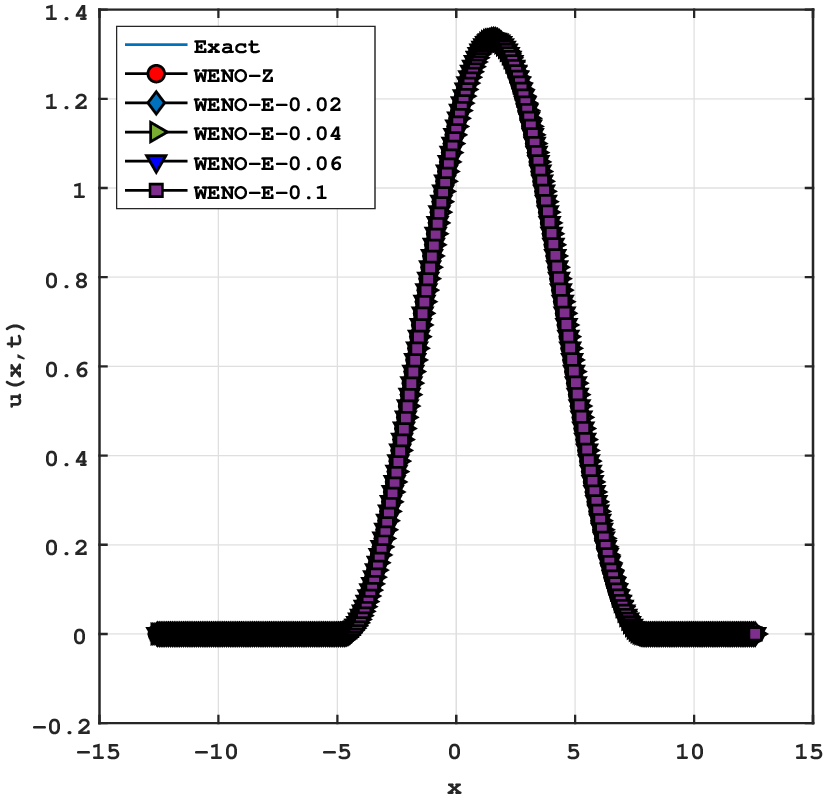}
\caption*{\normalsize{(a) Numerical solution}}
\endminipage \hfill
\minipage{0.47\textwidth}
\includegraphics[trim= 0cm 0cm 0cm 0cm, clip=true,width=\linewidth]{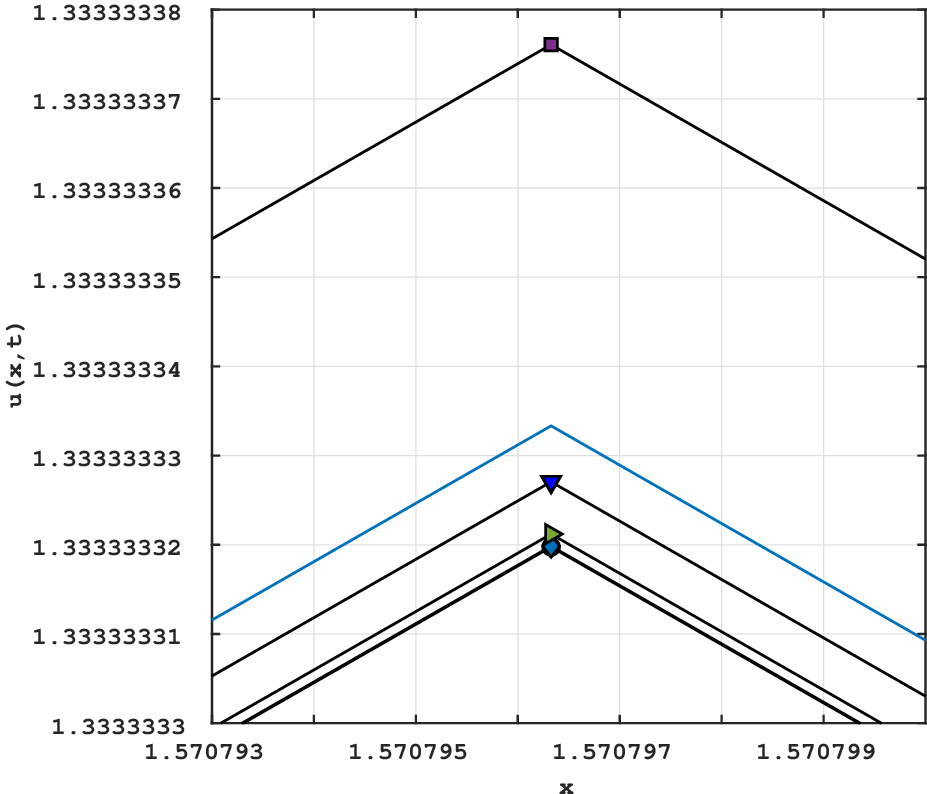}
\caption*{\normalsize{(b) Zoom in near smooth extrema}}
  \endminipage
\end{center}
\caption{ Comparison of numerical solutions obtained from using WENO-Z and WENO-E schemes with $N=320$ at $T=\pi/2$ in $x \in [-4\pi,4\pi]$, for Example \ref{example:4a}}\label{eqn:F_4b}
\end{figure}

\begin{figure}[ht!]
\begin{center}
\minipage{0.45\textwidth}
  \includegraphics[width=\linewidth]{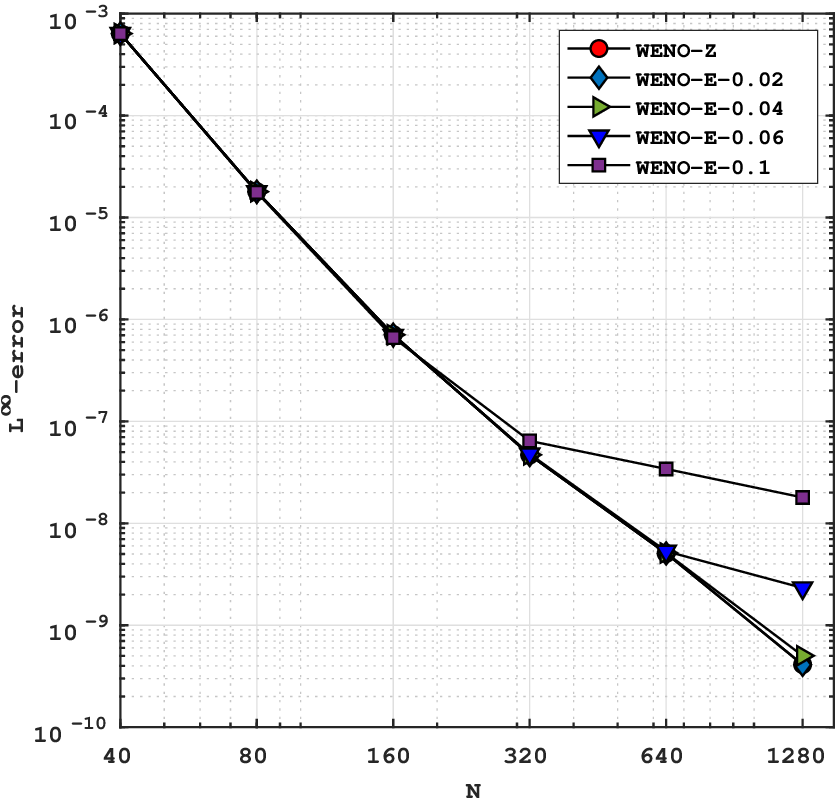}
  \subcaption*{\normalsize{\centering (a) Convergence plot: $L^{\infty}$- error}}
\endminipage\hfill
\minipage{0.45\textwidth}%
  \includegraphics[width=\linewidth]{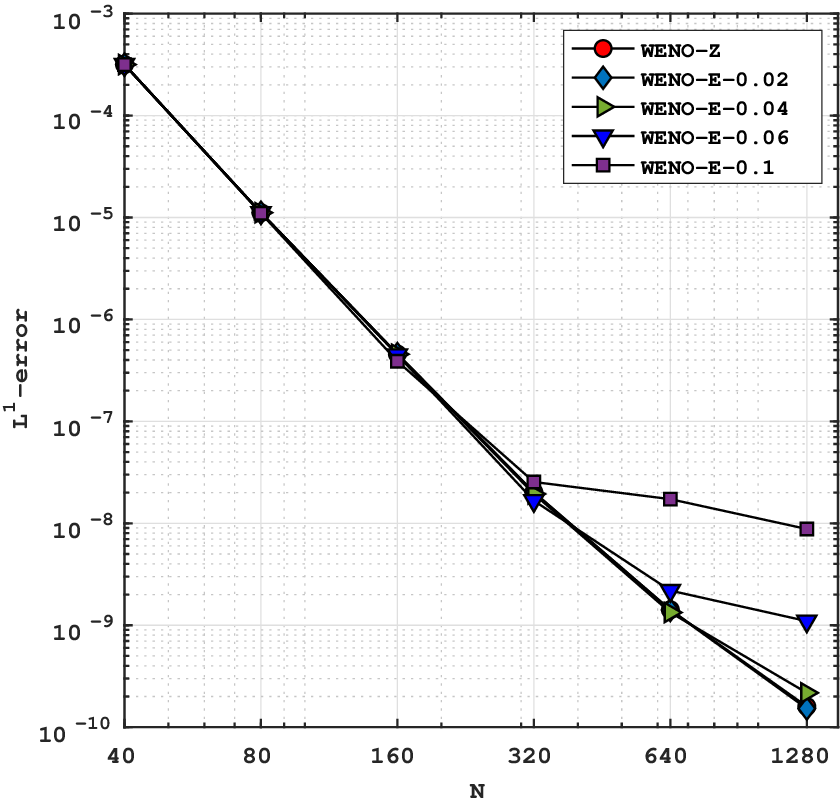}
  \subcaption*{\normalsize{\centering (b) Convergence plot: $L^{1}$- error}}
\endminipage
\caption{Comparison WENO-Z and WENO-E schemes in terms of $L^1$ and $L^{\infty}$ errors (in $\log_{10}$ scale) for Example \ref{example:4a} at $T=\pi/2$.}
\label{eqn:F_4a}
\end{center}
\end{figure}

\begin{figure}[h!]
\captionsetup[subfigure]{justification=centering}
    \centering
      \begin{subfigure}{0.3\textwidth}
        \includegraphics[trim=0cm 0cm 0cm 0cm, clip=true,width=\linewidth]{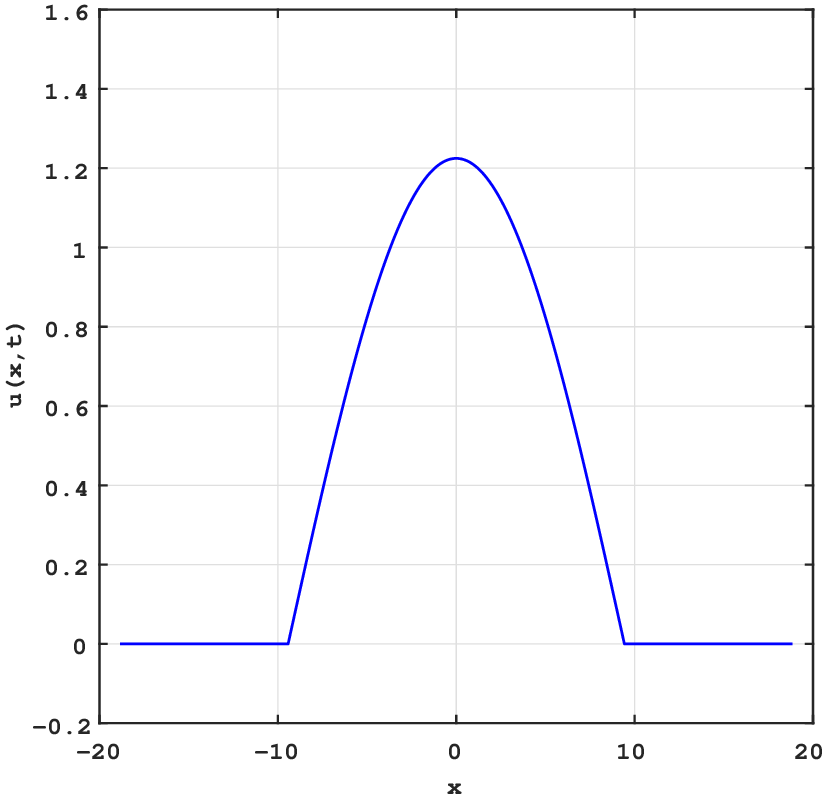} 
          \caption*{\normalsize{(a) $T=0$}}
      \end{subfigure}\hfill
      \begin{subfigure}{0.3\textwidth}
         \includegraphics[trim=0cm 0cm 0cm 0cm, clip=true,width=\linewidth]{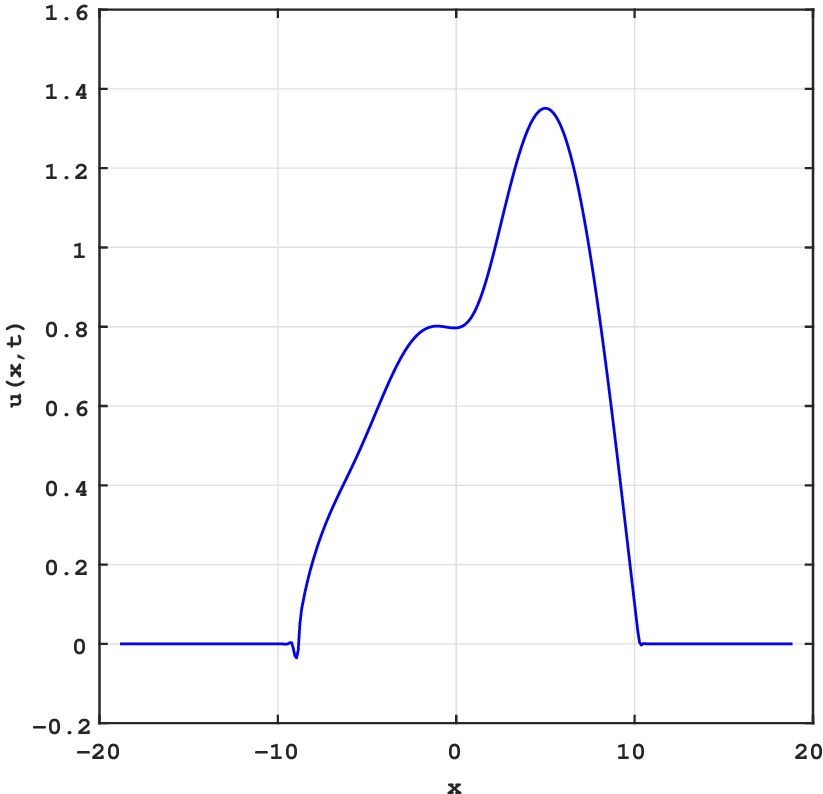}
          \caption*{\normalsize{(b) $T=2$}}
      \end{subfigure}\hfill
      \begin{subfigure}{0.3\textwidth}
         \includegraphics[trim= 0cm 0cm 0cm 0cm, clip=true,width=\linewidth]{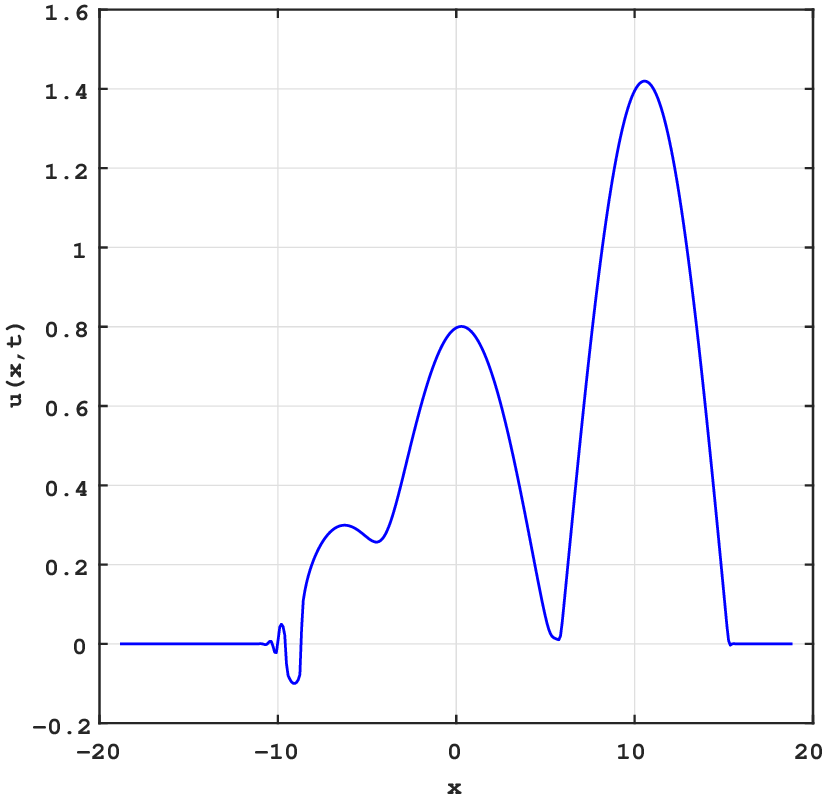} 
          \caption*{\normalsize{(c) $T=6$}}
      \end{subfigure}\hfill
      \begin{subfigure}{0.3\textwidth}
         \includegraphics[trim= 0cm 0cm 0cm 0cm, clip=true,width=\linewidth]{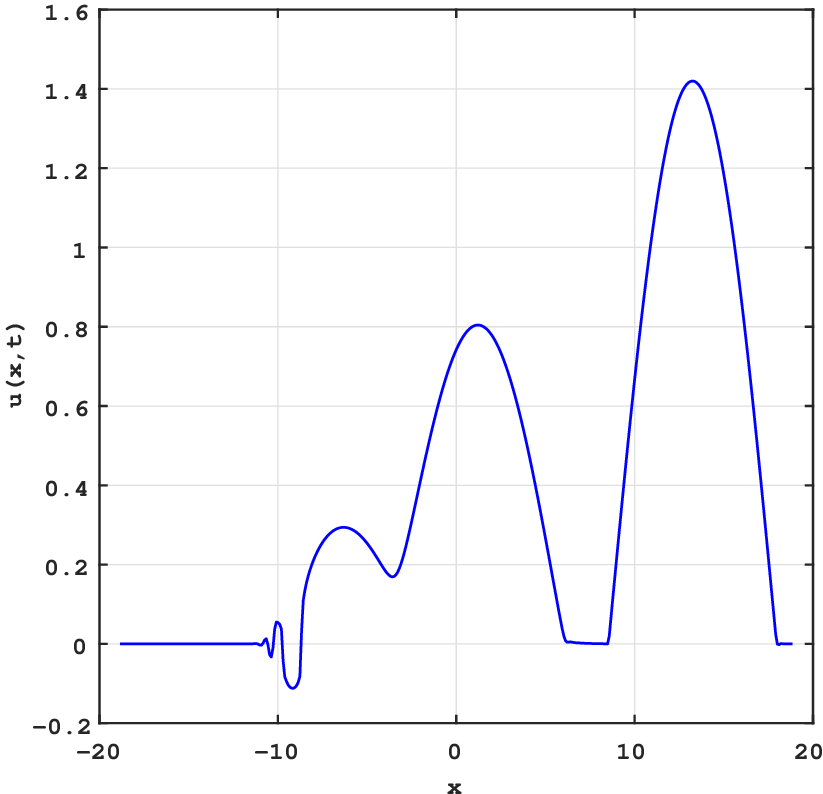} 
          \caption*{\normalsize{(d) $T=8$}}
      \end{subfigure}\hfill
            \begin{subfigure}{0.3\textwidth}
       \includegraphics[trim= 0cm 0cm 0cm 1cm, clip=true,width=\linewidth]{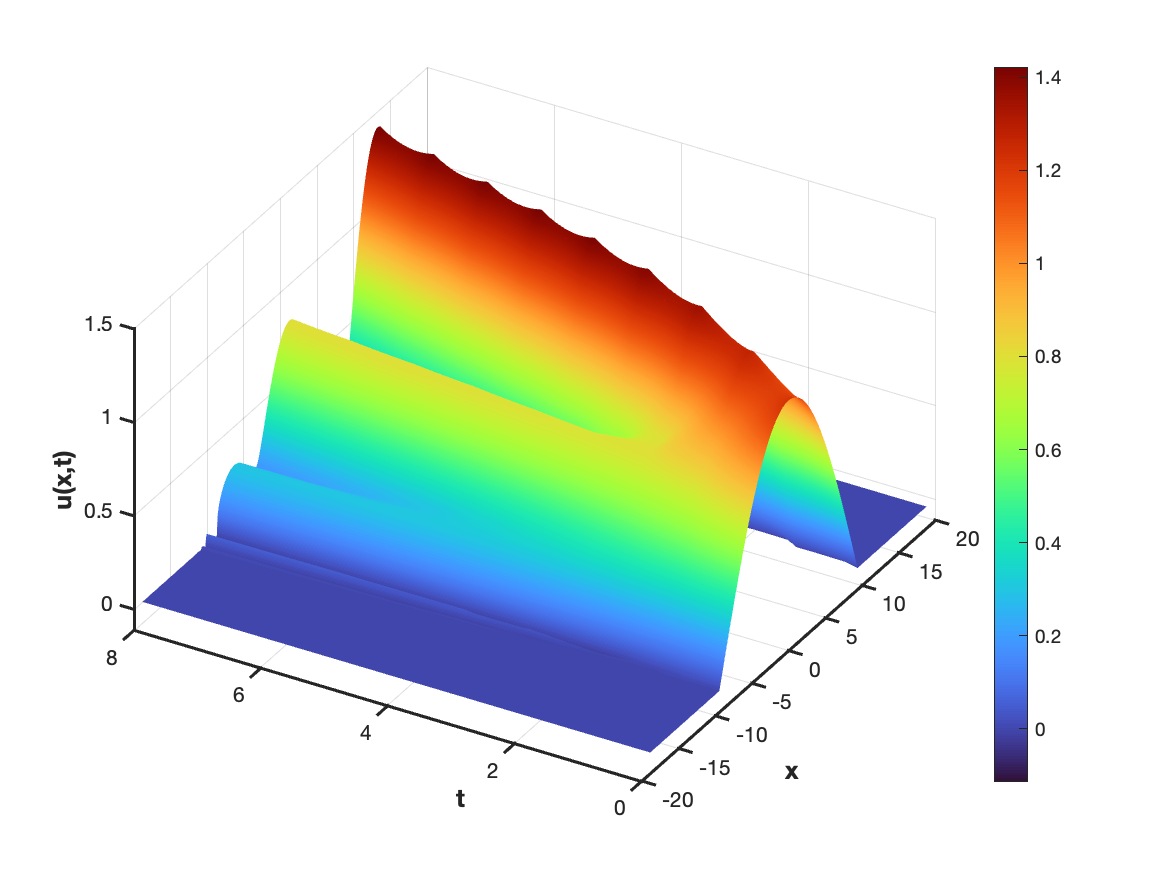} 
          \caption*{\normalsize{(e) Wave simulation across various T}}
      \end{subfigure}\hfill
           \begin{subfigure}{0.3\textwidth}
         \includegraphics[trim= 0cm 0cm 0cm 1cm, clip=true,width=\linewidth]{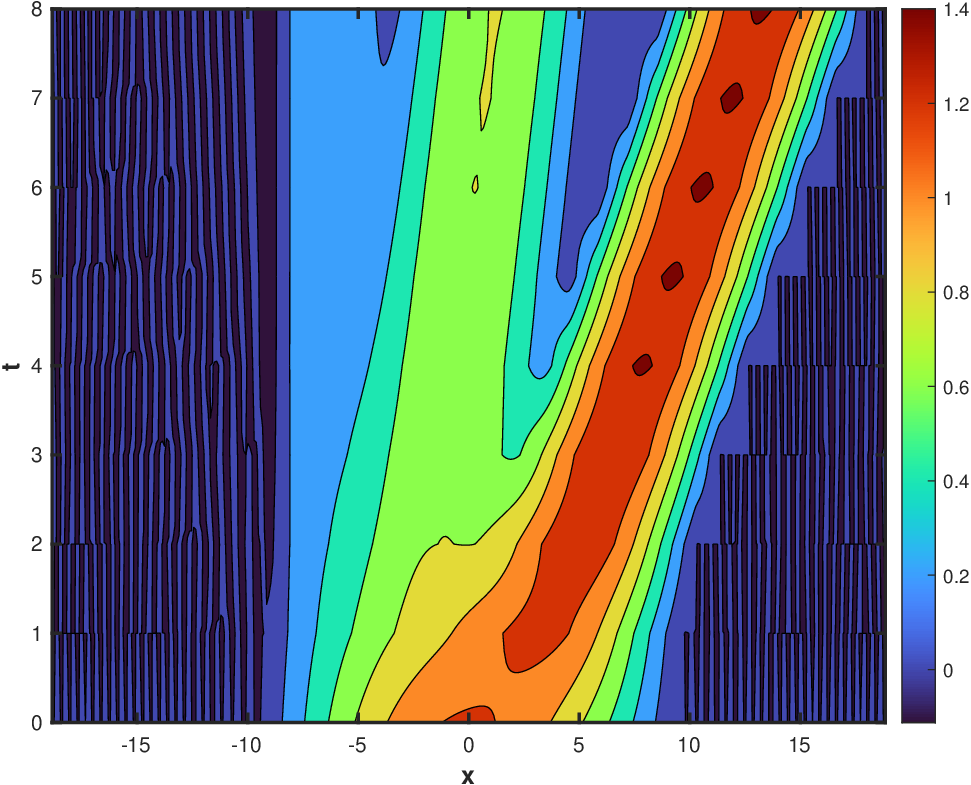} 
          \caption*{\normalsize{(f) Contour for wave simulation}}
      \end{subfigure}
       \caption{Numerical solution of K(3,3) equation with $N = 400$ at $T=0,2,6,8$ and $\lambda \Delta x =0.02$ WENO-E in $x \in [-6\pi, 6\pi]$ for initial condition (\ref{eqn:5b}) of Example \ref{example:5}}\label{eqn:F_5b}
\end{figure} 
To observe the behavior of the decomposition of compacton into compacton- anticompacton pairs, we take the initial data
\begin{equation}
u(x,0)=\begin{cases}
			\frac{4}{3} \cos^2\left(\frac{x}{8}\right), & \text{if} \quad -4\pi \le x \le 4\pi, \\
			0, & \text{else.}  
		\end{cases} 
  \label{example:4b}
 \end{equation}
 In Figure \ref{eqn:F_4c}, wave motion is simulated at $T=0,10, 25,$ and $50$ in the domain of $ [-5\pi,25\pi]$ with $N=400$ cells. The compactons are observed to split from one initial compacton and move to the right over time. A small residue is also developed at the left interface, similar to the LDG scheme \cite{LSY, AQ}, but no Gibbs oscillations occur in the non-smooth interface (i.e., the edges of the compactons). The small residue, which appears to be a compacton-anti-compacton pair on the left side of the compacton-packet, is believed to be physical and was also detected by the LDG scheme \cite{LSY, AQ}. Our scheme can differentiate between numerical and physical oscillation without any additional special treatment. 
\end{example}
%--------------------------------------------------------------------------------------
% Example 6
 \begin{example}\label{example:5}
 \normalfont
 Consider the $K(3,3)$ equation
\begin{equation}
u_t + (u^3)_x+(u^3)_{xxx} =0.
\end{equation}
In order to observe the interaction between three compactons, we take the initial data as
\begin{equation*}
u(x,0)=\begin{cases}
		\sqrt{3} \cos\left(\frac{x-10}{3}\right), & \text{if}           \quad |x-10| \le \frac{3\pi}{2}, \\
            \frac{3}{2} \cos\left(\frac{x-25}{3}\right), & \text{if} 
              \quad |x-25| \le \frac{3\pi}{2}, \\
            \sqrt{\frac{3}{2}} \cos\left(\frac{x-40}{3}\right), & \text{if} 
               \quad |x-40| \le \frac{3\pi}{2}, \\
            0, & \text{otherwise.}  
		\end{cases} 
 \end{equation*}
 The exact solution is given by 
 \begin{equation*}
u(x,t)=\begin{cases}
		\pm \sqrt{\frac{3c}{2}} \cos\left(\frac{x-ct}{3}\right), & \text{if}           \quad |x-ct| \le \frac{3\pi}{2}, \\
            0, & \text{otherwise.}  
		\end{cases} 
 \end{equation*}\\
 Figure \ref{eqn:F_5a} shows that as the compactons propagate to the right, three compactons with different speeds  ($c=10,25,40$) pass through each other during non-linear interaction while maintaining their coherent shapes after the collision. Although the compactons emerge intact from the collision, a minor residue  is reflected back from the collision on the left. Similar findings can also be found in \cite{LSY, AQ}. These compactons are discovered to be not fully elastic, as mentioned in \cite{MCCS} and identified in the original compacton study \cite{IT}. The wavelets of the residue become even more apparent with mesh refinement. Based on earlier compacton studies and numerical observations, we conclude that this phenomenon is not numerically induced.

\par Next, we solve (\ref{example:5}) subject to the initial data
\begin{equation}
u(x,0)=\begin{cases}
			\sqrt{\frac{3}{2}} \cos\left(\frac{x}{6}\right), & \text{if} \quad -3\pi \le x \le 3\pi, \\
			0, & \text{else.}  
		\end{cases} 
  \label{eqn:5b}
 \end{equation}
The results of our numerical simulations are shown in Figure \ref{eqn:F_5b}. As time evolves, a train of canonical compactons splits from the initial data and moves to the right. At the same time, a rapid oscillation develops at the left interface of the initial data. The oscillations that we get on the left side of the solution indicate that the oscillations are an integral part of the solution and not a numerical artifact.
\par After examining the six aforementioned examples, it becomes evident that there is no universally optimal value for the tension parameter, denoted as $\lambda \Delta x$. Its effectiveness is contingent upon the specific characteristics of the equation under consideration. In linear scenarios, regardless of whether the number of grid points $N$ is small or large, we observe improved accuracy within the range of $\lambda \Delta x$ between 0.02 and 0.04. Conversely, in non-linear cases, superior accuracy can be achieved by setting $\lambda \Delta x$ to 0.1 for low values of $N$ and maintaining the range of $\lambda \Delta x$ between 0.02 and 0.04 for high values of $N$. More details are provided in Appendix:A. In summary, it is advisable to employ the range, $\lambda \Delta x \in [0.02, 0.04]$, to ensure the optimal performance of the proposed scheme.
\end{example}
\section{Conclusion}
\label{sec:con}
In this study, a fifth-order WENO scheme based on exponential polynomials for solving the non-linear dispersion-type equations has been proposed. The proposed approximation space utilizes the exponential basis with a tension parameter that can be optimized to better fit the specific features of the characteristics of the initial data, resulting in improved results without the presence of spurious oscillations. Numerical tests were conducted to simulate several equations with various initial conditions, including the 1-D and 2-D linear and non-linear KdV equations, $K(2,2)$ equation, and $K(3,3)$ equation. In direct comparison to the WENO-Z method, our proposed method WENO-E, when executed with a optimized tension parameter, has superior capacity for resolving non-linear dispersion-type equations at higher resolutions.

\section*{Acknowledgements}
The author Samala Rathan is supported by NBHM, DAE, India (Ref. No. 02011/46/2021 NBHM(R.P.)/R $\&$ D II/14874) and IIPE, Visakhapatnam, India (IRG Grant No. IIPE/DORD/IRG/001). 

\section*{Conflict of interest}
The authors declare no potential conflict of interests.

\section*{Data availability}
The data that support the findings of this study are available upon reasonable request.

\section*{Appendix:A}
In this section, we provide a detailed tabulated values that enumerates the $L^1$-errors associated with a range of $\lambda\Delta x$ values in Table \ref{app:1}, \ref{app:2}, \ref{app:3}, and  \ref{app:4} for examples   \ref{example:1},  \ref{example:2D}, \ref{example:2}, and \ref{example:4a} respectively.   The chosen spectrum of $\lambda\Delta x$ values is carefully selected to ensure that there is a balanced and adequate representation of errors. This comprehensive overview allows us to effectively showcase the impact and significance of varying $\lambda\Delta x$. The optimal value of $\lambda\Delta x$ are highlighted for each example in bold.   To  see a universally applicable parameter, we've categorized $\lambda \Delta x$ values into three groups which are as follows:
\begin{itemize}
\item Category A (Cat-A) is $ \lambda \Delta x \in [0.02,0.04]$ which comprises optimal values that consistently yield accuracy across various grid point counts $N$.
\item Category B (Cat-B) is $ \lambda \Delta x \in  [0.06, 0.1]$ which  includes values closely resembling Cat-A at lower grid points $N$, but they neither improve nor maintain accuracy as $N$ increases, ultimately deteriorating.
\item Category C (Cat-C) is $ \lambda \Delta x >0.1$ which  encompasses values that exhibit inconsistency and substantial
deviation from the actual solution, numerically distant from Cat-A.
\end{itemize}
As per these observations, the choice of the tension parameter, $\lambda \Delta x$, does not possess a universally optimal value, as its effectiveness relies on the specific characteristics of the equation under consideration. In linear scenarios, whether the number of grid points $N$ is small or large, we observe enhanced accuracy within the $\lambda \Delta x$ range of $0.02$ to $0.04$. Conversely, in nonlinear cases, superior accuracy is achieved by setting $\lambda \Delta x$ to $0.1$ for low $N$ and adhering to the $\lambda \Delta x$ range of $0.02$ to $0.04$ for high $N$. In overall, one can utilize the Cat-A  i.e., $ \lambda \Delta x \in [0.02,0.04]$  to see the good performance of the proposed scheme.
 \begin{table}[htbp!]
\centering
\captionof{table}{WENO-E schemes in terms of $L^1$- errors for Example \ref{example:1} over the domain $\Omega = [0,2\pi]$ at time $T=1$.}\label{app:1}
\begin{tabular}{|^c|*{7}{_c|}}\hline
$\boldsymbol{L^{1}}$  & \multicolumn{7}{c|}{$\boldsymbol{N}$}\\
\cline{2-8}
\textbf{error}
&\makebox[3em]{10}&\makebox[3em]{20}&\makebox[3em]{40}
&\makebox[3em]{80}&\makebox[3em]{160}&\makebox[3em]{320}&\makebox[3em]{640}\\\hline
WENO-Z  &1.7456e-03&5.6856e-05&1.7815e-06&5.5652e-08&1.7390e-09&7.0742e-11 & 1.3677e-10\\\hline
WENO-E-0.01 & 1.7519e-03 & 5.7101e-05 & 1.7825e-06 & 5.5640e-08 & 1.7317e-09 & 6.7075e-11 & 1.3520e-10 \\\hline
\rowstyle{\bfseries}
WENO-E-0.02 & 1.7519e-03 & 5.7100e-05 & 1.7821e-06 & 5.5421e-08 & 1.6221e-09 & 1.3413e-11 & 1.0871e-10 \\\hline
WENO-E-0.03 & 1.7519e-03 & 5.7097e-05 & 1.7802e-06 & 5.4471e-08 & 1.1474e-09 & 2.2517e-10 & 3.5957e-11 \\\hline
WENO-E-0.04 & 1.7520e-03 & 5.7089e-05 & 1.7751e-06 & 5.1911e-08 & 1.3076e-10 & 8.6372e-10 & 3.3622e-10 \\\hline
WENO-E-0.05 & 1.7520e-03 & 5.7069e-05 & 1.7643e-06 & 4.6515e-08 & 2.8258e-09 & 2.2103e-09 & 1.0080e-09 \\\hline
WENO-E-0.06 & 1.7520e-03 & 5.7032e-05 & 1.7447e-06 & 3.6701e-08 & 7.7266e-09 & 4.6589e-09 & 2.2314e-09 \\\hline
WENO-E-0.07 & 1.7520e-03 & 5.6970e-05 & 1.7124e-06 & 2.0541e-08 & 1.5797e-08 & 8.6913e-09 & 4.2465e-09 \\\hline
WENO-E-0.08 & 1.7520e-03 & 5.6874e-05 & 1.6629e-06 & 4.2486e-09 & 2.8176e-08 & 1.4877e-08 & 7.3376e-09 \\\hline
WENO-E-0.09 & 1.7519e-03 & 5.6734e-05 & 1.5908e-06 & 4.0299e-08 & 4.6179e-08 & 2.3872e-08 & 1.1833e-08 \\\hline
WENO-E-0.1 & 1.7516e-03 & 5.6537e-05 & 1.4902e-06 & 9.0593e-08 & 7.1295e-08 & 3.6421e-08 & 1.8105e-08 \\\hline
WENO-E-0.2 & 1.7378e-03 & 4.7842e-05 & 2.8973e-06 & 2.2839e-06 & 1.1666e-06 & 5.8367e-07 & 2.9160e-07 \\\hline
WENO-E-0.3 & 1.6711e-03 & 1.0023e-05 & 2.1903e-05 & 1.1783e-05 & 5.9103e-06 & 2.9538e-06 & 1.4761e-06 \\\hline
WENO-E-0.4 & 1.4862e-03 & 9.1817e-05 & 7.3016e-05 & 3.7329e-05 & 1.8667e-05 & 9.3276e-06 & 4.6615e-06 \\\hline
WENO-E-0.5 & 1.0921e-03 & 3.0621e-04 & 1.8055e-04 & 9.1071e-05 & 4.5503e-05 & 2.2736e-05 & 1.4356e-03 \\\hline
WENO-E-0.6 & 3.7362e-04 & 6.9462e-04 & 3.7529e-04 & 1.8838e-04 & 9.4091e-05 & 4.7011e-05 & 1.4620e+32 \\\hline
\end{tabular}
\end{table}
\begin{table}[H]
\centering
\captionof{table}{WENO-E schemes in terms of $L^1$- errors for Example \ref{example:2D} over the domain $\Omega = (0,2\pi)\times(0,2\pi)$ at time $T=1$.}\label{app:2}
\begin{tabular}{|^c|*{5}{_c|}}\hline
$\boldsymbol{L^{1}}$  & \multicolumn{5}{c|}{$\boldsymbol{N_x \times N_y}$}\\
\cline{2-6}
\textbf{error}
&\makebox[3em]{10$\times$10}&\makebox[3em]{20$\times$20}&\makebox[3em]{40$\times$40}&\makebox[3em]{80$\times$80}&\makebox[3em]{160$\times$160}\\\hline
WENO-Z & 3.4560e-03 & 1.1104e-04 & 3.5318e-06 & 1.1080e-07 & 3.4732e-09 \\\hline
WENO-E-0.01 & 3.4272e-03 & 1.1149e-04 & 3.5338e-06 & 1.1077e-07 & 3.4586e-09 \\\hline	
WENO-E-0.02 & 3.4272e-03 & 1.1149e-04 & 3.5329e-06 & 1.1034e-07 & 3.2400e-09 \\\hline	
WENO-E-0.03 & 3.4273e-03 & 1.1149e-04 & 3.5292e-06 & 1.0845e-07 & 2.2923e-09 \\\hline
\rowstyle{\bfseries}	
WENO-E-0.04 & 3.4273e-03 & 1.1147e-04 & 3.5191e-06 & 1.0335e-07 & 2.5903e-10 \\\hline	
WENO-E-0.05 & 3.4273e-03 & 1.1143e-04 & 3.4977e-06 & 9.2607e-08 & 5.6387e-09 \\\hline	
WENO-E-0.06 & 3.4273e-03 & 1.1136e-04 & 3.4589e-06 & 7.3069e-08 & 1.5421e-08 \\\hline	
WENO-E-0.07 & 3.4272e-03 & 1.1124e-04 & 3.3949e-06 & 4.0895e-08 & 3.1531e-08 \\\hline	
WENO-E-0.08 & 3.4270e-03 & 1.1105e-04 & 3.2966e-06 & 8.4582e-09 & 5.6242e-08 \\\hline	
WENO-E-0.09 & 3.4267e-03 & 1.1078e-04 & 3.1537e-06 & 8.0231e-08 & 9.2177e-08 \\\hline	
WENO-E-0.1 & 3.4262e-03 & 1.1039e-04 & 2.9544e-06 & 1.8036e-07 & 1.4231e-07 \\\hline	
WENO-E-0.2 & 3.3973e-03 & 9.3414e-05 & 5.7438e-06 & 4.5470e-06 & 2.3286e-06 \\\hline	
WENO-E-0.3 & 3.2631e-03 & 1.9555e-05 & 4.3422e-05 & 2.3459e-05 & 1.1798e-05 \\\hline	
WENO-E-0.4 & 2.8944e-03 & 1.7936e-04 & 1.4476e-04 & 7.4321e-05 & 3.7262e-05 \\\hline	
WENO-E-0.5 & 2.1115e-03 & 5.9823e-04 & 3.5799e-04 & 1.8133e-04 & 9.0833e-05 \\\hline
WENO-E-0.6 & 6.8564e-04 & 1.3574e-03 & 7.4421e-04 & 3.7510e-04 & 1.8783e-04 \\\hline		
\end{tabular}
\end{table}
\begin{table}[h!]
\centering
\captionof{table} {WENO-E schemes in terms of $L^1$- errors for Example \ref{example:2} over the domain $\Omega =[-10,10]$ at time $T=0.5$.}\label{app:3}
\begin{tabular}{|^c|*{6}{_c|}}\hline
$\boldsymbol{L^{1}}$  & \multicolumn{6}{c|}{$\boldsymbol{N}$}\\
\cline{2-7}
\textbf{error}
&\makebox[3em]{40}&\makebox[3em]{80}
&\makebox[3em]{160}&\makebox[3em]{320}&\makebox[3em]{640}&\makebox[3em]{1280}\\\hline
WENO-Z & 3.1369e-02 & 2.3854e-02 & 5.9718e-05 & 1.9878e-06& 8.4259e-08& 2.5115e-08 \\\hline
WENO-E-0.01 & 3.0119e-02 & 2.1970e-03 & 6.8675e-05 & 2.0650e-06 & 8.4833e-08 & 2.5116e-08 \\\hline							
WENO-E-0.02 & 3.0119e-02 & 2.1970e-03 & 6.8676e-05 & 2.0648e-06 & 8.4750e-08 & 2.5080e-08 \\\hline							
WENO-E-0.03 & 3.0120e-02 & 2.1971e-03 & 6.8676e-05 & 2.0642e-06 & 8.4390e-08 & 2.4925e-08 \\\hline	
\rowstyle{\bfseries}						
WENO-E-0.04 & 3.0120e-02 & 2.1972e-03 & 6.8676e-05 & 2.0624e-06 & 8.3420e-08 & 2.4566e-08 \\\hline							
WENO-E-0.05 & 3.0120e-02 & 2.1973e-03 & 6.8672e-05 & 2.0584e-06 & 8.1376e-08 & 2.4272e-08 \\\hline							
WENO-E-0.06 & 3.0121e-02 & 2.1974e-03 & 6.8663e-05 & 2.0510e-06 & 7.7660e-08 & 2.5510e-08 \\\hline							
WENO-E-0.07 & 3.0122e-02 & 2.1975e-03 & 6.8645e-05 & 2.0389e-06 & 7.1545e-08 & 2.8296e-08 \\\hline							
WENO-E-0.08 & 3.0122e-02 & 2.1977e-03 & 6.8615e-05 & 2.0201e-06 & 6.2423e-08 & 3.2753e-08 \\\hline							
WENO-E-0.09 & 3.0123e-02 & 2.1978e-03 & 6.8571e-05 & 1.9928e-06 & 5.1088e-08 & 3.9343e-08 \\\hline							
WENO-E-0.1 & 3.0124e-02 & 2.1979e-03 & 6.8506e-05 & 1.9547e-06 & 4.4961e-08 & 4.8688e-08 \\\hline							
WENO-E-0.2 & 3.0136e-02 & 2.1967e-03 & 6.5444e-05 & 7.2544e-07 & 8.6967e-07 & 4.7576e-07 \\\hline							
WENO-E-0.3 & 3.0145e-02 & 2.1807e-03 & 5.1841e-05 & 7.1637e-06 & 4.5905e-06 & 2.3462e-06 \\\hline							
WENO-E-0.4 & 3.0135e-02 & 2.1291e-03 & 2.2468e-05 & 2.6852e-05 & 1.4602e-05 & 7.3805e-06 \\\hline							
WENO-E-0.5 & 3.0085e-02 & 2.0125e-03 & 7.3678e-05 & 6.8390e-05 & 3.5670e-05 & 1.7974e-05 \\\hline				
WENO-E-0.6 & 2.9965e-02 & 1.7948e-03 & 2.1361e-04 & 1.4370e-04 & 7.3836e-05 & 4.3984e-05 \\\hline													
\end{tabular}
\end{table}
 \begin{table}[h!]
\centering
\captionof{table}{WENO-E schemes in terms of $L^1$- errors for Example \ref{example:4a} over the domain $\Omega =[0,2\pi]$ at time $T=\pi/2$.}\label{app:4}
\begin{tabular}{|^c|*{6}{_c|}}\hline
$\boldsymbol{L^{1}}$  & \multicolumn{6}{c|}{$\boldsymbol{N}$}\\
\cline{2-7}
\textbf{error}
&\makebox[3em]{40}&\makebox[3em]{80}
&\makebox[3em]{160}&\makebox[3em]{320}&\makebox[3em]{640}&\makebox[3em]{1280}\\\hline
WENO-Z & 3.1608e-04 & 1.1146e-05 & 4.5757e-07 & 2.0152e-08& 1.4181e-09& 1.5958e-10 \\\hline
WENO-E-0.01 & 3.1654e-04 & 1.1150e-05 & 4.5757e-07 & 2.0149e-08& 1.4175e-09& 1.5902e-10 \\\hline
\rowstyle{\bfseries}			
WENO-E-0.02 & 3.1654e-04 & 1.1150e-05 & 4.5747e-07 & 2.0105e-08& 1.4092e-09& 1.5226e-10 \\\hline			
WENO-E-0.03 & 3.1654e-04 & 1.1149e-05 & 4.5702e-07 & 1.9910e-08& 1.3782e-09& 1.4107e-10 \\\hline			
WENO-E-0.04 & 3.1654e-04 & 1.1146e-05 & 4.5581e-07 & 1.9391e-08& 1.3336e-09& 2.1720e-10 \\\hline			
WENO-E-0.05 & 3.1653e-04 & 1.1141e-05 & 4.5326e-07 & 1.8373e-08& 1.4339e-09& 5.1372e-10 \\\hline			
WENO-E-0.06 & 3.1652e-04 & 1.1131e-05 & 4.4862e-07 & 1.6655e-08& 2.1872e-09& 1.0973e-09 \\\hline			
WENO-E-0.07 & 3.1649e-04 & 1.1115e-05 & 4.4098e-07 & 1.4819e-08& 3.9253e-09& 2.0677e-09 \\\hline			
WENO-E-0.08 & 3.1644e-04 & 1.1091e-05 & 4.2927e-07 & 1.3932e-08& 6.8062e-09& 3.5647e-09 \\\hline			
WENO-E-0.09 & 3.1638e-04 & 1.1055e-05 & 4.1223e-07 & 1.7225e-08& 1.1160e-08& 5.7471e-09 \\\hline			
WENO-E-0.1 & 3.1629e-04 & 1.1006e-05 & 3.8847e-07 & 2.5432e-08& 1.7251e-08& 8.7921e-09 \\\hline			
WENO-E-0.2 & 3.1212e-04 & 8.8400e-06 & 7.5150e-07 & 5.5608e-07& 2.8347e-07& 1.4158e-07 \\\hline			
WENO-E-0.3 & 2.9387e-04 & 2.2379e-06 & 5.4111e-06 & 2.8729e-06& 1.4366e-06& 7.1667e-07 \\\hline			
WENO-E-0.4 & 2.4463e-04 & 2.5784e-05 & 1.7949e-05 & 9.1032e-06& 4.5375e-06& 2.2632e-06 \\\hline			
WENO-E-0.5 & 1.4122e-04 & 7.8889e-05 & 4.4326e-05 & 2.2210e-05& 1.1061e-05& 5.5166e-06 \\\hline	
WENO-E-0.6 & 5.1057e-05 & 1.7505e-04 & 9.2092e-05 & 4.5940e-05& 2.2871e-05& 1.1407e-05 \\\hline			
\end{tabular}
\end{table}
\clearpage
\begin{landscape}
\thispagestyle{mylandscape}
\section*{Appendix:B}
The matrix $A$ is
\begin{eqnarray*}
  \left[\begin{smallmatrix}1 & -2\Delta x & \frac{17 \Delta x^2}{4} & \biggl(-\frac{4 e^{-2\lambda \Delta x} \sinh\left[\frac{\lambda \Delta x}{2}\right]}{\lambda^3 \Delta x^3}+\frac{4 e^{-2\lambda \Delta x}\cosh[\lambda \Delta x] \sinh\left[\frac{\lambda \Delta x}{2}\right]}{\lambda^3 \Delta x^3 } \biggr) & \biggl( \frac{e^{\frac{\lambda \Delta x}{2}}(-1+e^{\lambda \Delta x})}{\lambda^3 \Delta x^3}+\frac{e^{\frac{5\lambda \Delta x}{2}}(-1+e^{\lambda \Delta x})}{\lambda^3 \Delta x^3 }-\frac{4 e^{2\lambda \Delta x}\sinh\left[\frac{\lambda \Delta x}{2}\right]}{ \lambda^3 \Delta x^3}\biggr) & \biggl(\frac{8 \cos[2\lambda \Delta x]\sin{\left[\frac{\lambda \Delta x}{2}\right]}^3}{\lambda^3 \Delta x^3}\biggr) & \biggl(-\frac{8 \sin[2\lambda \Delta x]\sin{\left[\frac{\lambda \Delta x}{2}\right]}^3}{\lambda^3 \Delta x^3}\biggr)\\
1 & -\Delta x & \frac{5 \Delta x^2}{4} & \biggl(-\frac{4 e^{-\lambda \Delta x} \sinh\left[\frac{\lambda \Delta x}{2}\right]}{\lambda^3 \Delta x^3 }+\frac{4 e^{-\lambda \Delta x}\cosh[\lambda \Delta x] \sinh\left[\frac{\lambda \Delta x}{2}\right]}{\lambda^3 \Delta x^3 } \biggr) & \biggl( \frac{e^{-\frac{\lambda \Delta x}{2}}(-1+e^{\lambda \Delta x})}{\lambda^3 \Delta x^3 }+\frac{e^{\frac{3\lambda \Delta x}{2}}(-1+e^{\lambda \Delta x})}{\lambda^3 \Delta x^3}-\frac{4 e^{\lambda \Delta x}\sinh\left[\frac{\lambda \Delta x}{2}\right]}{\lambda^3 \Delta x^3}\biggr) & \biggl(\frac{8 \cos[\lambda \Delta x]\sin{\left[\frac{\lambda \Delta x}{2}\right]}^3}{\lambda^3 \Delta x^3}\biggr) & \biggl(-\frac{8 \sin[\lambda \Delta x]\sin{\left[\frac{\lambda \Delta x}{2}\right]}^3}{\lambda^3 \Delta x^3}\biggr)\\
1 & 0 & \frac{\Delta x^2}{4} & \biggl(-\frac{4 \sinh\left[\frac{\lambda \Delta x}{2}\right]}{\lambda^3 \Delta x^3 }+\frac{4\cosh[\lambda \Delta x] \sinh\left[\frac{\lambda \Delta x}{2}\right]}{\lambda^3 \Delta x^3 } \biggr) & \biggl( \frac{e^{\frac{-3\lambda \Delta x}{2}}(-1+e^{\lambda \Delta x})}{\lambda^3 \Delta x^3}+\frac{e^{\frac{\lambda \Delta x}{2}}(-1+e^{\lambda \Delta x})}{\lambda^3 \Delta x^3 }-\frac{4\sinh\left[\frac{\lambda \Delta x}{2}\right]}{\lambda^3 \Delta x^3 }\biggr) & \biggl(\frac{8 \sin{\left[\frac{\lambda \Delta x}{2}\right]}^3}{\lambda^3 \Delta x^3}\biggr) & 0\\
1 & \Delta x & \frac{5 \Delta x^2}{4} & \biggl(-\frac{4 e^{\lambda \Delta x} \sinh\left[\frac{\lambda \Delta x}{2}\right]}{\lambda^3 \Delta x^3 }+\frac{4 e^{\lambda \Delta x}\cosh[\lambda \Delta x] \sinh\left[\frac{\lambda \Delta x}{2}\right]}{\lambda^3 \Delta x^3 } \biggr) & \biggl( \frac{e^{\frac{-5\lambda \Delta x}{2}}(-1+e^{\lambda \Delta x})}{\lambda^3 \Delta x^3 }+\frac{e^{\frac{-\lambda \Delta x}{2}}(-1+e^{\lambda \Delta x})}{\lambda^3 \Delta x^3 }-\frac{4 e^{-\lambda \Delta x}\sinh\left[\frac{\lambda \Delta x}{2}\right]}{ \lambda^3 \Delta x^3}\biggr) & \biggl(\frac{8 \cos[\lambda \Delta x]\sin{\left[\frac{\lambda \Delta x}{2}\right]}^3}{\lambda^3 \Delta x^3}\biggr) & \biggl(\frac{8 \sin[\lambda \Delta x]\sin{\left[\frac{\lambda \Delta x}{2}\right]}^3}{\lambda^3 \Delta x^3}\biggr)\\
1 & 2\Delta x & \frac{17 \Delta x^2}{4} & \biggl(-\frac{4 e^{2\lambda \Delta x} \sinh\left[\frac{2\lambda \Delta x}{2}\right]}{\lambda^3 \Delta x^3 }+\frac{4 e^{2\lambda \Delta x}\cosh[2\lambda \Delta x] \sinh\left[\frac{\lambda \Delta x}{2}\right]}{ \lambda^3 \Delta x^3} \biggr) & \biggl( \frac{e^{\frac{-7\lambda \Delta x}{2}}(-1+e^{\lambda \Delta x})}{\lambda^3 \Delta x^3 }+\frac{e^{\frac{-3\lambda \Delta x}{2}}(-1+e^{\lambda \Delta x})}{\lambda^3 \Delta x^3 }-\frac{4 e^{-2\lambda \Delta x}\sinh\left[\frac{\lambda \Delta x}{2}\right]}{\lambda^3 \Delta x^3 }\biggr) & \biggl(\frac{8 \cos[2\lambda \Delta x]\sin{\left[\frac{\lambda \Delta x}{2}\right]}^3}{\lambda^3 \Delta x^3}\biggr) & \biggl(\frac{8 \sin[2\lambda \Delta x]\sin{\left[\frac{\lambda \Delta x}{2}\right]}^3}{\lambda^3 \Delta x^3}\biggr)\\
1 & 3\Delta x & \frac{37 \Delta x^2}{4} & \biggl(-\frac{4 e^{3\lambda \Delta x} \sinh\left[\frac{\lambda \Delta x}{2}\right]}{\lambda^3 \Delta x^3 }+\frac{4 e^{3\lambda \Delta x}\cosh[\lambda \Delta x] \sinh\left[\frac{\lambda \Delta x}{2}\right]}{\lambda^3 \Delta x^3 } \biggr) & \biggl( \frac{e^{\frac{-9\lambda \Delta x}{2}}(-1+e^{\lambda \Delta x})}{\lambda^3 \Delta x^3 }+\frac{e^{\frac{-5\lambda \Delta x}{2}}(-1+e^{\lambda \Delta x})}{\lambda^3 \Delta x^3 }-\frac{4 e^{-3\lambda \Delta x}\sinh\left[\frac{\lambda \Delta x}{2}\right]}{\lambda^3 \Delta x^3 }\biggr) & \biggl(\frac{8 \cos[3\lambda \Delta x]\sin{\left[\frac{\lambda \Delta x}{2}\right]}^3}{\lambda^3 \Delta x^3}\biggr) & \biggl(\frac{8 \sin[3\lambda \Delta x]\sin{\left[\frac{\lambda \Delta x}{2}\right]}^3}{\lambda^3 \Delta x^3}\biggr)\\
1 & 4\Delta x & \frac{65 \Delta x^2}{4} & \biggl(-\frac{4 e^{4\lambda \Delta x} \sinh\left[\frac{\lambda \Delta x}{2}\right]}{ \lambda^3 \Delta x^3}+\frac{4 e^{4\lambda \Delta x}\cosh[\lambda \Delta x] \sinh\left[\frac{\lambda \Delta x}{2}\right]}{\lambda^3 \Delta x^3 } \biggr) & \biggl( \frac{e^{\frac{-11\lambda \Delta x}{2}}(-1+e^{\lambda \Delta x})}{\lambda^3 \Delta x^3 }+\frac{e^{\frac{-7\lambda \Delta x}{2}}(-1+e^{\lambda \Delta x})}{ \lambda^3 \Delta x^3}-\frac{4 e^{-4\lambda \Delta x}\sinh\left[\frac{\lambda \Delta x}{2}\right]}{\lambda^3 \Delta x^3 }\biggr) & \biggl(\frac{8 \cos[4\lambda \Delta x]\sin{\left[\frac{\lambda \Delta x}{2}\right]}^3}{\lambda^3 \Delta x^3}\biggr) & \biggl(\frac{8 \sin[4\lambda \Delta x]\sin{\left[\frac{\lambda \Delta x}{2}\right]}^3}{\lambda^3 \Delta x^3}\biggr)\end{smallmatrix}\right],
\end{eqnarray*}
the vectors $B$ and $X$ are 
\begin{equation*}
B= \left[ g^+(u_{i-2}),  g^+(u_{i-1}),  g^+(u_{i}),  g^+(u_{i+1}),  g^+(u_{i+2}),  g^+(u_{i+3}),  g^+(u_{i+4}) \right]^T,
\end{equation*} 
\begin{equation*}
X= \left[ a_0, a_1,  a_2,  a_3,  a_4,  a_5,  a_6 \right]^T.
\end{equation*}
\end{landscape}
\clearpage
\bibliographystyle{plain}

\end{document}